\documentclass[a4paper,10pt]{article}
\usepackage{amssymb,amsmath,amsthm}
\usepackage{graphicx}

\usepackage{amsfonts}
\usepackage{latexsym}
\usepackage[english]{babel}
\usepackage{hyperref}

\numberwithin{equation}{section}

\newcommand{\Ker}{\mathrm{Ker}}

\newcommand{\la}{\langle}
\newcommand{\ra}{\rangle}

\newenvironment{pf}{\noindent{\sc Proof}.\enspace}{\rule{2mm}{2mm}\smallskip}

\newtheorem{theorem}{Theorem}[section]

\newtheorem{lemma}{Lemma}[section]
\newtheorem{corollary}{Corollary}[section]

\newtheorem{remark}{Remark}[section]
\newtheorem{remarks}{Remark}[section]
\newtheorem{definition}{Definition}[section]
\newcommand{\be}{\begin{equation}}
\newcommand{\ee}{\end{equation}}

\newcommand{\om}{\omega}

\newcommand{\e}{\varepsilon}

\newcommand{\ph}{\varphi}
\renewcommand{\th}{\vartheta}

\newcommand{\wtilde}{\widetilde}
\newcommand{\pp}{\mathbb{P}}
\newcommand{\R}{\mathbb R}
\newcommand{\C}{\mathbb C}
\newcommand{\mG}{\mathcal{G}}
\newcommand{\mP}{\mathcal{P}}
\newcommand{\Lm}{\Lambda}

\newcommand{\Z}{\mathbb Z}

\newcommand{\N}{\mathbb N}
\newcommand{\T}{\mathbb T}

\renewcommand{\a }{\alpha }
\renewcommand{\b }{\beta }
\newcommand{\s }{\sigma }
\newcommand{\ii }{{\rm i} }
\renewcommand{\d }{\delta }
\newcommand{\D }{\Delta}
\newcommand{\g }{\gamma}

\renewcommand{\l }{\lambda }
\renewcommand{\L }{\Lambda }
\newcommand{\m }{\mu }
\newcommand{\vphi}{\varphi }

\renewcommand{\t }{\tau }
\renewcommand{\o }{\omega }

\newcommand{\pa}{\partial}
\newcommand{\lm}{\lambda}
\newcommand{\inv}{^{-1}}
\newcommand{\p}{\pi}

\newcommand{\Lipg}{{\rm{Lip}(\g)}}
\newcommand{\lip}{{\rm lip}}

\newcommand{\mS}{\mathcal{S}}
\newcommand{\mL}{\mathcal{L}}
\newcommand{\mN}{\mathcal{N}}

\newcommand{\mD}{\mathcal{D}}
\newcommand{\mR}{\mathcal{R}}
\newcommand{\mA}{\mathcal{A}}

\newcommand{\ompaph}{\om \cdot \pa_{\ph}}
\newcommand{\ompath}{\om \cdot \pa_{\th}}

\newcommand{\muff}{m} 

\newcommand{\mathfraks}{\mathfrak{s}_0}

\begin{document}

\title{{\bf KAM for quasi-linear and fully nonlinear forced KdV}}

\date{}

\author{Pietro Baldi, Massimiliano Berti, Riccardo Montalto}

\maketitle

\noindent
{\bf Abstract:}
We prove the existence of quasi-periodic, small amplitude, solutions for 
quasi-linear and fully nonlinear forced perturbations of  KdV equations.
For Hamiltonian or reversible nonlinearities we also obtain the linear stability of the solutions. 
The proofs are based on a combination of different ideas and techniques: 
$(i)$ a Nash-Moser iterative scheme in Sobolev scales. 
$(ii)$ A regularization procedure, which conjugates the linearized operator to a differential operator with constant coefficients plus a bounded remainder. These transformations 
are obtained by changes of variables induced by diffeomorphisms of the torus and pseudo-differential operators.
$(iii)$ A 
reducibility KAM scheme, which completes the reduction to constant coefficients of the linearized operator, 
providing a sharp asymptotic expansion of the perturbed eigenvalues. 

\smallskip


\tableofcontents

\section{Introduction}

One of the most challenging and 
open questions in KAM theory concerns 
its possible extension to \emph{quasi-linear} and \emph{fully nonlinear} PDEs, namely 
partial differential equations whose nonlinearities contain derivatives of the 
same 
order as the linear operator. 
Besides its 
mathematical interest, this question is also relevant in view of applications to physical real world 
nonlinear models,  for example in fluid dynamics and elasticity.

\smallskip

The goal of this paper is to develop KAM theory for quasi-periodically forced KdV equations of the form 
\be\label{equation main} 
u_{t} + u_{xxx} + \e f(\o t , x , u, u_{x}, u_{xx}, u_{xxx} ) = 0 \, , \quad 
x \in \T := \R / 2\pi\Z \, . 
\ee 
First, we prove in Theorem \ref{thm:main}  an existence result 
of quasi-periodic solutions for a large class of quasi-linear nonlinearities $ f $.
Then for Hamiltonian or reversible nonlinearities, we also 
prove the linear stability of the solutions, see Theorems \ref{thm:mainH}, \ref{thm:mainrev}.  
Theorem \ref{thm:mainrev}  also holds for fully nonlinear perturbations. 
The precise meaning of stability is stated in Theorem \ref{cor:stab}. 
The key analysis is the reduction to constant coefficients of the 
linearized KdV equation, see Theorem \ref{thm:reducibility}. 
To the best of our knowledge, these are the first KAM results for quasi-linear or fully nonlinear PDEs.

\smallskip

Let us outline a short history of the subject. 
KAM and Nash-Moser theory for PDEs, which counts nowadays on a wide literature, 
started with the pioneering works of Kuksin \cite{Ku} and Wayne \cite{W1},  
and was developed in the 1990s by 
Craig-Wayne \cite{CW},  Bourgain \cite{Bo1}, \cite{B3}, P\"oschel \cite{Po2} 
(see also \cite{k1}, \cite{C} for more references).  
These papers concern wave and Schr\"odinger equations 
with bounded Hamiltonian nonlinearities.

\smallskip

The first KAM results for \emph{unbounded} perturbations 
have been obtained by Kuksin \cite{K2}, \cite{k1}, and, then, 
Kappeler-P\"oschel \cite{KaP}, 
for  Hamiltonian, analytic perturbations of KdV.   
Here the highest constant coefficients linear operator is  $\pa_{xxx}$ and the nonlinearity contains one space derivative $\pa_x$. 
Their approach has been recently improved by Liu-Yuan \cite{LY} and Zhang-Gao-Yuan \cite{ZGY} 
for $1$-dimensional derivative NLS (DNLS) and Benjamin-Ono equations, 
where the highest order constant coefficients linear operator is $ \pa_{xx}$ and the nonlinearity contains one derivative $\pa_x$.
These methods apply to dispersive PDEs with derivatives like KdV, DNLS, 
the  Duffing oscillator (see Bambusi-Graffi  \cite{Bambusi-Graffi}), but 
not to derivative wave equations (DNLW) which contain first order derivatives $\pa_x ,  \partial_t $ in the nonlinearity.

For DNLW, KAM theorems have been recently proved by Berti-Biasco-Procesi 
for both Hamiltonian \cite{BBiP1} and reversible  \cite{BBiP2} equations. 
The key ingredient is an asymptotic expansion of the perturbed eigenvalues 
that is sufficiently accurate to impose the second order Melnikov non-resonance conditions. 
In this way, the  scheme produces a constant coefficients normal form around the invariant torus (\emph{reducibility}), implying the linear stability of the solution.
This is achieved introducing the notion  of ``quasi-T\"oplitz''  vector field,
which is inspired to the concept of ``quasi-T\"oplitz" and ``T\"oplitz-Lipschitz''  Hamiltonians, 
developed, respectively,  in Procesi-Xu \cite{PX} and Eliasson-Kuksin \cite{EK}, \cite{EK1} 
(see also Geng-You-Xu \cite{GXY}, Gr\'ebert-Thomann \cite{GT}, Procesi-Procesi \cite{PP}).

\smallskip

Existence of quasi-periodic solutions of PDEs can also be proved by imposing only the first order Melnikov conditions. 
This approach  has been developed by Bourgain \cite{Bo1}-\cite{B5} 
extending the work of Craig-Wayne \cite{CW} for periodic solutions. 
It is especially convenient for PDEs in higher space dimension, because of the high multiplicity of the eigenvalues: see also the recent results by 
Wang \cite{Wang},  
Berti-Bolle \cite{BBo10}, \cite{BB12}
(and \cite{Berti-book}, \cite{Berti-Bolle-Procesi-AIHP-2010}, \cite{GP} for periodic solutions). 
This method does not provide informations about the stability of the quasi-periodic solutions, because
the linearized equations have variable coefficients.

\smallskip

All the aforementioned results concern ``semilinear'' PDEs,  
namely equations in which 
the nonlinearity contains \emph{strictly less} derivatives than the linear differential operator. 
For quasi-linear or fully nonlinear PDEs the perturbative effect is much stronger, 
and the possibility of extending KAM theory in this context is doubtful, 
see \cite{KaP},  \cite{C}, \cite{LY},  
because of the possible phenomenon of formation of singularities outlined in
Lax \cite{Lax}, Klainerman and Majda \cite{KM}. 
For example, Kappeler-P\"oschel \cite{KaP} (remark 3, page 19) wrote: 
``{\it It would be interesting to obtain perturbation results which also include terms of higher order, at least in the region where the KdV approximation is valid. However, results of this type are still out of reach, if true at all}''.
The study of this important issue is at its infancy.  

\smallskip

For quasi-linear and fully nonlinear PDEs, the  literature concerns, so far, 
only existence of \emph{periodic} solutions. 
We quote
the classical bifurcation results  
of Rabinowitz \cite{Rabinowitz-tesi-1967} for fully nonlinear forced wave equations
with a small dissipation term. 
More recently,  Baldi \cite{Baldi Kirchhoff} proved existence of periodic forced vibrations 
for  quasi-linear Kirchhoff equations.  
Here the quasi-linear perturbation term 
depends explicitly only on time. Both these results are proved via Nash-Moser methods.

For the water waves equations, which are  a fully nonlinear PDE, we mention 
the pioneering work of Iooss-Plotnikov-Toland \cite{Ioo-Plo-Tol}   
about the existence of time periodic standing waves, 
and of Iooss-Plotinikov 
\cite{IP09}, \cite{IP11} for 3-dimensional traveling water waves. 
The key idea is to 
use diffeomorphisms of the torus $\T^2$ and pseudo-differential operators, in order to conjugate
the linearized operator (at an approximate solution)
to a constant coefficients operator
plus a sufficiently regularizing remainder. 
This is enough to invert the whole linearized operator by Neumann series.  

Very recently Baldi \cite{Baldi-Benj-Ono} has further developed the techniques of \cite{Ioo-Plo-Tol},
proving the existence of periodic solutions for fully nonlinear autonomous, reversible Benjamin-Ono equations. 

These approaches  do not imply the linear stability of the solutions and, unfortunately, 
they do not work for quasi-periodic solutions, because stronger small divisors difficulties arise, 
see the observation \ref{obs4} below. 

\smallskip

We finally mention that,  for quasi-linear 
Klein-Gordon equations on  
spheres, Delort \cite{Delort-2009} has proved  long time existence results  
via Birkhoff normal form methods.  

\smallskip

In the present paper we combine different ideas and techniques.  
The key analysis concerns  the linearized KdV operator \eqref{linearized op} obtained at any step of the Nash-Moser iteration.
First, we use changes of variables, like quasi-periodic time-dependent 
diffeomorphisms of the space variable $ x $, 
a quasi-periodic reparametrization of time, 
multiplication operators and Fourier multipliers, 
in order to reduce the linearized operator to constant coefficients 
up to a bounded remainder, see \eqref{L6red}.  
These transformations, which are inspired to \cite{Baldi-Benj-Ono}, \cite{Ioo-Plo-Tol}, 
are very different from the usual KAM transformations. 
Then, we perform a quadratic KAM reducibility scheme {\it \`a la} Eliasson-Kuksin, 
which completely diagonalizes the linearized operator. 
For reversible or Hamiltonian KdV perturbations 
we get that the eigenvalues of this diagonal operator are purely imaginary, i.e. we prove the linear stability. 
In section \ref{sec:ideas} we present the main ideas of proof. 

We remark that the present  approach could be also applied to quasi-linear and fully nonlinear perturbations of dispersive PDEs like 1-dimensional NLS 
and Benjamin-Ono equations (but not to the wave equation, which is not dispersive).
For definiteness, we have developed all the  computations in KdV case. 

\smallskip

In the next subsection we state precisely our KAM results.
In order to highlight the main ideas,  
we consider the simplest setting of nonlinear perturbations of the Airy-KdV operator 
$ \pa_t + \pa_{xxx}$ and we look for small amplitude solutions.

\subsection{Main results} 

We consider 
problem \eqref{equation main} 
where $ \e > 0 $ is a small parameter, 
the nonlinearity is quasi-periodic in time with diophantine frequency vector
\be\label{omdio}
\om = \l \bar \om \in \R^{\nu} \, , \quad \l \in \L := \Big[ \frac12\,, \frac32 \Big], \quad 
|\bar \om \cdot l | \geq \frac{3 \g_0}{|l|^{\t_0}} \quad \forall l \in \Z^{\nu} \setminus \{ 0 \}, 
\ee
and $ f(\vphi, x, z )$,  $\ph \in \T^{\nu}$,  $ z := (z_0, z_1, z_2, z_3) \in \R^4 $, 
is a finitely many times differentiable function, namely
\be\label{f classe Cq}
f \in C^q ( \T^{\nu} \times \T \times \R^4; \R) 
\ee
for some  $ q \in \N $ large enough. For simplicity we fix in  \eqref{omdio}  the diophantine exponent $ \t_0 := \nu $. 
The only ``external'' parameter in \eqref{equation main} is $ \l $, 
which is the length of the frequency vector (this corresponds to a time scaling).

\smallskip

We consider the following questions:

\begin{itemize}
\item {\it For $ \e $ small enough, 
do there exist quasi-periodic solutions of \eqref{equation main} for positive measure sets of 
$ \l \in \L $? } 

\vspace{-4pt} 

\item {\it Are these solutions linearly stable?} 
\end{itemize}
Clearly, if $ f(\vphi,x, 0)$ is not identically zero, then $ u = 0 $ is not a solution of \eqref{equation main} for $ \e \neq  0 $.
Thus we look for non-trivial $ (2 \pi)^{\nu+1}$-periodic solutions $ u(\vphi,x) $ of 
\be\label{eq:invqp}
\o \cdot \partial_{\vphi} u + u_{xxx} 
+ \e f(\vphi , x , u, u_{x}, u_{xx}, u_{xxx} ) = 0 
\ee
in the Sobolev space
\begin{align} 
H^s & := H^s  ( \T^\nu \times \T; \R ) \label{Hs1}
\\ & := 
\Big\{ u(\vphi,x) = \sum_{(l,j) \in \Z^{\nu} \times \Z} 
u_{l,j} \, e^{\ii (l \cdot \vphi + jx)} \in \R,  
\ \
{\bar u}_{l,j} =  u_{-l,-j} \,, 
\ \  
\| u \|_s^2 := \sum_{(l,j) \in \Z^{\nu} \times \Z} \langle l, j \rangle^{2s} 
| u_{l,j} |^ 2 <  \infty   \Big\}  \nonumber 
\end{align}
where 
\[
\la l,j \ra := \max \{ 1, |l|, |j| \}. 
\]
From now on, we fix $ {\mathfrak s}_0 := (\nu + 2) / 2  > (\nu +1 ) / 2 $, so that for all $s \geq \mathfraks$ the Sobolev space $H^s$ is a Banach algebra, and it is continuously embedded $ H^s (\T^{\nu+1} ) \hookrightarrow  C(\T^{\nu+1} ) $. 
 
\smallskip

We need some assumptions on the nonlinearity. 
We consider {\it fully nonlinear} perturbations satisfying
\begin{itemize}
\item
{\sc Type (F)} 
\begin{equation} \label{type F}
\pa_{z_2} f = 0, 
\end{equation}
\end{itemize}
namely $ f $ is independent of $ u_{xx} $. 
Otherwise, we require that 
\begin{itemize}
\item
{\sc Type (Q)}  
\begin{equation}  \label{type Q}
\pa^2_{z_3 z_3} f = 0, \quad 
\pa_{z_2} f = \a(\ph) \Big( \pa^2_{z_3 x} f + z_1 \pa^2_{z_3 z_0} f + z_2 \pa^2_{z_3 z_1} f + z_3 \pa^2_{z_3 z_2} f \Big)  
\end{equation}
for some function $ \a(\ph) $ (independent on $ x $). 
\end{itemize}
If (Q) holds, then the nonlinearity $  f $ depends linearly on  $u_{xxx} $, 
namely equation \eqref{equation main} is {\it quasi-linear}.  
We note that the Hamiltonian nonlinearities, see \eqref{f Ham}, are a particular case of those 
satisfying  (Q),
see remark  \ref{rem:coeff in cases Q F}.
In comment  \ref{linea} after Theorem \ref{cor:stab} we explain the reason for assuming either condition (F) or (Q). 

\smallskip

The following theorem is an existence result of quasi-periodic solutions for quasi-linear KdV equations.

\begin{theorem}  \label{thm:main} {\bf (Existence)}
There exist $ s := s( \nu ) > 0$,  $ q := q( \nu) \in \N $, such that: 
\\[1mm]
For every quasi-linear nonlinearity $ f  \in C^q $ of the form 
\be\label{f = der g}
f = \partial_x \big( g(\om t, x, u, u_x, u_{xx})\big) 
\ee 
satisfying the (Q)-condition \eqref{type Q}, 
for all $\e \in (0, \e_0)$, where $\e_0 := \e_0 (f, \nu) $ is small enough, 
there exists  a Cantor set $ {\cal C}_\e \subset \L $ of asymptotically full Lebesgue measure, i.e.   
\be\label{Cmeas}
| {\cal C}_\e |  \to 1 \quad \text{as} \quad \e \to 0,  
\ee
such that, $ \forall \l \in  {\cal C}_\e  $ the perturbed KdV equation \eqref{eq:invqp} 
has a solution $  u( \e, \l) \in H^s $ with $ \| u(\e, \l) \|_s \to 0 $ as $ \e \to 0 $.
\end{theorem}
 
We may ensure the {\it linear stability} of the solutions requiring further conditions on the nonlinearity, see
Theorem \ref{cor:stab} for the precise statement.
The first case is that of {\it Hamiltonian} KdV equations
\be\label{Ham-KdV}
u_t = \partial_x \nabla_{L^2} H(t,x,u, u_x)  \, , \quad 
H(t,x,u,u_x) := \int_{\T} \frac{u_x^2}{2}\, + \e F(\om t,x,u,u_x) \, dx
\ee
which have the form \eqref{equation main}, \eqref{f = der g} with
\be  \label{f Ham}
f(\ph,x,u,u_x, u_{xx}, u_{xxx}) 
 =  - \partial_x \big\{ (\pa_{z_0} F)(\ph, x,  u , u_x) \big\} 
+ \pa_{xx} \big\{ (\pa_{z_1} F)(\ph, x,  u,  u_x ) \big\} \, .
\ee
The phase space of \eqref{Ham-KdV} is 
$$ 
H^1_0 (\T) := \Big\{  u(x) \in H^1(\T, \R) \, : \, \int_{\T} u(x) \, dx = 0  \Big\}
$$
endowed with the non-degenerate symplectic form
\be\label{KdV symplectic}
\Omega (u,v) := \int_{\T} (\partial_x^{-1} u)  v \, dx  \, , \quad \forall u, v \in H_0^1 (\T) \, , 
\ee
where 
$ \partial_x^{-1} u  $ is the periodic primitive of $ u $ 
with zero average, see \eqref{dx-1}.  
As proved in remark \ref{rem:coeff in cases Q F}, the Hamiltonian nonlinearity $  f $ in \eqref{f Ham} 
satisfies also the (Q)-condition \eqref{type Q}.
As a consequence,  Theorem \ref{thm:main} implies the existence of 
quasi-periodic solutions of  \eqref{Ham-KdV}. In addition, we  
also prove their linear stability. 

\begin{theorem}  \label{thm:mainH} {\bf (Hamiltonian KdV)}
For all Hamiltonian quasi-linear KdV equations \eqref{Ham-KdV} the quasi-periodic solution $u(\e,\lm)$
found in Theorem \ref{thm:main}  is {\sc linearly stable} (see Theorem \ref{cor:stab}). 
\end{theorem}

The  stability of the quasi-periodic solutions also follows by 
the {\it reversibility} condition
\begin{equation}  \label{parity f}
f (-\ph, -x, z_0, -z_1, z_2, -z_3) = - f(\ph, x, z_0, z_1, z_2, z_3). 
\end{equation}
Actually \eqref{parity f} implies that the infinite-dimensional non-autonomous dynamical system
$$
u_t = V(t, u ), \quad 
V(t,u ) := - u_{xxx} - \e f(\o t , x , u, u_{x}, u_{xx}, u_{xxx}) 
$$
is reversible with respect to the involution
\[
S : u(x) \rightarrow u(-x), \quad 
S^2 = I, 
\]
namely
\[
- S V(-t,u) = V(t,Su) \, .
\]
In this case it is natural to look for ``reversible" solutions of \eqref{eq:invqp},  
that is
\be\label{solPP}
u ( \vphi,x ) = u ( -\vphi , -x ) \, . 
\ee

\begin{theorem}  \label{thm:mainrev}
{\bf (Reversible KdV)} 
There exist $ s := s( \nu ) > 0$,  $ q := q( \nu) \in \N $, such that: 

\noindent   
For every nonlinearity $ f  \in C^q $ that satisfies 
\\[1mm]
\indent
$(i)$ the reversibility condition \eqref{parity f}, 
\\[1mm]
\noindent
and  
\\[1mm]
\indent
$(ii)$ either the (F)-condition \eqref{type F} or the (Q)-condition \eqref{type Q}, 
\\[1mm]
for all $\e \in (0, \e_0)$, where $\e_0 := \e_0 (f, \nu) $ is small enough, 
there exists  a Cantor set $ {\cal C}_\e \subset \L $ with Lebesgue measure satisfying \eqref{Cmeas}, 
such that for all $ \l \in  {\cal C}_\e $ the perturbed KdV equation \eqref{eq:invqp} 
has a solution $ u (\e, \l) \in H^s $ that satisfies \eqref{solPP}, with 
$ \| u (\e, \l) \|_s \to 0 $ as $ \e \to 0 $.
In addition, $u(\e,\lm)$ is {\sc linearly stable}. 
\end{theorem}

Let us make some  comments on the results. 

\begin{enumerate}
\item 
The previous theorems (in particular the Hamiltonian Theorem  \ref{thm:mainH})  
give a positive answer to the question that was posed by Kappeler-P\"oschel \cite{KaP}, page 19, Remark 3, about the possibility of KAM type results for quasi-linear perturbations of KdV.

\item 
In Theorem \ref{thm:main}  we do not have informations 
about the linear stability of the solutions because the nonlinearity $ f $ has no
special structure and it may happen that  some eigenvalues of the linearized operator have non zero real part
(partially hyperbolic tori). We remark that, in any case, we may compute
the eigenvalues (i.e. Lyapunov exponents) of the linearized operator  
with any order of accuracy. 
With further conditions on the nonlinearity---like reversibility or in the Hamiltonian case---the 
eigenvalues  are purely imaginary, 
and the torus is linearly stable.
The present situation   
is very different with respect to  
\cite{CW},  \cite{Bo1}-\cite{B5}, \cite{BBo10}-\cite{BB12} and also \cite{Ioo-Plo-Tol}-\cite{IP11}, \cite{Baldi-Benj-Ono},
 where the  lack of stability informations  is due to the fact that the linearized equation   
has variable coefficients, and it is not reduced as in  Theorem \ref{thm:reducibility} below.
\item 
One cannot expect the existence of quasi-periodic solutions of \eqref{eq:invqp} for {\it any} perturbation $ f $. 
Actually, if $ f = m \neq 0 $ is a constant, then, integrating \eqref{eq:invqp} in $ (\vphi,x) $ we find the contradiction $ \e m = 0 $.
This is a consequence of the fact that 
\be\label{Kernel}
\Ker (\om \cdot \partial_\vphi + \partial_{xxx}) = \R  
\ee
is non trivial.  
Both the condition \eqref{f = der g}  (which is satisfied by the Hamiltonian nonlinearities) 
and the reversibility condition \eqref{parity f} 
allow to overcome this obstruction, working in a space of functions with zero average.  
The degeneracy \eqref{Kernel} also reflects in the fact that 
the solutions of \eqref{eq:invqp} appear as a $1$-dimensional family 
$ c + u_c( \e, \l) $ parametrized by the ``average'' $ c \in \R $. 
We could also avoid this degeneracy by adding a ``mass'' term $ + m u $ in \eqref{equation main}, 
but it does not seem to have physical meaning.  

\item
In Theorem \ref{thm:main} we have not considered the case in which $ f $ is fully nonlinear and satisfies  
condition (F) in \eqref{type F}, because any nonlinearity of the form \eqref{f = der g} is automatically quasi-linear
(and so the first condition in \eqref{type Q} holds) and 
\eqref{type F}  
trivially implies the second condition in \eqref{type Q}
with $ \a (\vphi ) = 0 $.

\item 
The solutions $ u \in H^s $ have the same regularity in both variables $ (\vphi,x) $. 
This functional setting is convenient when using changes of variables that mix the time and space variables, like
the composition operators $\mA$, $\mathcal{T}$ in sections \ref{step-1}, \ref{step-4},

\item
In the Hamiltonian case \eqref{Ham-KdV}, 
the nonlinearity $f$ in \eqref{f Ham} satisfies the reversibility condition \eqref{parity f} 
if and only if 
$ F( -\ph, -x, z_0, -z_1) =  F( \ph, x, z_0, z_1) $. 
\end{enumerate}

Theorems \ref{thm:main}-\ref{thm:mainrev} are based on a Nash-Moser iterative scheme. 
An essential ingredient in the proof---which also implies the linear stability of the quasi-periodic solutions---is the {\it reducibility} of the linear operator 
\be\label{linearized op}
\mL := \mL (u) = \om \cdot \partial_\ph + (1 + a_3(\vphi,x)) \pa_{xxx} + a_2(\vphi,x) \pa_{xx} + a_1(\vphi,x) \pa_x + a_0 (\vphi,x) 
\ee
obtained linearizing \eqref{eq:invqp} at any approximate (or exact) solution $ u $, namely
the coefficients $ a_i (\vphi, x) $ are defined in \eqref{ai formula}.
Let $ H^s_x := H^s (\T) $ denote the usual Sobolev spaces of functions of $ x \in \T $ only (phase space).

\begin{theorem}\label{thm:reducibility} {{\bf (Reducibility)}} 
There exist $ \bar \s >  0 $, $ q \in \N $, depending on $ \nu $, such that:
\\[1mm]  
For every nonlinearity $ f  \in C^q $ that satisfies the hypotheses of Theorems \ref{thm:main} or \ref{thm:mainrev},   
for all $\e \in (0, \e_0)$, where $\e_0 := \e_0 (f, \nu) $ is small enough, 
for all $u$ in the ball $\| u \|_{ { \mathfrak s}_0  + \bar \s} \leq 1$, 
there exists a Cantor like set $  \L_\infty (u)  \subset \L $ 
such that,  for all $ \l \in \L_\infty (u) $: 
\\[1mm]
{i)} for all $ s \in ({\mathfrak s}_0, q - \bar \s) $, if  $\| u \|_{ s  + \bar \s} < + \infty $ then  
there exist linear invertible bounded operators $ W_1 $, $ W_2  : H^s (\T^{\nu+1})\to H^s ( \T^{\nu+1} ) $  
with bounded inverse, that semi-conjugate
the linear operator $ {\cal L}(u) $ in \eqref{linearized op} to the diagonal operator
$ {\cal L}_\infty $, namely 
\be\label{semicon}
{\cal L}(u)  = W_1  {\cal L}_\infty  W_2^{-1}  \, , \quad {\cal L}_\infty := \ompaph + {\cal D}_\infty
\ee
where 
\be\label{thm:diag}
{\cal D}_\infty := {\rm diag}_{j \in \Z} \{ \mu_j \}, 
\quad \mu_j := \ii (-m_3 j^3 + m_1 j) + r_j \, , 
\quad m_3, m_1 \in \R \, ,  
\quad \sup_j |r_j | \leq C \e \, . 
\ee
{ii)} For each $ \vphi \in \T^\nu $ the operators $ W_i  $ are also bounded 
linear bijections of the phase space (see notation \eqref{notationA})
$$
W_i  ( \vphi ) \, ,  W_i^{-1}  ( \vphi ) : H^s_x \to H^s_x \, , \quad i = 1,2 \, .
$$
A curve $ h(t) = h(t, \cdot ) \in H^{s}_x $ is a solution of  the quasi-periodically forced linear KdV equation 
\be\label{KdV:lin}
\pa_t h +  (1 + a_3(\om t,x)) \pa_{xxx}h + a_2(\om t,x) \pa_{xx}h + a_1(\om t,x) \pa_xh + a_0 (\om t,x)h = 0 
\ee
if and only if the transformed curve 
$$ 
v(t) := v(t, \cdot ) := W_2^{-1} ( \om t ) [h(t)] \in H^{s}_x 
$$ 
is a solution of the constant coefficients dynamical system
\be\label{Lin: Red}
\pa_t v + {\cal D}_\infty v = 0 \,  , \quad  
{\dot v}_j = - \mu_j v_j \, , \ \ \forall j \in \Z \, . 
\ee
In the reversible or Hamiltonian case  all the $ \mu_j \in \ii \R $ are purely imaginary. 
\end{theorem}

The exponents $ \mu_j $ can be effectively computed. 
All the solutions of \eqref{Lin: Red} are
$$
v(t) = \sum_{j \in \Z}  v_j(t)  e^{\ii j x} \, , \quad   v_j(t) = e^{- \mu_j t } v_j(0) \, . 
$$
If the $ \mu_j $ are purely imaginary -- as in the reversible or the Hamiltonian cases -- 
all the solutions of \eqref{Lin: Red} are almost periodic in time (in general) and the Sobolev norm
\be\label{constant v}
\| v(t) \|_{H^s_x} = \Big( \sum_{j \in \Z}  |v_j(t)|^2 \langle j \rangle^{2s}\Big)^{1/2} =  
\Big( \sum_{j \in \Z}  |v_j(0)|^2 \langle j \rangle^{2s}\Big)^{1/2} =
\| v(0) \|_{H^s_x} 
\ee
is constant in time.  As a consequence we have:

\begin{theorem}\label{cor:stab} {\bf (Linear stability)} Assume 
the hypothesis of Theorem \ref{thm:reducibility} and, in addition, that  $ f $ is Hamiltonian (see \eqref{f Ham}) 
or it satisfies the reversibility condition \eqref{parity f}. 
Then,   $ \forall s \in ( \mathfrak{s}_0, q - \bar \s - \mathfrak s_0) $, 
$ \| u  \|_{s+ \mathfrak s_0 + \bar \s } < + \infty $,  there exists $ K_0 > 0 $ such that 
for all $\l \in \Lambda_\infty(u) $, $\e \in (0,\e_0)$, all the solutions of \eqref{KdV:lin} satisfy 
\be \label{stability s}
\| h(t)\|_{H^s_x} \leq K_0 \| h(0)\|_{H^s_x} \,  
\ee
and, for some   $ \mathtt a  \in (0,1) $, 
\be \label{stability epsilon}
\| h(0)\|_{H^s_x} - \e^{\mathtt a}  K_0 \| h(0)\|_{H^{s+1}_x} 
\leq \| h(t)\|_{H^s_x} \leq 
\| h(0)\|_{H^s_x} + \e^{\mathtt a} K_0 \| h(0)\|_{H^{s+1}_x} \, . 
\ee
\end{theorem}

Theorems  \ref{thm:main}-\ref{cor:stab} are proved in section \ref{sec:proof} collecting all the informations of 
sections \ref{sec:2}-\ref{sec:NM}. 

\subsection{Ideas of proof}\label{sec:ideas}

The proof of Theorems  \ref{thm:main}-\ref{thm:mainrev}  is based 
on a Nash-Moser iterative scheme in the scale of Sobolev spaces $ H^s $. 
The main issue concerns the invertibility of the linearized KdV operator $ {\cal L} $  in  \eqref{linearized op}, at each step of the iteration,  and the proof of the tame estimates \eqref{L-1alta} for its right inverse.
This information  is obtained in Theorem \ref{inversione linearizzato} 
by conjugating $ {\cal L} $ to constant coefficients.
This is also the key  which implies the stability results for the Hamiltonian and reversible nonlinearities,
see Theorems \ref{thm:reducibility}-\ref{cor:stab}. 

\smallskip

We now explain the main ideas of the reducibility scheme. 
The term of $ {\cal L} $ that produces the strongest perturbative effects to the spectrum (and eigenfunctions) is
$ a_3 (\vphi,x) \partial_{xxx} $, and, then $ a_2 (\vphi,x) \partial_{xx} $.
The usual KAM transformations are not able to deal with these terms because they are ``too close"
to the identity. Our strategy is the following. 
First, we conjugate the operator $ \mL  $ in \eqref{linearized op}
to a constant coefficients third order differential operator plus a zero order remainder
\be\label{L6red}
\mL_5 = \om \cdot \partial_\ph + m_3 \partial_{xxx} + m_1 \partial_x + {\cal R}_0, 
\quad m_3 = 1 + O(\e), 
\ m_1 = O(\e ) \, , \ m_1, m_3 \in \R \, , 
\ee
(see \eqref{mL5}),  via changes of variables induced by diffeomorphisms of the torus, 
reparametrization of time, and pseudo-differential operators. 
This is the goal  of section \ref{sec:regu}. 
All these transformations could be composed into one map, but we find it
more convenient to split the regularization procedure into separate steps (sections \ref{step-1}-\ref{step-5}), 
both 
to highlight the basic  ideas, and, especially, in order to derive estimates on the coefficients, section 
\ref{subsec:mL0 mL5}.  Let us make some comments on this procedure. 
\begin{enumerate}
\item  \label{comment2}
In order to  eliminate the space variable dependence of the 
highest order perturbation $ a_3 (\vphi,x) \partial_{xxx} $ (see \eqref{mL1})
we use,  in section \ref{step-1}, $\vphi $-dependent changes of variables
 like
$$
({\cal A} h)(\vphi, x) := h(\vphi, x + \beta (\vphi, x)) \, . 
$$
These transformations converge pointwise to the identity if $ \b \to 0 $ but not in operatorial norm.
If $ \b $ is odd, $\mA$ preserves 
the reversible structure, see remark \ref{reversibilitˆ step 1}.
On the other hand for the Hamiltonian KdV \eqref{Ham-KdV} we use the modified transformation
\begin{equation}\label{operatore1 simplettico}
({\cal A}h)(\vphi,x):= (1+ \b_x(\vphi, x)) \, h(\vphi, x + \b (\vphi, x)) = 
\frac{d}{dx} \big\{ ({\partial_x}^{-1} h )(\vphi, x+ \b(\vphi,x)) \big\} 
\end{equation}
for all  $ h(\vphi, \cdot ) \in H^1_0 (\T) $.   
This map is canonical, for each $ \vphi \in \T^\nu $, with respect to
the KdV-symplectic form \eqref{KdV symplectic}, see remark  \ref{rem: Ham0}.
Thus \eqref{operatore1 simplettico} preserves the Hamiltonian structure and
also eliminates the term of order $ \partial_{xx} $, see 
remark  \ref{rem: Ham1}. 
\item In the second step of section \ref{step-2} we eliminate the time dependence  of the coefficients of the 
highest order spatial derivative operator $ \partial_{xxx} $
by a quasi-periodic time re-parametrization. 
This procedure preserves the reversible and the Hamiltonian structure, see remark 
\ref{reversibilitˆ step 2} and \ref{rem: Ham2}. 
\item \label{linea}
Assumptions (Q) (see \eqref{type Q}) or (F) (see \eqref{type F})
allow  to eliminate terms like  $ a (\vphi, x) \partial_{xx} $ along this reduction procedure, see  \eqref{mL3}. 
This is possible, by a conjugation with multiplication operators (see \eqref{cambio3}), 
if  
(see \eqref{viapez}) 
\begin{equation}  \label{zero mean in the intro}
\int_\T \frac{a_2(\ph,x)}{1 + a_3(\ph,x)} \, dx = 0 \, . 
\end{equation}
If  (F) holds, then the coefficient $ a_2(\ph,x) = 0 $ and \eqref{zero mean in the intro} is satisfied.
If (Q) holds, then an easy computation shows that  
$ a_2(\ph,x) = \a(\ph) \, \pa_x a_3(\ph,x) $ 
(using the explicit expression of the coefficients in \eqref{ai formula}), and so
\[
\int_\T  \frac{a_2(\ph,x)}{1 + a_3(\ph,x)} \, dx  =  
\int_\T \a(\ph) \, \pa_x \big( \log[ 1+a_3(\ph,x)] \big) \, dx  = 0  \, .
\]
In both cases (Q) and (F), condition  \eqref{zero mean in the intro} is satisfied. 

In the Hamiltonian case there is no need of this step because
the symplectic transformation \eqref{operatore1 simplettico} also eliminates 
the term of order $ \partial_{xx} $, see remark \ref{rem: Ham2}.

We note  that  without assumptions (Q) or (F) we may always reduce $\mL$ 
to a time dependent operator with $ a (\vphi ) \partial_{xx} $. 
If $ a(\vphi ) $ were a constant, then this term would even simplify the analysis, 
killing the small divisors. 
The pathological situation that we want to eliminate 
assuming  (Q) or (F)  is when $ a(\vphi ) $ changes sign. 
In such a case, this term acts as a friction when $ a(\vphi) < 0 $ and 
as an amplifier when $ a(\vphi) > 0 $.

\item 
In  sections \ref{step-4}-\ref{step-5}, we 
are finally able to conjugate the linear operator to another one with a coefficient in front of $ \partial_x $
which is  constant, i.e.  obtaining  \eqref{L6red}. 
In this step we use 
a transformation of the form
$  I +  w(\vphi,x) \partial_x^{-1} $, see \eqref{mS}.  
In the Hamiltonian case we use the symplectic map
$ e^{\pi_0 w(\vphi, x) \partial_x^{-1}} $, see remark \ref{rem:Ham5}. 

\item \label{obs4}
We can iterate the regularization procedure at any {\it finite} order $ k = 0, 1, \ldots $, 
conjugating 
$ {\cal L} $ to an operator of the form ${\mathfrak D} + {\cal R }$, where 
$$ 
{\mathfrak D} =  \om \cdot \partial_\ph + \mD, 
\quad
\mD = m_3 \partial_{x}^3 +  m_1 \partial_x + \ldots
+ m_{-k} \partial_x^{-k} \, ,
\quad m_{i} \in \R \, ,  
$$ 
has constant coefficients, and the rest $ {\cal R } $ is arbitrarily regularizing in space, namely
\be\label{Lmany}
\pa_x^{k} \circ \mR = \text{bounded} \, .  
\ee
However, one cannot iterate this regularization infinitely many times, because 
it is not a quadratic scheme, and therefore, 
because of the small divisors, it does not converge. 
This regularization procedure is sufficient to prove the invertibility of $ {\cal L} $,
giving tame estimates for the inverse, in the periodic case, 
but it does not work for quasi-periodic solutions. 
The reason is the following. 
In order to use Neumann series, one needs that 
${\mathfrak D} \inv \mR = ({\mathfrak D}\inv \pa_x^{-k}) (\pa_x^{k} \mR)$ is bounded, 
namely, in view of \eqref{Lmany},  that $ {\mathfrak D} \inv \pa_x^{-k} $ is bounded.  
In the region where the eigenvalues $(\ii \om \cdot l + \mD_j)$ of ${\mathfrak D} $ are small, 
space and time derivatives are related, 
$|\om \cdot l| \sim |j|^3$, 
where $l$ is the Fourier index of time, $j$ is that of space, 
and $\mD_j = - \ii m_3 j^3 + \ii m_1 j + \ldots$ are the eigenvalues of $\mD$.
Imposing the first order Melnikov conditions $|\ii \om \cdot l + \mD_j| > \g |l|^{-\t}$,   
in that region, $( {\mathfrak D} \inv \pa_x^{-k})$ has eigenvalues 
\[
\Big| \frac{1}{(\ii \om \cdot l + \mD_j) j^{k}}\, \Big| 
< \frac{|l|^\t}{\g |j|^{k}} \, < \frac{C |l|^\t}{|\om \cdot l|^{k/3}} \,.
\]
In the periodic case, $\om \in \R$, $l \in \Z$, $|\om \cdot l| = |\om| |l|$, and this determines the order of regularization that is required by the procedure: $ k \geq 3 \t$. 
In the quasi-periodic case, instead, $|l|$ is not controlled by $|\om \cdot l|$, and the argument fails.
\end{enumerate}

Once \eqref{L6red} has been obtained, we implement a quadratic reducibility KAM scheme 
to diagonalize $ {\cal L}_5 $, namely to conjugate 
$  {\cal L}_5  $  to the diagonal operator $ {\cal L}_\infty $ in \eqref{semicon}. 
Since we work with finite regularity, we perform a Nash-Moser smoothing regularization 
(time-Fourier truncation).
We use standard KAM transformations, in order to decrease, quadratically at each step, 
the size of the perturbation $\mR$, see section \ref{the-reducibility-step}. 
This iterative scheme converges (Theorem \ref{thm:abstract linear reducibility}) 
because the initial remainder $ {\cal R}_0 $ is a bounded operator (of the space variable $x$),
and this property is preserved along the iteration. 
This is the reason for performing 
the regularization procedure of sections \ref{step-1}-\ref{step-5}. 
We manage to impose the second order Melnikov non-resonance conditions 
\eqref{Omgj}, which are required by the reducibility scheme,  
thanks to 
the good control of the eigenvalues
$ \mu_j =  - \ii m_3(\e,\l) j^3 +  \ii m_1(\e,\l) j + r_j (\e,\l) $, 
where $  \sup_j |r_j (\e,\l)|  = O(\e ) $. 

Note that the eigenvalues $ \mu_j $  
could be not purely imaginary,  i.e. $ r_j $ could have a non-zero real part which depends on the 
nonlinearity (unlike the reversible or Hamiltonian case, where $ r_j \in \ii \R $). 
In such a case, the invariant torus could be 
(partially)  hyperbolic. 
Since we  
do not control the real part of $ r_j $ (i.e. the hyperbolicity may vanish), 
we perform the measure estimates 
proving the diophantine lower bounds   
of the imaginary part of the small divisors.

\smallskip

The final comment concerns the dynamical consequences of Theorem \ref{thm:reducibility}-$ii$).
All the above transformations   (both the changes of variables of 
sections \ref{step-1}-\ref{step-5} as well as the KAM matrices  of the reducibility scheme)
are  time-dependent  quasi-periodic maps of the phase space 
(of functions of $ x$ only), see section \ref{sec: dyn redu}.
It is thanks to this  ``T\"oplitz-in-time" structure that   
 the linear KdV equation \eqref{KdV:lin} is transformed into the dynamical system 
 \eqref{Lin: Red}. Note that in \cite{Ioo-Plo-Tol} (and also \cite{B5}, \cite{BBo10},\cite{BB12})  
the analogous transformations
have not this T\"oplitz-in-time structure and
stability informations are not obtained.

 \smallskip

\emph{Acknowledgements.} 
We warmly thank W. Craig for many discussions about the reduction approach of the
linearized operators and the reversible structure, and
P. Bolle for deep observations about the Hamiltonian case. We also thank
T. Kappeler, 
M. Procesi for many useful comments. 
 
\section{Functional setting}\label{sec:2} 

For a function $f : \Lm_o \to E$, $\lm \mapsto f(\lm)$, where $(E, \| \ \|_E)$ is a Banach space and 
$ \L_o $ is a subset of $\R$, we define the sup-norm and the Lipschitz semi-norm
\be \label{def norma sup lip}
\| f \|^{\sup}_E 
:= \| f \|^{\sup}_{E,\L_o} 
:= \sup_{ \lm \in \Lm_o } \| f(\lm) \|_E \, , 
\quad
\| f \|^{\lip}_E 
:= \| f \|^{\lip}_{E,\Lm_o}  
:= \sup_{\begin{subarray}{c} \lm_1, \lm_2 \in \Lm_o \\ \lm_1 \neq \lm_2 \end{subarray}} 
\frac{ \| f(\lm_1) - f(\lm_2) \|_E }{ | \lm_1 - \lm_2 | }\,,
\ee
and, for $ \g > 0 $, the Lipschitz norm
\be \label{def norma Lipg}
\| f \|^{\Lipg}_E  
:= \| f \|^\Lipg_{E,\Lm_o}
:= \| f \|^{\sup}_E + \g \| f \|^{\lip}_E  \, . 
\ee
If $ E = H^s $ we simply denote $ \| f \|^{\Lipg}_{H^s} := \| f \|^{\Lipg}_s $.

As a notation,  we write  
$$
a \leq_s b \quad \ \Longleftrightarrow \quad a \leq C(s) b
$$ 
for some constant $ C(s) $.  For $ s = \mathfrak s_0 := (\nu+2) \slash 2 $  we only write $ a \lessdot b $. 
More in general the notation
$ a \lessdot b  $ means $ a \leq C b $ 
where
the constant $ C $ may depend on the data of the problem, namely the nonlinearity  
$ f $,  the number $ \nu $ of frequencies, the diophantine vector $ \bar \om $, 
the diophantine exponent $ \tau >  0 $  in the non-resonance conditions in \eqref{Omegainfty}.
Also the small constants $ \d $ in the sequel depend on the data of the problem.

\subsection{Matrices with off-diagonal decay}

Let $ b \in \N $ and consider the exponential basis $\{ e_i : i \in \Z^b \} $ 
of $L^2(\T^b) $,
so that $L^2(\T^b)$ is the vector space $\{ u = \sum u_i e_i$, $\sum |u_i |^2 < \infty \}$.
Any linear operator $A : L^2 (\T^b) \to L^2 (\T^b) $ can be represented by the infinite dimensional matrix
\[
( A_{i}^{i'} )_{i, i' \in \Z^b}, \quad 
A_{i}^{i'} := ( A e_{i'}, e_{i})_{L^2(\T^b)}, \quad  
A u = \sum_{i, i'} A_{i}^{i'} u_{i'} e_{i}. 
\]
We now define  
the $ s $-norm (introduced in \cite{BBo10}) of an infinite dimensional matrix. 
\begin{definition}\label{def:norms}
The $s$-decay norm of an infinite dimensional matrix $ A := (A_{i_1}^{i_2} )_{i_1, i_2 \in \Z^b } $ is  
\begin{equation} \label{matrix decay norm}
\left| A \right|_{s}^2 := 
\sum_{i \in \Z^b} \left\langle i \right\rangle^{2s} 
\Big( \sup_{ \begin{subarray}{c} i_{1} - i_{2} = i 
\end{subarray}}
| A^{i_2}_{i_1}| \Big)^{2} \, .
\end{equation}
For parameter dependent matrices $ A := A(\l) $, $\l \in \L_o \subseteq \R$, 
the definitions \eqref{def norma sup lip} and \eqref{def norma Lipg} become 
\[
| A |^{\sup}_s  := \sup_{ \lm \in \Lm_o } | A(\lm) |_s \, , 
\quad
| A |^{\lip}_s := \sup_{\lm_1 \neq \lm_2} 
\frac{ | A(\lm_1) - A(\lm_2) |_s }{ | \lm_1 - \lm_2 | }\,,
\quad
| A |^{\Lipg}_s := | A |^{\sup}_s + \g | A |^{\lip}_s  \,. 
\]
\end{definition}
Clearly, the matrix decay norm \eqref{matrix decay norm} is increasing with respect to the index $ s  $, 
namely 
$$
| A |_s \leq | A |_{s'} \, ,  \quad \forall s < s'.
$$
The $ s $-norm is designed to estimate the polynomial  off-diagonal decay of matrices, actually it implies
$$
|A_{i_1}^{i_2}| \leq \frac{|A|_s}{\langle i_1 - i_2 \rangle^s} \, , \quad \forall i_1, i_2 \in \Z^b \, , 
$$
and, on the diagonal elements,  
\be\label{Aii}
|A_i^i | \leq |A|_0 \, , \quad |A_i^i |^{\rm lip} \leq |A|_0^{\rm lip} \, . 
\ee
We now list some properties of the matrix decay norm proved in \cite{BBo10}.
\begin{lemma}\label{lem:multi} {\bf (Multiplication operator)}
Let $ p = \sum_i  p_i e_i \in H^s(\T^b)$. 
The multiplication operator $ h \mapsto p h$ is represented by the T\"oplitz matrix 
$ T_i^{i'} = p_{i - i'} $ and 
\be\label{multiplication}
|T|_s = \| p \|_s.
\ee
Moreover, if $p = p(\l)$ is a Lipschitz family of functions, 
\be\label{multiplication Lip}
|T|_s^\Lipg = \| p \|_s^\Lipg\,.
\ee
\end{lemma}

The $s$-norm  satisfies classical algebra and interpolation inequalities. 
 
\begin{lemma} \label{prodest}
{\bf (Interpolation)} 
For all $s \geq s_0 > b/2 $ there are $ C(s) \geq C(s_0) \geq 1 $ such that
\be \label{interpm}
| A B|_{s} \leq  C(s) |A|_{s_0} |B|_s + C(s_0) |A|_s |B|_{s_0} \, .
\ee
In particular, the algebra property holds
\be \label{algebra}
|A B |_s \leq  C(s) |A|_s |B|_s \, .
\ee
If  $A = A(\lambda)$ and $B = B(\lambda)$ depend in a Lipschitz way on the parameter $\lambda \in 
\Lambda_o \subset \R$, then 
\begin{align} \label{algebra Lip}
|A B |_s^{\Lipg} 
& \leq  C(s) |A|_s^{\Lipg} |B|_s^{\Lipg} \, , 
\\
\label{interpm Lip}
|A B|_{s}^{\Lipg} 
& \leq C(s) |A|_{s}^{\Lipg} |B|_{s_0}^{\Lipg} 
+ C(s_0) |A|_{s_0}^{\Lipg} |B|_{s}^{\Lipg} .
\end{align}
\end{lemma}

For all $n \geq 1$, using \eqref{algebra} with $ s = s_0 $, we get  
\be\label{Mnab}
| A^n |_{s_0} \leq [C(s_0)]^{n-1} | A |_{s_0}^n \qquad
\text{and}  \qquad   
| A^n |_{s} \leq  n [ C(s_0) |A|_{s_0} ]^{n-1} C(s) | A |_{s} \, , \ \forall s \geq s_0 \, .
\ee
Moreover \eqref{interpm Lip} implies that \eqref{Mnab} also holds for Lipschitz norms $| \ |_s^\Lipg$. 

The $ s $-decay norm controls the Sobolev norm, also for Lipschitz families:
\be\label{interpolazione norme miste}
\| A h \|_s \leq C(s) \big(|A|_{s_0} \| h \|_s +  |A|_{s} \| h \|_{s_0} \big), 
\ \
\| A h \|_s^\Lipg 
\leq C(s) \big(|A|_{s_0}^\Lipg \| h \|_s^\Lipg + |A|_{s}^\Lipg \| h \|_{s_0}^\Lipg \big).
\ee
\begin{lemma}\label{lem:inverti}
Let $ \Phi = I + \Psi $ with $\Psi := \Psi(\l)$, depending in a Lipschitz way on the parameter $\l \in \L_o \subset \R $, 
such that  $ C(s_0) | \Psi |_{s_0}^{\Lipg} \leq 1/ 2 $. 
Then $ \Phi $ is invertible and, for all $ s \geq s_0 > b / 2  $, 
\be\label{PhINV}
| \Phi^{-1} - I |_s \leq C(s) | \Psi |_s \, , \quad
| \Phi^{-1} |_{s_0}^{\Lipg} \leq 2 \, , \quad 
| \Phi^{-1} - I |_{s}^{\Lipg} \leq C(s) | \Psi |_{s}^{\Lipg} \, . 
\ee
If $ \Phi_i = I + \Psi_i $, $  i = 1,2 $,  satisfy $ C(s_0) | \Psi_i |_{s_0}^{\Lipg} \leq 1/ 2 $, 
then 
\begin{equation}\label{derivata-inversa-Phi}
\vert \Phi_2^{-1} - \Phi_1^{-1} \vert_{s} 
\leq C(s) 
\big(  \vert \Psi_2 - \Psi_1 \vert_{s}  
+ \big( \vert \Psi_1 \vert_s + \vert \Psi_2 \vert_s \big)  
\vert \Psi_2 - \Psi_1 \vert_{s_0} \big) \, . 
\end{equation}
\end{lemma}

\begin{pf} 
Estimates \eqref{PhINV} follow by Neumann series and \eqref{Mnab}.
To prove \eqref{derivata-inversa-Phi}, observe that 

\[
\Phi_2^{-1} - \Phi_1^{-1} = \Phi_1^{-1} (\Phi_1 - \Phi_2) \Phi_2^{-1} = 
\Phi_1^{-1} (\Psi_1 - \Psi_2) \Phi_2^{-1}
\]
and use \eqref{interpm}, \eqref{PhINV}.  
\end{pf}

\subsubsection{T\"oplitz-in-time matrices}

Let now $ b := \nu + 1 $ and 
$$
e_i (\vphi, x) := e^{\ii (l \cdot \vphi + j x)},   \quad i := (l, j) \in \Z^b , \quad l \in \Z^\nu, \quad j \in \Z \, . 
$$
An important sub-algebra of matrices  is formed by the matrices T\"oplitz in time defined by
\be\label{Topliz matrix}
 A^{(l_2, j_2)}_{(l_1, j_1)}  := A^{j_2}_{j_1}(l_1 - l_2 )\,  ,
\ee
whose  decay norm \eqref{matrix decay norm} is
\be\label{decayTop}
|A|_s^2 =  \sum_{j \in \Z, l \in \Z^\nu} \sup_{j_1 - j_2 = j} |A_{j_1}^{j_2}(l)|^2  \langle l,j \rangle^{2 s} \, . 
\ee
These matrices are identified with the $ \vphi $-dependent family
of operators
\be\label{Aphi}
A(\vphi) := \big( A_{j_1}^{j_2} (\vphi)\big)_{j_1, j_2 \in \Z} \, , \quad 
A_{j_1}^{j_2} (\vphi) := \sum_{l \in \Z^\nu} A_{j_1}^{j_2}(l) e^{\ii l \cdot \vphi}
\ee
which act on functions of the $x$-variable 
as
\be\label{notationA}
A(\vphi) : h(x) = \sum_{j \in \Z} h_j e^{\ii jx} \mapsto  
A(\vphi) h(x) = \sum_{j_1, j_2 \in \Z}  A_{j_1}^{j_2} (\vphi) h_{j_2} e^{\ii j_1 x} \, .  
\ee
We still denote by $ | A(\vphi) |_s $ the $ s $-decay norm of the matrix in \eqref{Aphi}.

\begin{lemma}\label{Aphispace}
Let $ A $ be a T\"oplitz matrix as in \eqref{Topliz matrix}, and $\mathfrak s_0 := (\nu + 2)/2$ (as defined above). Then
$$
|A(\vphi)|_{s} \leq C(\mathfrak s_0) |A|_{s+ \mathfrak s_0} \, ,  \quad \forall  \vphi \in \T^\nu \, . 
$$
\end{lemma}

\begin{pf}
For all $ \vphi \in \T^\nu $ we have
\begin{eqnarray*}
|A(\vphi)|_{s}^2 & := & \sum_{j \in \Z} \langle j \rangle^{2 s}  \sup_{j_1 - j_2 = j} |A_{j_1}^{j_2}(\vphi)|^2 
 \lessdot  
 \sum_{j \in \Z} \langle j \rangle^{2 s} 
 \sup_{j_1 - j_2 = j} \sum_{l \in \Z^\nu} |A_{j_1}^{j_2}(l)|^2 \langle l \rangle^{2 {\mathfrak s}_0}  \\
 & \lessdot &  \sum_{j \in \Z} \sup_{j_1 - j_2 = j} \sum_{l \in \Z^\nu} |A_{j_1}^{j_2}(l)|^2  \langle l,j \rangle^{2 (s + {\mathfrak s}_0)}
 \lessdot  \sum_{j \in \Z, l \in \Z^\nu} \sup_{j_1 - j_2 = j} |A_{j_1}^{j_2}(l)|^2  \langle l,j \rangle^{2 (s + {\mathfrak s}_0)}\\ 
 & \stackrel{\eqref{decayTop}} \lessdot  & |A|_{s + {\mathfrak s}_0}^2 ,
\end{eqnarray*}
whence the lemma follows.
\end{pf}

Given $ N \in \N $, we define the smoothing operator $\Pi_N$ as
\be\label{SM}
\big(\Pi_N A \big)^{(l_2, j_2)}_{(l_1, j_1)} := 
\begin{cases}
A^{(l_2, j_2)}_{(l_1, j_1)} \qquad \, {\rm if} \  | l_1 - l_2| \leq N  \\
0 \quad \qquad \qquad {\rm otherwise.} 
\end{cases}
\ee
\begin{lemma} 
The operator $ \Pi_N^\bot := I - \Pi_N $ satisfies 
\be\label{smoothingN}
| \Pi_N^\bot A |_{s} \leq N^{- \b} |  A |_{s+\b} \, , \quad 
| \Pi_N^\bot A |_{s}^{\Lipg} \leq N^{- \b} |  A |_{s+\b}^{\Lipg} \, , 
\quad \b \geq 0,
\ee
where in the second inequality $A := A(\l)$ is a Lipschitz family $\l \in \L$. 
\end{lemma}

\subsection{Dynamical reducibility}\label{sec: dyn redu} 

All the transformations that we construct in sections \ref{sec:regu} and \ref{sec:redu} act on 
functions $ u(\vphi, x ) $ (of time and space). They  
can also be seen as: 
\begin{itemize}
\item[$(a)$] transformations of the phase space $H^s_x$ that depend quasi-periodically on time 
(sections \ref{step-1}, \ref{step-3}-\ref{step-5} and \ref{sec:redu}); 
\item[$(b)$] quasi-periodic reparametrizations of time (section \ref{step-2}).
\end{itemize}
This observation allows to interpret the conjugacy procedure from a dynamical point of view. 

Consider a quasi-periodic linear dynamical system
\be\label{SD}
\partial_t u = L(\om t)  u. 
\ee
We want to describe how \eqref{SD} changes under the action of a transformation of type $(a)$ or $(b)$. 

Let $A(\om t)$ be  of type $(a)$, and let $u = A(\om t)v$. 
Then \eqref{SD} is transformed  
into the linear system 
\begin{equation}\label{sistematrasformato}
\partial_{t} v = L_+(\om t)v \quad {\rm where} \quad
L_{+}(\om t) = A(\om t)^{-1} L(\om t) A(\om t) - A(\om t)^{-1} \partial_t A(\om t) \, . 
\end{equation}
The transformation $A(\om t)$ 
may be regarded to act on functions $ u(\vphi, x) $ as
\be\label{trasfo spazio}
({\tilde A} u)(\vphi,x) := \big(A(\vphi)u(\vphi, \cdot )\big) (x)  := A(\vphi)u(\vphi, x) 
\ee
and one can 
check that  $ ({\tilde A}^{-1} u)(\vphi,x) = A^{-1}(\vphi) u(\vphi, x) $.  
The operator associated to \eqref{SD} (on quasi-periodic functions) 
\be\label{associa}
{\cal L}  := \om \cdot \partial_\vphi - L(\vphi) 
\ee
transforms under the action of $ {\tilde A} $ into 
$$
{\tilde A}^{-1} {\cal L} {\tilde A} = \om \cdot \partial_\vphi - L_+(\vphi), 
$$
which is exactly the linear system in \eqref{sistematrasformato}, 
acting on quasi-periodic functions. 

\smallskip

Now consider a transformation of type $(b)$, namely a change of the time variable
\be\label{time repar}
\t := t + \a(\om t) \  \Leftrightarrow \ 
t = \t + \tilde \a (\om \t);  \quad (Bv)(t) := v(t + \a(\om t)),
\ \ 
(B^{-1} u)(\t) = u(\t + \tilde \a(\om \t)),
\ee
where $\a = \a(\ph)$, $\ph \in \T^\nu$, is a $2\p$-periodic function of $\nu$ variables
(in other words, $ t \mapsto  t + \a(\om t) $ is the diffeomorphisms of $\R$ induced by the transformation $B$).
If $ u(t) $ is a solution of \eqref{SD}, then $ v(\tau) $, defined by 
$ u = Bv$,  
solves
\be\label{timeqpDS}
\partial_\tau v(\t) = L_+ (\om \tau) v (\tau) \, , \quad    
L_+ (\om \tau) := \Big( \frac{L(\om t)}{1+ (\ompaph \a) (\om t)} \Big)_{|t= \t + \tilde \a (\om \t)} \,.  
\ee
We may regard the associated transformation on quasi-periodic functions defined by
$$
(\tilde B h)(\vphi,x) := h( \vphi + \om \a (\vphi ), x) \, , \quad 
(\tilde B^{-1} h)(\vphi,x) := h( \vphi + \om \tilde \a (\vphi ), x) \, , 
$$
as in step \ref{step-2}, where we calculate 
$$
B^{-1} {\cal L} B = \rho(\vphi) {\cal L}_+ \, , 
\quad  \rho(\vphi) := B^{-1} (1+ \om \cdot \partial_\vphi \a) \, , 
$$
\be\label{timeqp}
{\cal L}_+ =  \om \cdot \partial_\vphi - L_+(\vphi) \, , \ \ 
L_+(\vphi)  := \frac{1}{\rho(\vphi)} L(\vphi + \om {\tilde \a}(\vphi )) \, . 
\ee
\eqref{timeqp} is nothing but the linear system \eqref{timeqpDS}, acting on quasi-periodic functions. 

\subsection{Real, reversible and Hamiltonian operators}

We consider the space of {\it real} functions
\begin{equation}\label{funzionireali}
Z := \{ u(\vphi,x) = \overline{u(\vphi,x)} \}, 
\end{equation}
and of 
 {\it even} (in space-time), respectively  
 {\it odd}, functions 
 \begin{equation}\label{funzionipari}
X := \{ u(\vphi,x) = u(-\vphi,-x) \} \, , \quad Y := \{ u(\vphi,x) = -u(-\vphi,-x) \} \, .
\end{equation}

\begin{definition}\label{def:RR}
An operator $ R $ is
\begin{enumerate}
\item {\sc real} if $ R : Z \to Z $
\item {\sc reversible} if $ R : X \to Y $
\item {\sc reversibility-preserving} if $ R : X \to X $, $ R : Y \to Y $.
\end{enumerate}
\end{definition}
The composition of a reversible and a reversibility-preserving operator is reversible. 

The above properties may be  characterized in terms of matrix elements. 
\begin{lemma} \label{lem:PR} We have
$$ 
R : X \to Y \  \Longleftrightarrow  \  R^{-j}_{-k}(-l) = - R^j_{k}(l) \, , \qquad 
R : X \to X \  \Longleftrightarrow  \ R^{-j}_{-k}(-l) =  R^j_k (l) \, , 
$$
$$
R : Z \to Z \quad \Longleftrightarrow  \quad \overline{R^j_{k}(l)} = R^{-j}_{-k}(-l) \, .  
$$
\end{lemma} 

For the Hamiltonian KdV the phase space is $ H^1_0 := \{ u \in H^1 (\T) \, : \, \int_{\T} u(x) dx = 0 \} $
and it is more convenient the dynamical systems perspective. 

\begin{definition}
A time dependent linear vector field  $ X(t) : H_0^1 \to H_0^1$ is \textsc{Hamiltonian} if 
$ X(t) = \partial_x G(t) $ for some  
real linear operator $ G(t) $ which is self-adjoint with respect to the $ L^2 $ scalar product.  

If $ G(t) = G(\om t)$ is quasi-periodic in time, we say that the associated operator 
$ \om \cdot \partial_{\vphi} - \partial_x G( \vphi ) $ (see \eqref{associa}) is Hamiltonian. 
\end{definition}

\begin{definition}
A map $ A : H_0^1 \to H_0^1$ is  
 \textsc{symplectic} if 
 \be\label{mappa simplettica}
 \Omega(A u, A v) = \Omega (u, v) \, , \quad \forall u,v \in H_0^1 \, , 
\ee 
where the symplectic 2-form $ \Omega $ is defined in \eqref{KdV symplectic}. 
Equivalently $ A^T \partial_x^{-1} A = \partial_x^{-1}  $.

If $ A (\vphi ) $, $ \forall \vphi \in \T^\nu $,  is a family of symplectic maps we say that 
the corresponding operator  in \eqref{trasfo spazio} 
is symplectic. 

\end{definition}

Under a time dependent family of symplectic transformations
 $ u = \Phi(t) v   $
the linear  Hamiltonian 
equation
$$
u_t = \partial_x G(t) u  \quad {\rm with \ Hamiltonian} \quad H(t, u) := \tfrac12 \, \big(G(t)u ,u \big)_{L^2}
 $$
transforms into the 
equation 
\[
v_t = \pa_x E(t) v, \quad
E(t) := \Phi(t)^T G(t) \Phi(t) - \Phi(t)^T \pa_x^{-1} \Phi_t(t) 
\]
with Hamiltonian 
\be\label{transformed KdV}
K(t,v) =  \tfrac12\, \big( G(t) \Phi (t)  v ,  \Phi(t) v \big)_{L^2} 
- \tfrac12\, \big( \partial_x^{-1} \Phi_t(t)v,  \Phi (t) v \big)_{L^2}  \, .
\ee
Note that $E(t)$ is self-adjoint with respect to the $L^2$ scalar product because 
$\Phi^T \pa_x^{-1} \Phi_t + \Phi_t^T \pa_x^{-1} \Phi = 0$. 

\section{Regularization of the linearized operator}\label{sec:regu}

Our existence proof is based on a Nash-Moser iterative scheme. 
The main step concerns the invertibility of the linearized operator (see \eqref{linearized op}) 
\begin{equation} \label{mL}
\mL h = \mL(\lm,u,\e) h := \ompaph  h + (1 + a_3) \pa_{xxx} h + a_2  \pa_{xx} h + a_1 \pa_{x} h + a_0 h
\end{equation}
obtained linearizing \eqref{eq:invqp} at any approximate (or exact) solution $ u $. The 
coefficients $a_i = a_i(\ph,x) = a_i(u,\e)(\ph,x)$ are periodic functions of $(\ph,x)$, depending on $u,\e$. They are 
explicitly obtained from the partial derivatives of $\e f(\ph,x,z)$ as
\begin{equation}  \label{ai formula}
a_i(\ph,x) = \e (\pa_{z_i} f)\big( \ph, x, u(\ph,x), u_x(\ph,x), u_{xx}(\ph,x), u_{xxx}(\ph,x) \big), 
\quad i=0,1,2,3.
\end{equation}
The operator $\mL$ depends on $\lm$ because $\om = \lm \bar\om$. 
Since $\e$ is a (small) fixed parameter, we simply write $\mL(\lm,u)$ instead of $\mL(\lm,u,\e)$, 
and $a_i(u)$ instead of $a_i(u,\e)$.  
We emphasize that the coefficients $a_i$ do not depend explicitely on the parameter $\lm$ (they 
depend on $ \l $ only through  $ u(\l) $).

\smallskip

In the Hamiltonian case \eqref{f Ham} the linearized KdV operator \eqref{mL} has the form
$$
{\cal L}h = \ompaph  h + \partial_{x} \Big( \partial_x \big\{ A_1 (\vphi,x) \partial_x h \big\}
- A_0 (\vphi,x) h \Big)
$$
where
$$
A_1 (\vphi,x) := 1 + \e (\partial_{z_1 z_1} F) (\vphi,x,u,u_x)\, , 
\quad 
A_0 (\vphi,x) := - \e \pa_x \{ (\pa_{z_0 z_1} F)(\ph,x,u,u_x) \} + \e (\partial_{z_0 z_0} F) (\vphi,x,u,u_x)
$$
and it is generated by the quadratic Hamiltonian
$$
H_L(\vphi, h) := \frac{1}{2} \int_{\T} \Big( A_0 (\vphi, x) h^2 + A_1 (\vphi, x) h_x^2 \Big) \, dx \,, 
\quad  h \in H^1_0 \,.
$$
\begin{remark}
In the reversible case, i.e.  the nonlinearity $ f$ satisfies 
\eqref{parity f} and $ u \in X $ (see \eqref{funzionipari}, \eqref{solPP}) the coefficients $ a_i $ satisfy the parity
\be\label{a3a1a2a0}
a_3, a_1 \in X, \quad a_2, a_0 \in Y,
\ee
and  $\mL$ maps $X$ into $Y$, namely $\mL$ is reversible, see Definition \ref{def:RR}. 
\end{remark}

\begin{remark} \label{rem:coeff in cases Q F}
In the Hamiltonian case \eqref{f Ham}, assumption (Q)-\eqref{type Q} is automatically satisfied (with $ \a (\vphi ) = 2 $) because  
$$
f(\ph,x,u,u_x, u_{xx}, u_{xxx}) = a(\ph, x, u, u_x) 
+ b(\ph, x, u, u_x) u_{xx}  + c(\ph, x, u, u_x) u_{xx}^2  + d(\ph, x, u, u_x) u_{xxx}   
$$
where 
$$ 
b =  2 (\partial_{z_1 z_1 x}^3 F) + 2 z_1 (\partial_{z_1 z_1 z_0}^3 F),  
\qquad c  = \partial_{z_1}^3 F,   \qquad d = \partial_{z_1}^2 F, 
$$ 
and so
$$ 
\partial_{z_2} f =  b + 2 z_2 c = 2(d_x + z_1 d_{z_0} + z_2 d_{z_1}) = 
2 \Big(  \pa^2_{z_3 x} f + z_1 \pa^2_{z_3 z_0} f + z_2 \pa^2_{z_3 z_1} f + z_3 \pa^2_{z_3 z_2} f  \Big) \, . 
$$ 
\end{remark}

The coefficients $a_i$, together with their derivative $\pa_u a_i(u)[h]$ with respect to $u$ in the direction $h$,  
satisfy tame estimates:

\begin{lemma}  \label{lemma:stime ai} 
Let $ f \in C^q $, see \eqref{f classe Cq}. For all 
$ \mathfrak s_{0} \leq s \leq q - 2 $,  $  \| u \|_{\mathfrak s_0 + 3} \leq 1 $, 
we have, for all $i = 0,1,2,3$, 
\begin{align} 
\label{stima coeff ai 1}
\| a_i(u) \|_s 
& \leq \e \, C(s) \big( 1 + \| u \|_{s+3} \big),  
\\ 
\label{stima coeff ai 2}
\| \pa_{u} a_i(u)[h] \|_{s} 
& \leq \e \, C(s) \big( \| h \|_{s+3} + \| u \|_{s+3} \| h \|_{\mathfrak s_0+3} \big) \, .  
\end{align}
If, moreover, $ \l \mapsto u(\l) \in H^s $ is Lipschitz family satisfying  $ \| u \|_{\mathfrak s_0 + 3}^{\Lipg} \leq 1 $ (see \eqref{def norma Lipg}), then 
\be
\label{stima coeff ai 3}
\| a_i \|_{s}^{\Lipg}  \leq \e \, C(s) \big( 1 + \| u \|_{s+3}^{\Lipg} \big) \, . 
\ee
\end{lemma}

\begin{pf} 
The tame estimate  
\eqref{stima coeff ai 1} follows by 
Lemma \ref{lemma:composition of functions, Moser}$(i)$ applied to the function $\pa_{z_i}f$, $i=0,\ldots,3 $, 
which is valid for $s+1 \leq q$. The tame bound   
\eqref{stima coeff ai 2} for  
\[
\pa_u a_i(u) [h]  
\stackrel{\eqref{ai formula}} = \e \sum_{k=0}^3 (\pa^2_{z_k z_i} f)\big( \ph, x, u, u_x, u_{xx}, u_{xxx} \big) \, \pa_x^k h, 
\quad i = 0, \ldots, 3, 
\]
follows 
by \eqref{asymmetric tame product} and 
applying 
Lemma \ref{lemma:composition of functions, Moser}$(i)$  to the functions $\pa^2_{z_k z_i}f$, 
which  gives 
\[
\| (\pa^2_{z_k z_i} f)\big( \ph, x, u, u_x, u_{xx}, u_{xxx} \big) \|_s
\leq C(s) \| f \|_{C^{s+2}} (1 + \| u \|_{s+3}),
\] 
 for $s+2 \leq q$. The Lipschitz bound  \eqref{stima coeff ai 3} follows similarly.
\end{pf}

\subsection{Step 1. Change of the space variable} \label{step-1}

We consider a $ \vphi $-dependent family of  diffeomorphisms of the $ 1 $-dimensional torus 
$ \T $ of the form 
\begin{equation}\label{cambio1}
y  = x + \beta(\vphi,x), 
\end{equation}
where $ \beta $ is a (small) real-valued function, $2\pi$ periodic in all its arguments.
The change of variables (\ref{cambio1}) induces on the space of functions the linear operator 
\begin{equation}\label{operatore1}
({\cal A}h)(\vphi,x):= h(\vphi, x + \b (\vphi, x)). 
\end{equation}
The operator $ {\cal A} $ is invertible, with inverse
\begin{equation}\label{inverse}
({\cal A}^{-1} v)(\vphi,y) = v(\vphi, y  + {\tilde \beta}(\vphi,y) ),
\end{equation}
where $ y \mapsto y  + {\tilde \beta}(\vphi,y) $ is the inverse diffeomorphism of \eqref{cambio1}, namely
\be \label{INVDIF}
x = y  + {\tilde \beta}(\vphi,y)  \quad \Longleftrightarrow \quad  y  = x + \beta(\vphi,x).  
\ee
\begin{remark}\label{rem: Ham0}
In the Hamiltonian case \eqref{f Ham} we use, instead of \eqref{operatore1},  the modified 
change of variable 
\eqref{operatore1 simplettico}
which is symplectic, for each $ \vphi \in \T^\nu $.  
Indeed, setting $ U := \partial_x^{-1} u $ (and neglecting to write the $ \vphi $-dependence)
\begin{align*}
\Omega ({\cal A}u, {\cal A}v) 
& =  
\int_{\T} \partial_{x}^{-1} \Big( \pa_x  
\big\{ U(x+ \b (x) ) \big\} \Big) \, (1+ \b_x (x) ) v(x+ \b (x) ) \, dx 
\\ 
& = \int_{\T}  U(x+ \b (x))  (1+ \b_x (x) ) v(x+ \b (x) ) dx 
- c \int_{\T}  (1+ \b_x (x) ) v(x+ \b (x) ) dx 
\\
& = \int_{\T}   U(y) v(y ) dy = \Omega (u,v) \, , \quad v \in H^1_0 \, , 
\end{align*}
where $ c $ is the average of $ U(x+ \b (x) ) $ in $ \T $.
The inverse operator of \eqref{operatore1 simplettico} is 
$ ({\cal A}^{-1} v) (\vphi, y) = (1+ {\tilde \beta}_y (\vphi, y)) v( y + \tilde \beta (\vphi, y)) $ 
which  is also symplectic.
\end{remark}

\noindent
Now
we calculate the conjugate $ {\cal A}^{-1} {\cal L} {\cal A} $ of the linearized operator $\mL$ in \eqref{mL} 
with $ {\cal A} $ in \eqref{operatore1}.

The conjugate $ {\cal A} \inv a {\cal A} $ of any multiplication operator $a : h(\ph,x) \mapsto a(\ph,x) h(\ph,x)$ is the multiplication operator $( {\cal A} \inv a)$ that maps $v(\ph,y) \mapsto ( {\cal A} \inv a)(\ph,y) \, v(\ph,y)$.
By conjugation, the differential operators become 
\begin{align*}
{\cal A} \inv \ompaph {\cal A} & = \ompaph + \{ {\cal A} \inv(\ompaph \b) \} \, \pa_y, 
\\
{\cal A} \inv \pa_x {\cal A} & = \{ {\cal A} \inv(1 + \b_x) \} \, \pa_y, 
\\
{\cal A} \inv \pa_{xx} {\cal A} & = \{ {\cal A} \inv (1+\b_x)^2 \} \, \pa_{yy} + \{ {\cal A} \inv (\b_{xx}) \} \, \pa_y, 
\\
{\cal A} \inv \pa_{xxx} A  & = \{ {\cal A} \inv (1+\b_x)^3 \} \, \pa_{yyy} 
+ \{ 3 {\cal A} \inv[ (1+\b_x) \b_{xx}] \} \, \pa_{yy} 
+ \{ {\cal A} \inv (\b_{xxx}) \} \, \pa_y,
\end{align*}
where all the coefficients $\{ A\inv (\ldots) \}$ are periodic functions of $(\ph,y)$. 
Thus (recall \eqref{mL})
\begin{equation} \label{inv1s}
\mL_1 
:= {\cal A}^{-1} \mL {\cal A} 
= \ompaph + b_3(\ph,y) \pa_{yyy} 
+ b_2(\ph,y) \pa_{yy} 
+ b_1(\ph,y) \pa_{y} 
+ b_0(\ph,y)
\end{equation} 
where 
\begin{alignat}{2}
\label{b1 b3}
b_3 & = {\cal A} \inv[(1+a_3) (1+\b_x)^3], 
& \qquad 
b_1 & = {\cal A} \inv[\ompaph \b + (1+a_3) \b_{xxx} + a_2 \b_{xx} + a_1 (1+\b_x)], 
\\
b_0 & = {\cal A} \inv(a_0),
& \qquad 
b_2 & = {\cal A} \inv[(1+a_3) 3 (1+\b_x) \b_{xx} + a_2 (1+\b_x)^2].
\label{b0 b2}
\end{alignat}
We look for $\b(\ph,x)$ such that the coefficient $b_3(\ph,y)$ of the highest order derivative $\partial_{yyy}$ in \eqref{inv1s} does not depend on $y$, namely 
\be \label{b(ph) 1}
b_3(\ph,y) \stackrel{\eqref{b1 b3}} = {\cal A}^{-1} [(1+a_3) (1+\b_x)^3] (\ph,y) = b(\vphi)
\ee
for some function $b(\ph)$ of $\vphi$ only. 
Since ${\cal A}$ changes only the space variable, ${\cal A}b = b$ for every function $b(\ph)$ that is independent on $y$. 
Hence \eqref{b(ph) 1} is equivalent to 
\be\label{eq:ste1}
\big( 1 + a_3(\ph,x) \big) \big( 1 + \b_x(\ph,x) \big)^3 = b(\ph),
\ee
namely 
\begin{equation} \label{primaequazione}
\beta_{x} = \rho_0, \qquad 
\rho_0(\vphi,x) := b(\ph)^{1/3} \big( 1 + a_3(\ph,x) \big)^{-1/3} - 1. 
\end{equation} 
The equation \eqref{primaequazione} has a solution $\b$, periodic in $ x $, if and only if
$ \int_{\T}{\rho_0(\vphi,x) \, dx} = 0 $. This condition uniquely determines 
\begin{equation} \label{c}
b(\vphi) = \left( \frac{1}{2\pi}\int_{\T} \big( 1 + a_3(\vphi,x) \big)^{-\frac13} \, dx \right)^{-3}.
\end{equation}
Then we fix the solution (with zero average)  
of \eqref{primaequazione}, 
\be\label{defb1b0}
\beta(\vphi,x) := \, (\pa_x\inv \rho_0)(\ph,x) \, , 
\ee
where $ \partial_x^{-1} $ is defined by linearity as
\be\label{dx-1}
\partial_x^{-1} e^{\ii j x} := \frac{ e^{\ii j x} }{\ii j}\,  \quad \forall j \in \Z \setminus \{ 0 \}, \qquad \partial_x^{-1} 1 = 0.
\ee
In other words, $\pa_x\inv h$ is the primitive of $h$ with zero average in $x $.

With this choice of $ \b $, we get (see \eqref{inv1s}, \eqref{b(ph) 1}) 
\begin{equation}  \label{mL1}
\mL_1 = {\cal A}\inv \mL {\cal A}  
= \ompaph + b_3(\ph) \pa_{yyy} + b_2(\ph,y) \pa_{yy} + b_1(\ph,y) \pa_y + b_0(\ph,y),
\end{equation} 
where $ b_3(\ph) := b(\ph) $ is defined in $ \eqref{c} $.

\begin{remark}\label{reversibilitˆ step 1} 
In the reversible case, 
$ \b \in Y $ because $a_3 \in X$, see \eqref{a3a1a2a0}. 
Therefore the operator 
$ A $ in \eqref{operatore1}, as well as $ {\cal A}^{-1} $ in \eqref{inverse}, 
maps $ X \to X $ and $ Y \to Y $, namely it is  reversibility-preserving, see Definition \ref{def:RR}.
By \eqref{a3a1a2a0} the coefficients of $\mL_1$ (see \eqref{b1 b3}, \eqref{b0 b2}) 
have parity 
\be\label{b3b1b2b0}
b_3, b_1 \in X, \qquad b_2, b_0 \in Y,
\ee
and $\mL_1$ maps $X \to Y$, namely it is reversible. 
\end{remark}

\begin{remark}\label{rem: Ham1}
In the Hamiltonian case \eqref{f Ham} 
the resulting operator $ {\cal L}_1 $ in \eqref{mL1} is  Hamiltonian and $ b_2 (\vphi, y) = 2 \partial_y b_3 (\vphi) \equiv 0 $.
Actually, by  \eqref{transformed KdV}, the corresponding Hamiltonian has the form  
\be\label{Ham K}
K(\vphi, v) = \frac{1}{2} \int_{\T} b_3(\vphi ) v_y^2 + B_0 (\vphi ,y) v^2\, dy \, ,   
\ee
for some function $ B_0  (\vphi ,y) $. 
\end{remark}

\subsection{Step 2. Time reparametrization} \label{step-2}

The goal of this section is to make constant the coefficient of the highest 
order spatial derivative operator $\partial_{yyy}$ of $ {\cal L}_1 $ in \eqref{mL1}, 
by a quasi-periodic reparametrization of time. 
We consider a diffeomorphism of the torus  $ \T^{\nu} $
of the form
\begin{equation}\label{cambio2}
\vphi \mapsto \vphi + \omega \a(\vphi),   \quad  \vphi \in \T^{\nu}, \quad \a (\vphi ) \in \R \, ,  
\end{equation}
where $ \a $ is a (small) {\it real} valued function, $ 2\pi $-periodic in all its arguments.
The induced  linear operator on the space of functions is
\begin{equation}\label{operatore2}
(Bh)(\vphi,y):= h \big( \vphi + \omega \a(\vphi), \,y \big) 
\end{equation}
whose  inverse is
\be\label{B-1}
(B^{-1} v)(\th,y):= v \big( \th + \omega {\tilde \a}(\th), \,y \big) 
\ee
where $ \vphi = \th + \omega {\tilde \a}(\th) $ is the inverse diffeomorphism of $ \th = \vphi + \omega \a(\vphi) $. 
By conjugation, the differential operators become 
\begin{equation}  \label{anche def rho}
B\inv \ompath B = \rho(\th)\,  \ompath ,
\quad 
B\inv \pa_y B = \pa_y, 
\quad 
\rho := B\inv(1 + \ompaph \a).
\end{equation}
Thus, see  \eqref{mL1}, 
\begin{equation}\label{L1B}
B\inv \mL_1 B = \rho \,  \ompath 
+ \{ B\inv b_3 \} \, \pa_{yyy} 
+ \{ B\inv b_2 \} \, \pa_{yy} 
+ \{ B\inv b_1 \} \, \pa_{y} 
+ \{ B\inv b_0 \} .
\end{equation}
We look for $\a(\ph)$ such that the (variable) coefficients of the highest order derivatives ($\ompath$ and $\pa_{yyy}$) are proportional, namely
\be \label{B-1b3}
\{ B\inv b_3\}(\th) 
= \muff_3 \rho(\th) 
= \muff_3  \{ B\inv(1 + \ompaph \a)\}(\th)
\ee
for some constant $\muff_3 \in \R$. 
Since $ B $ is invertible, this is equivalent to require that
\begin{equation} \label{proportional}
b_3(\ph) = \muff_3 \big( 1 + \ompaph \a(\ph) \big).
\end{equation}
Integrating on $\T^\nu$ determines the value of the constant $\muff_3$, 
\begin{equation}  \label{mu 3}
\muff_3 := \frac{1}{(2\p)^\nu} \, \int_{\T^\nu} b_3(\ph) \, d\ph.
\end{equation} 
Thus we choose the unique solution of \eqref{proportional}
with zero average 
\begin{equation}  \label{alpha}
\a(\ph) := \frac{1}{\muff_3} \, (\ompaph)\inv (b_3 - \muff_3)(\ph) 
\end{equation} 
where $ (\ompaph)\inv $ is defined by linearity 
$$
(\ompaph)\inv e^{\ii l \cdot \vphi} := \frac{e^{\ii l \cdot \vphi}}{\ii \om \cdot l} \, , \ l \neq 0 \, , \quad 
(\ompaph)\inv 1 = 0 \, . 
$$
With this choice of $ \a $ we get (see \eqref{L1B}, \eqref{B-1b3})
\begin{equation}\label{mL2}
B\inv \mL_1 B = \rho \, \mL_2, 
\qquad 
\mL_2 := \ompath + \muff_3 \, \pa_{yyy} 
+ c_2(\th,y) \, \pa_{yy} 
+ c_1(\th,y) \, \pa_{y} 
+ c_0(\th,y),
\end{equation}
where
\begin{equation}\label{coefficienti mL2}
c_i := \frac{B\inv b_i}{\rho}\,, \quad i = 0,1,2.
\end{equation}

\begin{remark}\label{reversibilitˆ step 2}
In the reversible case, $\a$ is odd because $b_3$ is even (see \eqref{b3b1b2b0}), 
and $ B $ is reversibility preserving.  Since $\rho $ (defined in \eqref{anche def rho}) is even, 
the coefficients 
$ c_3, c_1 \in X$, $ c_2, c_0 \in Y $
and 
$\mL_2 : X \to Y $ is reversible. 
\end{remark}

\begin{remark}\label{rem: Ham2}
In the Hamiltonian case, the operator $ \mL_2 $ is still Hamiltonian
(the new Hamiltonian  is the old one at the new time, 
divided by the factor $ \rho $).
The coefficient $ c_2 (\vartheta, y) \equiv 0 $ because $ b_2 \equiv 0 $, see remark \ref{rem: Ham1}.
\end{remark}

\subsection{Step 3. Descent method: step zero}\label{step-3}

The aim of this section is to eliminate the term of order $\pa_{yy}$ from $ {\cal L}_2 $ in \eqref{mL2}.

Consider the multiplication operator 
\begin{equation}\label{cambio3}
{\cal M} h :=  v(\th , y) h 
\end{equation}
where the function $v$ is periodic in all its arguments. 
  Calculate the difference 
\begin{equation}
{\cal L}_2 \, {\cal M}   - {\cal M}  \, (\o \cdot \partial_{\th} + m_3 \partial_{yyy}) 
= T_2 \partial_{yy} + T_{1} \partial_{y} + T_{0},  
\label{carmelino}
\end{equation}
where 
\begin{equation} \label{T2} 
T_{2} := 3 \muff_3 v_{y} + c_{2} v,  
\quad 
T_{1} := 3 \muff_3 v_{yy} + 2 c_{2} v_{y} + c_{1} v, \quad 
T_{0} := \o\cdot\partial_{\th} v + \muff_3 v_{yyy}  + c_{2} v_{yy} +  c_{1} v_{y} + c_0 v .
\end{equation}
To eliminate the factor $T_2$, we need 
\be\label{sted}
3 \muff_3 v_{y} + c_{2} v = 0.
\ee
Equation \eqref{sted} has the periodic solution 
\begin{equation} \label{descent ordine zero}
v(\th,y) = \exp \Big\{ - \frac{1}{3\muff_3} \, (\pa_y\inv c_2)(\th,y) \Big\}
\end{equation}
provided that 
\be\label{zeromean}
\int_\T c_2(\th,y) \, dy = 0.
\ee
Let us prove \eqref{zeromean}. 
By \eqref{coefficienti mL2}, \eqref{anche def rho}, for each $\th = \ph + \om \a(\ph)$ we get
\[
\int_\T c_2(\th,y) \, dy 
= \frac{1}{ \{ B\inv(1 
+ \ompaph \a)\}(\th)} \, 
\int_\T (B\inv b_2)(\th,y) \, dy 
= \frac{1}{ 1 + \ompaph \a(\ph) } \, 
\int_\T b_2(\ph,y) \, dy.
\]
By the definition \eqref{b0 b2} of $b_2$ and 
 changing variable $ y = x + \b(\ph,x) $ in the integral (recall \eqref{operatore1}) 
\begin{align}
\int_\T b_2(\ph,y) \, dy  
& \stackrel{\eqref{b0 b2}} = \int_\T \Big( (1+a_3) 3 (1+\b_x) \b_{xx} + a_2 (1+\b_x)^2 \Big)   \, (1 + \b_x) \, dx \nonumber
\\ 
& \stackrel{ \eqref{eq:ste1}} = b(\ph) \Big\{ 3 \int_\T \frac{ \b_{xx}(\ph,x)}{1 + \b_x(\ph,x)} \, dx 
+ \int_\T \frac{ a_2(\ph,x) }{ 1 + a_3(\ph,x) } \, dx \Big\}. \label{viapez}
\end{align}
The first integral in \eqref{viapez} is zero because $\b_{xx} / (1 + \b_x) = \pa_x \log (1 + \b_x)$. 
The second one is zero because of assumptions (Q)-\eqref{type Q} or (F)-\eqref{type F}, see \eqref{zero mean in the intro}.
As a consequence \eqref{zeromean} is proved, and 
 \eqref{sted} has the periodic solution $v$ defined in \eqref{descent ordine zero}. 
Note that  $ v $ 
is close to $ 1 $ for  $ \e $ small. 
Hence the multiplication operator $ {\cal M}  $ defined in \eqref{cambio3} is invertible and $ {\cal M} ^{-1} $ is the multiplication
operator for $ 1 / v $. By \eqref{carmelino} and since $T_2 = 0$, we deduce
\begin{equation} \label{mL3}
\mL_3 := {\cal M} ^{-1} \mL_2 {\cal M}  
= \ompath + m_3 \pa_{yyy} 
+ d_{1}(\th , y)\partial_y + d_{0}(\th , y) , 
\qquad d_i := \frac{T_i}{v},\quad i=0,1.
\end{equation}

\begin{remark}\label{reversibilitˆ step 3}
In the reversible case, since $c_2$ is odd (see Remark \ref{reversibilitˆ step 2} ) the function $v$ is even, then 
$ {\cal M}  $, $ {\cal M} ^{-1}$ are reversibility preserving and by \eqref{T2} and \eqref{mL3} $d_1 \in X$ and $d_0 \in Y$, which implies that $\mL_3 : X \rightarrow Y$.
\end{remark}

\begin{remark}\label{rem: Ham3}
In the Hamiltonian case, there is no need to perform this step because $ c_2  \equiv 0 $,
see remark \ref{rem: Ham2}.
\end{remark}

\subsection{Step 4. Change of space variable (translation)}\label{step-4}

Consider the change of the space variable 
$$
z = y + p(\th)
$$
which induces the operators 
\be\label{MM-1}
{\cal T} h(\th,y) := h(\th, y + p(\th)), \quad 
{\cal T}\inv v(\th,z) := v(\th, z - p(\th)).
\ee
The differential operators become
\[
{\cal T}\inv \ompath {\cal T} = \ompath + \{ \ompath p(\th) \} \, \pa_z, 
\qquad 
{\cal T}\inv \pa_y {\cal T} = \pa_z.
\]
Thus, by \eqref{mL3}, 
\[
\mL_4 := {\cal T}\inv \mL_3 {\cal T} 
= \ompath + \muff_3 \pa_{zzz} + e_1(\th,z) \, \pa_z + e_0(\th,z)
\]
where 
\begin{equation}\label{e0 e1}
e_1(\th,z) := \ompath p(\th) + ({\cal T}\inv d_1) (\th,z), 
\quad 
e_0(\th,z) := ({\cal T} \inv d_0)(\th,z).
\end{equation}
Now we look for $p(\th)$ such that the average 
\be \label{condiz fastidio}
\frac{1}{2\p} \, \int_\T e_1(\th,z) \, dz = \muff_1 \, , \quad \forall \vartheta \in \T^\nu \, , 
\ee
for some constant $m_1 \in \R$ (independent of $ \th $). Equation \eqref{condiz fastidio} is equivalent to
\begin{equation}\label{equazione omologica step 4}
\o \cdot \partial_{\th} p  = \muff_1 - \int_{\T} d_{1}(\th , y) \, dy =: V(\th).
\end{equation}
The equation \eqref{equazione omologica step 4} has a periodic solution $p(\th)$ 
if and only if  $\int_{\T^{\nu}}V(\th) \, d \th = 0$. Hence we have to define
\begin{equation}\label{m-p}
\muff_1   :=  \frac{1}{(2\pi)^{\nu+1}} \, \int_{\T^{\nu + 1}} d_{1}(\th , y) \, d \th d y 
\end{equation} 
and 
\be\label{def:p} 
p(\th)  := (\o \cdot \partial_{\th})^{-1}V(\th) \, .
\ee
With this choice of $p$, after renaming the space-time variables $z = x$ and $\th = \ph$, we have 
\begin{equation} \label{mL4}
\mL_4 = \ompaph + \muff_3 \pa_{xxx} + e_1(\ph,x) \, \pa_x + e_0(\ph,x),
\qquad 
\frac{1}{2\p} \, \int_{\T} e_1(\ph,x) \, dx = \muff_1 \, , \ \ \forall \ph \in \T^\nu \, .
\end{equation}

\begin{remark}\label{reversibilitˆ step 4}
By \eqref{equazione omologica step 4}, \eqref{def:p} 
and since $ d_1 \in X $ (see remark \ref{reversibilitˆ step 3}), the function $p$ is odd. 
Then $ {\cal T} $ and $ {\cal T}^{-1}$ defined in \eqref{MM-1} 
are reversibility preserving and the coefficients $ e_1, e_0 $ defined in \eqref{e0 e1} 
satisfy $e_1 \in X$, $e_0 \in Y$. Hence $\mL_4 : X \rightarrow Y$ is reversible.
\end{remark}

\begin{remark}\label{rem: Ham4}
In the Hamiltonian case the operator $ \mL_4 $ is Hamiltonian, because 
the operator $ {\cal T} $ in \eqref{MM-1} is symplectic (it is a particular case
of the change of variables \eqref{operatore1 simplettico} with $ \b (\vphi,x) = p( \vphi)  $). 
\end{remark}

\subsection{Step 5. Descent method: conjugation by pseudo-differential operators}\label{step-5}

The goal of this section is to conjugate $ \mL_4 $ in \eqref{mL4} to an operator of the form 
$ \o \cdot \partial_{\vphi} + \muff_3 \partial_{xxx} + \muff_1 \partial_{x} + \mR $ where 
the constants $\muff_3$, $\muff_1$ are defined in \eqref{mu 3}, \eqref{m-p},  
and $\mR$ is a pseudo-differential operator of order $0$. 

Consider an operator of the form
\begin{equation}\label{mS}
\mS := I + w(\ph,x)  \partial_{x}^{-1} 
\end{equation}
where $w : \T^{\nu + 1}\rightarrow \R$ and the operator $\partial_{x}^{-1}$ is defined in \eqref{dx-1}. 
Note that 
$\partial_x^{-1} \partial_x =  \partial_x \partial_x^{-1} = \pi_0$,  
where $ \pi_0 $ is the $ L^2 $-projector on the subspace 
$ H_0 := \{ u(\ph,x) \in L^2 (\T^{\nu+1})\, : \, \int_{\T} u(\ph, x) \,  dx = 0 \} $.

A direct computation shows that the difference
\begin{equation} \label{carmelino2}
\mL_4 \mS - \mS (\o \cdot \partial_{\vphi} + \muff_3 \partial_{xxx} + \muff_1 \partial_{x}) 
= r_1 \partial_{x} + r_0 + r_{-1} \partial_{x}^{-1} 
\end{equation}
where (using $ \partial_x \pi_0 = \pi_0 \partial_x = \pa_x $, $ \partial_x^{-1} \partial_{xxx} = \partial_{xx} $)
\begin{eqnarray}
r_{1} & := & 3\muff_3 w_{x} + e_{1}(\vphi,x) - \muff_1 \label{r1}\\
r_{0}& := & e_{0} + \big( 3\muff_3 w_{xx} + e_{1}w - \muff_1 w \big)\pi_{0}\label{r0}\\
r_{-1}& := & \o \cdot \partial_{\vphi}w + \muff_3 w_{xxx} + e_{1} w_{x}\,.\label{r-1}
\end{eqnarray}
We look for a periodic function $ w (\vphi, x )$ such that $ r_1 = 0$.  By \eqref{r1} and \eqref{condiz fastidio} we take
\begin{equation}\label{w}
w = \frac{1}{3\muff_3}\partial_{x}^{-1}[\muff_1 - e_{1}].
\end{equation}
For $ \e $ small enough the operator $ {\cal S} $ is invertible and 
we obtain, by \eqref{carmelino2}, 
\be\label{mL5}
\mL_5 := \mS^{-1} \mL_4 \mS 
= \o \cdot \partial_{\vphi} + \muff_3 \partial_{xxx} + \muff_1 \partial_{x} + {\cal R},
\qquad 
{\cal R} := \mS^{-1} ( r_{0} + r_{-1} \partial_{x}^{-1} ).
\end{equation}

\begin{remark}\label{reversibilitˆ step 5}
In the reversible case,  the function $w \in Y$, because $e_1 \in X$, see remark \ref{reversibilitˆ step 4}. 
Then $\mS$, $\mS^{-1}$ are reversibility preserving.  
By \eqref{r0} and \eqref{r-1}, $r_0 \in Y$ and $r_{-1} \in X$. Then 
the operators $ \mR, \mL_5 $ defined in \eqref{mL5} are reversible, namely $\mR, \mL_5 : X \rightarrow Y$.
\end{remark}

\begin{remark}\label{rem:Ham5}
In the Hamiltonian case,  we consider, instead of \eqref{mS}, the modified operator
\be\label{modifiedmS}
\mS  := e^{\pi_0 w(\vphi,x) \partial_{x}^{-1}} := I + \pi_0 w(\vphi,x) \partial_{x}^{-1} + \ldots  
\ee
which, for each $ \vphi \in \T^\nu $, is symplectic. Actually $ \mS $ is the time one flow map of the 
Hamiltonian vector field
$  \pi_0 w(\vphi, x) \partial_{x}^{-1} $ 
which is generated by the Hamiltonian
$$
H_\mS(\vphi, u) := - \frac12\, \int_{\T} w(\vphi, x) \big( \partial_x^{-1} u \big)^2 dx \, \, , \quad u \in H^1_0 \, . 
$$ 
The corresponding $ \mL_5 $ in \eqref{mL5}  is Hamiltonian.  
Note that the operators \eqref{modifiedmS} and \eqref{mS} differ only for 
 pseudo-differential smoothing 
operators of order $ O( \partial_{x}^{-2} ) $ and of smaller size $ O( w^2 ) = O(\e^2) $. 
\end{remark}

\subsection{Estimates on $\mL_5$}  \label{subsec:mL0 mL5}

Summarizing the steps performed in the previous sections \ref{step-1}-\ref{step-5},  we have (semi)-conjugated the operator 
$ \mL  $ defined in \eqref{mL} to the operator $\mL_5 $ defined in \eqref{mL5}, namely  
\begin{equation} \label{Phi 1 2 def}
\mL = \Phi_1 \mL_5 \Phi_2\inv, \qquad 
\Phi_1 := {\cal A} B \rho {\cal M}  {\cal T} \mS, \quad 
\Phi_2 := {\cal A} B {\cal M}  {\cal T} \mS
\end{equation}
(where $ \rho $ means the multiplication operator for the function $ \rho $ defined in \eqref{anche def rho}). 

In the next lemma we give tame estimates for $\mL_5$ and  $\Phi_1, \Phi_2$.
We define the constants
\be\label{costanti lemma mostro}
\s := 2\t_0 + 2 \nu + 17, \quad \s' := 2\t_0 + \nu + 14 
\ee
where $ \t_0 $ is  defined in \eqref{omdio} and $ \nu $ is the number of frequencies.

\begin{lemma}  \label{lemma:mostro} 
Let $ f \in C ^q $, see \eqref{f classe Cq}, and  
 $ \mathfrak s_0 \leq s \leq q - \s $. 
There exists $\d > 0 $ such that, if $  \e \g_0 \inv < \d $ 
(the constant $ \g_0 $ is defined in \eqref{omdio}), then, for all 
\begin{equation}  \label{palla di sicurezza}
\| u \|_{\mathfrak s_0 + \s} \leq 1 \, ,  
\end{equation}
$(i)$ the transformations ${\Phi}_1, {\Phi}_2$ defined in \eqref{Phi 1 2 def} are invertible operators of $H^s(\T^{\nu+1})$, and satisfy 
\be \label{stima Phi 12 nel lemma} 
\| \Phi_i h \|_s + \| \Phi_i\inv h \|_s 
\leq C(s) \big( \| h \|_s + \| u \|_{s+\s} \| h \|_{\mathfrak s_0} \big),  
\ee
for $i = 1, 2$.
Moreover, if $u(\lm)$, $h(\lm)$ are Lipschitz families with 
\begin{equation}  \label{palla Lip di sicurezza}
\| u \|_{\mathfrak s_0 + \s}^{\Lipg} \leq 1, 
\end{equation}
then 
\begin{equation}
\label{stima Lip Phi 12 nel lemma}
\| \Phi_i h \|_s^{\Lipg} + \| \Phi_i\inv h \|_s^{\Lipg} 
\leq C(s) \big( \| h \|_{s+3}^{\Lipg} + \| u \|_{s+\s}^{\Lipg} \| h \|_{\mathfrak s_0+3}^{\Lipg} \big), 
\quad i = 1,2.
\end{equation}
$(ii)$ 
The constant coefficients $\muff_3, \muff_1$ of $\mL_5$ defined in \eqref{mL5} satisfy
\begin{align} \label{coefficienti costanti 1} 
| \muff_3 - 1| + |\muff_1| 
& \leq \e C \, ,
\\
\label{coefficienti costanti 2} 
| \pa_u \muff_3(u)[h]| + | \pa_u \muff_1(u)[h]| 
& \leq \e C \| h \|_{\s} \, . 
\end{align}
Moreover, if $u(\lm)$ is a  Lipschitz family  satisfying 
\eqref{palla Lip di sicurezza}, 
then 
\begin{equation} \label{Lete 3} 
| \muff_3 - 1 |^{\Lipg} + | \muff_1 |^{\Lipg} 
\leq \e C .
\end{equation}
$(iii)$  
The operator $\mR$ defined in \eqref{mL5} satisfies:
\begin{align} \label{stima R 1} 
| \mR |_s 
& \leq \e C(s) (1 + \| u \|_{s + \s}),  
\\ 
\label{stima R 2} 
| \pa_u \mR(u)[h] \, |_{s} 
& \leq \e C(s) \big(  \| h \|_{s + \s'} 
  + \| u \|_{s + \s} \| h \|_{\mathfrak s_0 + \s'} \big)  \,  ,
\end{align}
where $ \s > \s' $ are defined in \eqref{costanti lemma mostro}.
Moreover, if $u(\lm)$ is a Lipschitz family satisfying \eqref{palla Lip di sicurezza}, 
then 
\begin{equation} \label{stima R 3} 
| \mR |_s^{\Lipg} 
\leq \e C(s) (1 + \| u \|_{s + \s}^{\Lipg}).  
\end{equation} 
Finally, in the reversible case,  the maps $\Phi_i, \Phi_i^{-1}$, $i=1,2$ are reversibility preserving and $\mR, \mL_5 : X \rightarrow Y$
are reversible. In the Hamiltonian case the operator $ \mL_5 $ is Hamiltonian. 
\end{lemma}

\begin{pf} 
In section \ref{sec:proofs}. 
\end{pf}

\begin{lemma}  \label{lemma:stime stabilita Phi 12}
In the same hypotheses of Lemma \ref{lemma:mostro}, 
for all $\ph \in \T^\nu$, the operators 
${\cal A}(\ph)$, $ {\cal M} (\ph)$, $ {\cal T} (\ph)$, $\mS(\ph)$ 
are invertible operators of the phase space $H^s_x := H^s(\T)$, 
with 
\begin{align} 
\label{A(ph)}
\| {\cal A}^{\pm 1}(\ph) h \|_{H^s_x} 
& \leq 
C(s) \big( \| h \|_{H^s_x} + \| u \|_{s + \mathfrak s_0 + 3} \| h \|_{H^1_x} \big),
\\
\label{A(ph)-I}
\| ({\cal A}^{\pm 1}(\ph) - I) h \|_{H^s_x} 
& \leq 
\e C(s) \big( \| h \|_{H^{s+1}_x} + \| u \|_{s + \mathfrak s_0 + 3} \| h \|_{H^2_x} \big), 
\\
\label{Phi mM mS (ph)}
\| ({\cal M} (\ph) {\cal T}(\ph) \mS(\ph))^{\pm 1} h \|_{H^s_x} 
& \leq C(s) \big( \| h \|_{H^s_x} 
+ \| u \|_{s + \s} \| h \|_{H^1_x} \big),
\\
\label{Phi mM mS (ph) - I}
\| (( {\cal M} (\ph) {\cal T}(\ph) \mS(\ph) )^{\pm 1} - I) h \|_{H^s_x} 
& \leq 
\e \g_0^{-1} C(s) 
\big( \| h \|_{H^{s+1}_x} 
+ \| u \|_{s + \s} \| h \|_{H^1_x} \big).
\end{align}
\end{lemma}

\begin{pf} 
In section \ref{sec:proofs}. 
\end{pf}

\section{Reduction of the linearized operator to constant coefficients}\label{sec:redu}

The goal of this section is to diagonalize the linear operator $ \mL_5 $  obtained in \eqref{mL5}, and therefore to complete the reduction of 
$ {\cal L}  $ in \eqref{mL} into constant coefficients.  For $ \t > \t_0 $ (see \eqref{omdio}) we define the constant
\be \label{defbq}
\beta := 7 \tau + 6 \, . 
\ee

\begin{theorem}\label{teoremadiriducibilita}
Let $ f \in C^q $,  see \eqref{f classe Cq}.
Let $ \g \in (0,1) $ and 
$ \mathfrak s_0 \leq s \leq q - \s - \b $
where $ \s $ is defined in \eqref{costanti lemma mostro}, and $ \b $ in \eqref{defbq}. 
Let  $u(\lm) $ be a family of functions depending on the parameter $\lm \in 
\Lm_o \subset \Lm := [1/2, 3/2]$ in a Lipschitz way, with 
\be\label{norma bassa u riducibilitˆ}
\Vert u \Vert_{\mathfrak s_0 + \s + \b, \L_o}^{\Lipg} \leq 1.
\ee
Then there exist $ \delta_{0} $,  $ C $ (depending on the data of the problem) such that, if
\begin{equation}\label{condizione-kam}
\e \gamma^{-1}   \leq \delta_{0} \, ,
\end{equation}
then:
\\[1mm]
$(i)$ {\bf (Eigenvalues)}
$\forall \lm \in \Lm $ there exists a sequence 
\be\label{espressione autovalori}
\mu_j^\infty(\lm) := \mu_j^\infty(\lm, u) = {\tilde \mu}^{0}_j(\lm) + r_j^\infty(\lm) \, , \ {\tilde \mu}^0_j(\lm) := 
 \ii \big( - {\tilde m}_3 ( \lm) j^3 + {\tilde m}_1(\lm)  j \big) \, , \ j \in \Z \, , 
\ee
where $ {\tilde m}_3, {\tilde m}_1$ coincide with the coefficients of $ {\cal L}_5 $ in \eqref{mL5} for all $ \l \in  \Lambda_o $,
and the corrections $r_j^\infty$ satisfy
\begin{align}
\label{autofinali}
| {\tilde \muff}_3 - 1 |^{\Lipg} + | {\tilde \muff}_1 |^{\Lipg} + | r^{\infty}_j |^{\Lipg}_{\Lm}
& \leq \e C \, , \ \ \forall j \in \Z \, .
\end{align}
Moreover, in the reversible case (i.e. \eqref{parity f} holds) or Hamiltonian case
(i.e. \eqref{f Ham} holds), all the eigenvalues $\mu_j^{\infty}$ are purely imaginary.

$(ii)$ {\bf (Conjugacy)}.
For all $\lm$ in 
\begin{equation}  \label{Omegainfty}
\Lambda_\infty^{2\g} := \Lambda_\infty^{2\g} (u) := \Big\{ \lm \in \L_o \, : \,
| \ii \lm \bar\om \cdot l + \mu^{\infty}_j (\lm) - \mu^{\infty}_{k} (\lm) |
\geq 2 \gamma | j^{3} - k^{3} | \langle l \rangle^{-\tau}, 
\ \forall l \in \Z^{\nu}, \, j ,k \in \Z \Big\} 
\end{equation}
there is a bounded, invertible linear operator $\Phi_\infty(\lm) : H^s \to H^s$, with bounded inverse 
$\Phi_\infty\inv(\lm)$, that conjugates $\mL_5$ in \eqref{mL5} to constant coefficients, namely
\be\label{Lfinale}
{\cal L}_{\infty}(\lm)
:= \Phi_{\infty}\inv(\lm) \circ \mL_5(\lm) \circ  \Phi_{\infty}(\lm)
= \lm \bar \om \cdot \partial_{\vphi} + {\cal D}_{\infty}(\lm),  \quad
{\cal D}_{\infty}(\lm)
:= {\rm diag}_{j \in \Z}  \mu^{\infty}_{j}(\lm) \, .
\ee
The transformations $\Phi_\infty, \Phi_\infty\inv$ are close to the identity in matrix decay norm,
with estimates
\begin{equation} \label{stima Phi infty}
| \Phi_{\infty} (\lm) - I |_{s,\Lm_\infty^{2\g}}^\Lipg
+ | \Phi_{\infty}^{- 1} (\lm) - I |_{s,\Lm_\infty^{2\g}}^\Lipg 
\leq \e \gamma^{-1} C(s) \big( 1 + \| u \|_{s + \s + \b ,\Lm_o }^\Lipg \big).
\end{equation}
For all $\ph \in \T^\nu$, the operator 
$\Phi_\infty(\ph) : H^s_x \to H^s_x $  
is invertible  
(where $H^s_x := H^s(\T)$) with inverse $ (\Phi_\infty(\ph))^{-1} = \Phi_\infty^{-1}(\ph)$, 
and 
\begin{align} 
\label{Phi infty (ph) - I}
\| (\Phi_\infty^{\pm 1}(\ph) - I) h \|_{H^s_x} 
& \leq 
\e \g^{-1} C(s) \big( \| h \|_{H^s_x} + \| u \|_{s + \s + \b + \mathfrak s_0} \| h \|_{H^1_x} \big). 
\end{align}
In the reversible case $\Phi_{\infty}, \Phi_{\infty}^{-1} : X  \rightarrow X $, $Y  \rightarrow Y$ are reversibility preserving, and 
$\mL_\infty : X \rightarrow Y$ is reversible. In the Hamiltonian case the final  $\mL_\infty $ is Hamiltonian.
\end{theorem}

An important point of  Theorem \ref{teoremadiriducibilita} is to require {\it only} the bound 
\eqref{norma bassa u riducibilitˆ} for  the low norm of $ u $,   
but it provides the estimate for $ \Phi_\infty^{\pm 1} - I $ in \eqref{stima Phi infty} also for the higher norms $ | \cdot |_s $, 
depending also on the high norms of $ u $.  
From Theorem \ref{teoremadiriducibilita}  we shall deduce 
tame estimates for the inverse linearized operators in
Theorem \ref{inversione linearizzato}. 

\smallskip

Note also that the set $ \L_{\infty}^{2 \g} $ in \eqref{Omegainfty} depends only 
of the final eigenvalues, and it is not defined inductively as in usual KAM theorems. 
This  characterization of the set of parameters which fulfill  
all the required Melnikov non-resonance conditions (at any step of the iteration) 
was first observed in \cite{BB4}, \cite{BB10}
in an analytic setting. Theorem \ref{teoremadiriducibilita} extends this property also in a differentiable setting. 
A main  advantage of this formulation is that it allows to discuss the measure estimates only 
once and not inductively: the Cantor set  $  \L_{\infty}^{2 \g} $  in \eqref{Omegainfty}  could be 
empty (actually its measure $ |\L_{\infty}^{2 \g} | = 1  - O(\g) $ as $ \g \to 0 $) 
but the functions $  \mu^\infty_j (\l) $ are anyway well defined for all $ \l \in \L $, see \eqref{espressione autovalori}. 
In particular we shall perform the measure estimates only along the nonlinear iteration, see section \ref{sec:NM}. 

\smallskip

Theorem \ref{teoremadiriducibilita} is deduced from  
the following iterative Nash-Moser reducibility theorem for a linear operator of the form 
\be\label{L0}
{\cal L}_{0} = \o \cdot \partial_{\vphi} + {\cal D}_{0} + {\cal R}_{0} \, , 
\ee
where $\om = \lm \bar\om $,    
\be\label{D0R0} 
{\cal D}_{0} := m_3 (\l,u(\l)) \partial_{xxx} + m_1(\l,u(\l)) \partial_{x} \, , \quad  \mR_0(\l,u(\l)) :=  \mR(\l,u(\l))  \, , 
\ee
the  $ m_3(\l,u(\l)), m_1 (\l,u(\l))  \in \R $ and 
$ u(\l) $ is defined for $ \l \in \L_o \subset \L $. 
Clearly $\mL_5$ in \eqref{mL5} has the form \eqref{L0}. 
Define
\be\label{defN}
N_{-1} := 1 \, ,  \quad 
N_{\nu} := N_{0}^{\chi^{\nu}} \  \forall \nu \geq 0 \, ,  \quad  
\chi := 3 /2  
\ee
(then $ N_{\nu+1} = N_{\nu}^\chi $, $ \forall \nu \geq 0 $) and 
 \begin{equation}\label{alpha-beta}
 \alpha := 7 \tau + 4, \quad \s_2 := \s + \b  
 \end{equation}
 where $\s$ is defined in \eqref{costanti lemma mostro} and $\b$ is defined in \eqref{defbq}.

\begin{theorem}{{\bf (KAM reducibility)}} \label{thm:abstract linear reducibility}  
Let $ q >  \s +  \mathfrak s_0 + \b $. 
There exist $C_0 > 0 $,   $ N_{0}  \in \N $ large,  such that, if 
\begin{equation}\label{piccolezza1}
N_{0}^{C_0}  |{\cal R}_{0} |_{\mathfrak s_{0} + \beta}^{\Lipg} \g^{-1} \leq 1,
\end{equation}
then, for all $ \nu \geq 0 $:

\begin{itemize}
\item[${\bf(S1)_{\nu}}$] 
There exists an operator 
\be\label{def:Lj}
{\cal L}_\nu := \o \cdot \partial_{\vphi} + {\cal D}_\nu + {\cal R}_\nu \quad 
where \quad {\cal D}_\nu = {\rm diag}_{j \in \Z} \{ \mu^{\nu}_{j}(\l ) \}
\ee
\be\label{mu-j-nu}
\mu_j^{\nu}(\l) 
= \mu_j^0(\l) + r_j^{\nu}(\l),\quad 
\mu_j^0(\l) := - \ii \big( m_3(\l,u(\l)) j^3 - m_1(\l,u(\l)) j \big) , \ \  j \in \Z \, , 
\ee
defined for all $ \l \in \L_{\nu}^{\g}(u)$, where 
$\L_{0}^{\g}(u) := \L_o $ (is the domain  of $ u $), and, for $\nu \geq 1$, 
\be\label{Omgj}
\L_{\nu}^{\g} := 
\L_{\nu}^{\gamma}(u):= \Big\{\l \in \L_{\nu - 1}^{\gamma} : \left| \ii \o \cdot l + \mu^{\nu-1}_{j}(\l) -
\mu^{\nu - 1}_{k}(\l) \right| \geq \gamma \frac{ |j^{3} - k^{3}|}{\left\langle l\right\rangle^{\tau}}
\  \forall \left|l\right| \leq N_{ \nu-1}, \ j, k \in \Z \Big\}.
\ee
For $\nu \geq 0$, $r_j^{\nu} = \overline{r_{-j}^{\nu}}$, equivalently  $ \mu_j^{\nu} = \overline{\mu_{-j}^{\nu}}$, and 
\begin{equation}  \label{rjnu bounded}
|r_j^{\nu}|^{\Lipg} := |r_j^{\nu}|^{\Lipg}_{\L_{\nu}^\g} \leq \e C \, . 
\end{equation}
The remainder $  {\cal R}_\nu  $ is real (Definition \ref{def:RR}) and, $ \forall s \in [ \mathfrak s_0, q - \s - \b ] $,  
\be\label{Rsb} 
\left|{\cal R}_{\nu}\right|_{s}^{\Lipg} \leq  \left|{\cal R}_{0}\right|_{s+\beta}^{\Lipg}N_{\nu - 1}^{-\alpha} \, , \quad 
\left|{\cal R}_{\nu}\right|_{s + \beta}^{\Lipg} \leq \left|{\cal R}_{0}\right|_{s+\beta}^{\Lipg}\,N_{\nu - 1}\,.
\ee
Moreover, for $ \nu \geq 1 $, 
\be	\label{Lnu+1}
{\cal L}_{\nu} = \Phi_{\nu-1}^{-1} {\cal L}_{\nu-1} \Phi_{\nu-1} \, , \quad \Phi_{\nu-1} := I + \Psi_{\nu-1} \, , 
\ee
where the map $ \Psi_{\nu-1} $  is real,  T\"oplitz in time $ \Psi_{\nu-1} := \Psi_{\nu-1}(\vphi) $ 
(see \eqref{Aphi}), 
 and satisfies 
\be\label{Psinus} 
\left|\Psi_{\nu-1} \right|_{s}^{\Lipg} \leq 
|{\cal R}_{0} |_{s+\beta}^{\Lipg} \gamma^{-1} N_{\nu-1}^{2 \t+1} N_{\nu - 2}^{- \a}  \, . 
\ee
In the reversible case, $\mR_{\nu} : X \rightarrow Y$, $\Psi_{\nu - 1}, \Phi_{\nu - 1}, \Phi_{\nu - 1}^{-1} $ are reversibility preserving. 
Moreover, all the $ \mu^\nu_{j}(\l) $ are purely imaginary  
and $ \mu^{\nu}_{j} = - \mu^{\nu}_{-j} $, $ \forall j \in \Z $.

\item[${\bf(S2)_{\nu}}$] 
For all $ j \in \Z $,  there exist Lipschitz extensions 
$ \widetilde{\mu}_{j}^{\nu}(\cdot ): \L \to \R $ of  $ \mu_{j}^{\nu}(\cdot ) : \L_\nu^\g \to \R $ 
satisfying,  for $\nu \geq 1$, 
\be\label{lambdaestesi}  
|\widetilde{\mu}_{j}^{\nu} -  \widetilde{\mu}_{j}^{\nu-1} |^{\Lipg}  \leq | {\cal R}_{\nu-1} |^{\Lipg}_{\mathfrak s_0} \,.
\ee
\item[ ${\bf (S3)_{\nu}}$] 
Let $ u_1(\l)$, $ u_2(\l)$, be Lipschitz families of Sobolev functions, 
defined for $\l \in \L_o$ and such that conditions  \eqref{norma bassa u riducibilitˆ},
\eqref{piccolezza1} hold
with $ {\cal R}_0 := {\cal R}_0 ( u_i) $, $ i = 1,2  $, see \eqref{D0R0}. 

Then, for $ \nu \geq 0 $,   $\forall \l \in \L_{\nu}^{\g_1}(u_1) \cap  \L_{\nu}^{\g_2}(u_2)$, 
with $ \g_1, \g_2 \in [\g/2, 2\g]$, 
\begin{equation}\label{derivate-R-nu}
|{\cal R}_{\nu}(u_2) - \mR_{\nu}(u_1)|_{\mathfrak s_{0}}\leq\e N_{\nu - 1}^{-\alpha} \Vert u_1 - u_2 \Vert_{\mathfrak s_0 + \s_2},\,\,|{\cal R}_{\nu}(u_2) - \mR_{\nu}(u_1)|_{\mathfrak s_{0}+\beta} \leq \e  N_{\nu - 1} \Vert u_1 - u_2 \Vert_{\mathfrak s_0 + \s_2} \, .
\end{equation}
Moreover, for $\nu \geq 1$, $ \forall s \in [\mathfrak s_{0},\mathfrak s_{0}+\b] $, 
$\forall j \in\Z $, 
\be\label{deltarj12}
\big|\big(r_{j}^{\nu}(u_2) - r_{j}^{\nu}(u_1)\big) - \big(r_{j}^{\nu-1}(u_2) - r_{j}^{\nu-1}(u_1)\big) \big|\leq  \vert {\cal R}_{\nu -1}(u_2) - \mR_{\nu -1}(u_1) \vert_{\mathfrak s_0} \,, 
\ee
\be\label{Delta12 rj}
| r_j^{\nu}(u_2) - r_j^{\nu}(u_1) | \leq  \e C \|  u_1 - u_2 \|_{\mathfrak s_0 + \s_2} \, . 
\ee 

\item[ ${\bf (S4)_{\nu}}$]
Let $u_1, u_2$ like in $({\bf S3})_\nu $ and $ 0 < \rho < \g / 2 $. For all $ \nu \geq 0 $ such that  
\begin{equation}\label{legno}
\e C N_{\nu - 1}^{\tau} \Vert u_1 - u_2 \Vert_{\mathfrak s_0 + \s_2}^{\rm sup} \leq \rho 
\quad \Longrightarrow \quad
 \Lambda_{\nu }^{\gamma}(u_1) \subseteq 
 \Lambda_{\nu}^{\gamma - \rho}(u_2) \, . 
\end{equation}
\end{itemize}
\end{theorem}

\begin{remark}\label{rem:Ham6}
In the Hamiltonian case $ \Psi_{\nu-1}$ is Hamiltonian and, 
instead of \eqref{Lnu+1} we consider the symplectic map 
\be\label{modified Phi nu}
\Phi_{\nu-1}  := \exp(\Psi_{\nu-1}) \, . 
\ee 
The corresponding operators $ \mL_\nu $, ${\cal R}_\nu $  are Hamiltonian.  
Note that the operators \eqref{modified Phi nu} and \eqref{Lnu+1} differ for an operator of order 
$\Psi_{\nu - 1}^2$. 
\end{remark}

The proof of Theorem \ref{thm:abstract linear reducibility} is postponed in Subsection 
\ref{subsec:proof of thm abstract linear reducibility}. We first give some consequences.

\begin{corollary}\label{lem:convPhi}
{\bf (KAM transformation)} 
$ \forall \l \in  \cap_{\nu \geq 0} \L_{\nu}^{\g} $ 
the sequence 
\be\label{Phicompo}
\widetilde{\Phi}_{\nu} := \Phi_{0} \circ \Phi_1 \circ \cdots\circ \Phi_{\nu} 
\ee
converges in $ |\cdot |_{s}^{\Lipg}$ to an operator $\Phi_{\infty} $
and 
\be\label{Phinftys}
\left|\Phi_{\infty} - I \right|_{s}^{\Lipg} + \left|\Phi_{\infty}^{-1} - I \right|_{s}^{\Lipg} \leq 
C(s) \left|{\cal R}_{0}\right|_{s + \beta}^{\Lipg} \gamma^{-1} \, . 
\ee
In the reversible case $\Phi_\infty$ and $\Phi_{\infty}^{-1}$ are reversibility preserving.
\end{corollary}

\begin{pf}
To simplify notations we write $|\cdot|_s $ for $|\cdot|_s^{\Lipg}$. 
For all $ \nu \geq 0 $ we have
$  \widetilde{\Phi}_{\nu + 1} = \widetilde{\Phi}_{\nu}\circ \Phi_{\nu + 1} = 
 \widetilde{\Phi}_{\nu} + \widetilde{\Phi}_{\nu}\Psi_{\nu + 1}  $ (see \eqref{Lnu+1}) and so
 \be\label{lasopra}
 |\widetilde{\Phi}_{\nu + 1}|_{\mathfrak s_{0}} \stackrel{\eqref{algebra Lip}} \leq 
 |\widetilde{\Phi}_{\nu} |_{\mathfrak s_{0}} + C 
 |\widetilde{\Phi}_{\nu} |_{\mathfrak s_{0}}\left|\Psi_{\nu + 1}\right|_{\mathfrak s_{0}}
 \stackrel{\eqref{Psinus}} \leq  |\widetilde{\Phi}_{\nu} |_{\mathfrak s_{0}} ( 1+  \e_\nu )  
 \ee
where $ \e_\nu := C' |{\cal R}_{0} |_{\mathfrak s_0 +\beta}^{\Lipg} \gamma^{-1} N_{\nu+1}^{2 \t+1} N_{\nu}^{- \a} $.  
Iterating \eqref{lasopra} we get, for all $ \nu $,  
\be\label{bassa}
|\widetilde{\Phi}_{\nu + 1}|_{\mathfrak s_{0}} \leq  | \widetilde{\Phi}_0  |_{\mathfrak s_{0}}   
\Pi_{\nu \geq 0} (1+ \e_\nu) \leq  | \Phi_0 |_{\mathfrak s_0} e^{C |{\cal R}_{0} |_{\mathfrak s_0+\beta}^{\Lipg} \gamma^{-1}} 
\leq 2  
\ee
using \eqref{Psinus} (with $ \nu =1 $, $ s = \mathfrak s_0  $) to estimate $ | \Phi_0 |_{\mathfrak s_0} $ and \eqref{piccolezza1}.
The high norm of $  \widetilde{\Phi}_{\nu + 1} = \widetilde{\Phi}_{\nu} + \widetilde{\Phi}_{\nu}\Psi_{\nu + 1}  $ 
is estimated 
by \eqref{interpm Lip}, \eqref{bassa} (for $ {\wtilde \Phi}_\nu $), as 
\begin{eqnarray*} 
 |\widetilde{\Phi}_{\nu + 1}|_s 
 & \leq & | \widetilde{\Phi}_{\nu} |_s ( 1 + C(s)\left|\Psi_{\nu + 1}\right|_{\mathfrak s_{0}}  )  
+ C(s)  \left|\Psi_{\nu + 1}\right|_s \nonumber  \\
& \stackrel{\eqref{Psinus}, \eqref{alpha-beta}}\leq &| 
\widetilde{\Phi}_{\nu} |_s ( 1 + \e_{\nu}^{(0)}) +  \e_{\nu}^{(s)} \, , 
\ \e_{\nu}^{(0)} := |{\cal R}_0|_{\mathfrak s_0+\b} \g^{-1} N_\nu^{-1} \, , 
\ \e_{\nu}^{(s)} := |{\cal R}_0|_{s+\b} \g^{-1} N_\nu^{-1} \, .
\end{eqnarray*}
Iterating the above inequality and, using $ \Pi_{j \geq 0} (1+ \e_j^{(0)}) \leq 2 $, 
we get
\be\label{trieste}
 |\widetilde{\Phi}_{\nu + 1}|_s \leq_s \sum_{j=0}^\infty \e_j^{(s)} +  
   |\widetilde{\Phi}_0 |_s \leq C(s) \big(1+ | {\cal R}_0|_{s+ \b} \g^{-1} \big)  
\ee
using $   |\Phi_0 |_s \leq 1+ C(s) | {\cal R}_0|_{s+ \b} \g^{-1}  $. 
Finally, the $\widetilde{\Phi}_{j}$   a Cauchy sequence in  norm $ | \cdot |_s $ because 
\begin{eqnarray}
| \widetilde{\Phi}_{\nu+m} - \widetilde{\Phi}_{\nu}  |_{s} \! \!\! \! \!\! 
& \! \!\! \leq \! \!\! & \! \!\! \! \!\! 
\sum_{j =\nu}^{\nu+m-1} |\widetilde{\Phi}_{j + 1} - \widetilde{\Phi}_{j} |_{s}
 \stackrel{ \eqref{interpm Lip}} 
 {\leq_s} \sum_{j = \nu}^{\nu+m-1}\left( |\widetilde{\Phi}_j |_s |\Psi_{j + 1} |_{\mathfrak s_{0}}
 + |\widetilde{\Phi}_j |_{\mathfrak s_{0}} |\Psi_{j + 1} |_{s}\right) \nonumber \\
& \! \!\! \! \! \stackrel{\eqref{trieste},  \eqref{Psinus}, \eqref{bassa}, \eqref{piccolezza1}} {\leq_s} \! \!\! \! \! &  \sum_{j \geq \nu}  \left|{\cal R}_{0}\right|_{s + \beta}  \gamma^{-1}  N_j^{-1} \leq_s \left|{\cal R}_{0}\right|_{s + \beta}  \gamma^{-1} N_\nu^{-1} \, . \label{quasifi} 
\end{eqnarray}
Hence $ \widetilde{\Phi}_{\nu} \stackrel{\left|\cdot\right|_ s}{\rightarrow} \Phi_{\infty} $. 
The  bound for $ \Phi_\infty -  I $ in \eqref{Phinftys} follows  
by \eqref{quasifi} with $ m = \infty $, $ \nu = 0 $ and 
$ |\wtilde \Phi_0 - I |_s = $ $ |\Psi_0|_s  \lessdot  \g^{-1} | {\cal R}_0|_{s+\b} $. 
Then the estimate for $ \Phi_\infty^{-1} -  I $ follows by \eqref{PhINV}. 

In the reversible case all the $\Phi_\nu $ are reversibility preserving and so
$\widetilde{\Phi}_\nu $, $\Phi_{\infty}$ are reversibility preserving.
\end{pf}

\begin{remark}\label{rem:Ham7}
In the Hamiltonian case, the transformation $\widetilde{\Phi}_\nu$ in \eqref{Phicompo} 
is symplectic, because $\Phi_\nu$ is symplectic for all $\nu$ (see Remark \ref{rem:Ham6}). 
Therefore $\Phi_\infty$ is also symplectic. 
\end{remark}

Let us define for all $j \in \Z$
$$
\mu^{\infty}_{j}(\l) = \lim_{\nu \to + \infty} \widetilde{\mu}_{j}^{\nu}(\l) = \tilde{\mu}_j^0 + r_j^{\infty}(\l), \quad r_j^{\infty}(\l) := \lim_{\nu \to + \infty} \tilde{r}_j^{\nu}(\l) \quad \forall \l \in \L. 
$$
It could happen that $ \L_{\nu_0}^\g = \emptyset $ (see \eqref{Omgj}) for some $ \nu_0 $. In such a case
the iterative process of Theorem \ref{thm:abstract linear reducibility} stops  after finitely many steps. 
However, we can always set $ \widetilde{\mu}_{j}^{\nu}  := \widetilde{\mu}_{j}^{\nu_0} $, $ \forall \nu \geq \nu_0 $,
and the functions $ \mu^{\infty}_{j} : \L \to \R $  are always well defined.

\begin{corollary} {\bf (Final eigenvalues)} For all $ \nu \in \N $, $ j \in \Z $
\be\label{autovcon}
| { \mu }_{j}^{\infty} - {\widetilde \mu }^{\nu}_{j} |_\L^{\Lipg} = 
| r_{j}^{\infty} - {\widetilde r }^{\nu}_{j} |^{\Lipg}_\L 
\leq 
C \left|{\cal R}_{0}\right|_{\mathfrak s_{0}+\beta}^{\Lipg} N_{\nu-1}^{-\alpha}  \, , \ \ 
| { \mu }_{j}^{\infty} - {\widetilde \mu }^{0}_{j}|_\L^{\Lipg} = | r_j^{\infty} |_\L^{\Lipg}
 \leq C \left|{\cal R}_{0}\right|_{\mathfrak s_{0}+\beta}^{\Lipg} \, . 
\ee
\end{corollary}

\begin{pf}
The bound \eqref{autovcon} follows by \eqref{lambdaestesi} and \eqref{Rsb} by summing the telescopic series. 
\end{pf}

\begin{lemma} {\bf (Cantor set)}
\be\label{cantorinclu}
\L_{\infty}^{2 \g} \subset \cap_{\nu \geq 0} \L_{\nu}^\g \, . 
\ee 
\end{lemma}

\begin{pf} 
Let $ \lambda \in \Lambda_{\infty}^{2\g} $. By definition $ \Lambda_{\infty}^{2\g} \subset  \L_0^\g :=  \L_o $. Then 
for all $ \nu > 0 $, $ | l | \leq N_{\nu} $, $ j \neq k $ 
\begin{eqnarray*}
\left| \ii \o \cdot l + { \mu}_{j}^{\nu} - { \mu}_k^{\nu}\right| & \geq & 
\left| \ii \o \cdot l + \mu_j^{\infty} - 
\mu_k^{\infty}\right| - \left| 
{ \mu}_j^{\nu} -   \mu_j^{\infty}\right| - 
\left| {\mu}_k^{\nu} - \mu_k^{\infty}\right| \\
& \stackrel{\eqref{Omegainfty}, \eqref{autovcon}} \geq & 
2\gamma \left |j^{3} - k^{3}\right|\left\langle l\right\rangle^{-\tau} - 2 C | {\cal R}_0|_{\mathfrak s_0+ \b} 
N_{\nu-1}^{-\alpha} \geq \gamma \left|j^{3} - k^{3}\right|\left\langle l\right\rangle^{-\tau}
\end{eqnarray*}
because
$ \gamma |j^{3} - k^{3} | \langle l \rangle^{-\tau} \geq 
\gamma N_\nu^{-\tau} \stackrel{\eqref{piccolezza1}} \geq  2 C | {\cal R}_0|_{\mathfrak s_0+ \b} N_{\nu-1}^{-\alpha} $. 
\end{pf}

\begin{lemma}\label{lemma42} 
For all $\l \in \L_{\infty}^{2\g} (u) $ , 
\be\label{realtˆ autovalori finali}
\mu_j^{\infty}(\l) = \overline{\mu_{-j}^{\infty}(\l)}, \quad r_j^{\infty}(\l) = \overline{r_{-j}^{\infty}(\l)}\,, 
\ee
and in the reversible case 
\be\label{reversibilitˆ autovalori finali}
\mu_j^{\infty}(\l) = - \mu_{-j}^{\infty}(\l), \quad r_j^{\infty}(\l) = - r_{-j}^{\infty}(\l)\,.
\ee
Actually in the reversible case $ \mu_j^{\infty}(\l)  $ are purely imaginary for all $ \l \in \L $.
\end{lemma}

\begin{pf}
Formula \eqref{realtˆ autovalori finali} and \eqref{reversibilitˆ autovalori finali} follow because,  
for all  $ \l \in \L_\infty^{2\g} \subseteq \cap_{\nu\geq 0} \L_\nu^\g $ (see \eqref{cantorinclu}), we have
$ \mu_j^{\nu} = \overline{\mu_{-j}^{\nu}} $, $ r_j^{\nu} = \overline{r_{-j}^{\nu}} $, and, in the reversible case, 
the $ \mu_j^{\nu} $ are purely imaginary and  $ \mu_j^{\nu} = - \mu_{-j}^{\nu} $, $ r_j^{\nu} = - r_{-j}^{\nu} $. 
The final statement follows because,  in the reversible case, the $ \mu_j^\nu (\l) \in \ii \R $ 
as well as its extension $ {\wtilde \mu}_j^\nu (\l)  $. 
\end{pf}

\begin{remark}\label{r0=0}
In the reversible case,  \eqref{reversibilitˆ autovalori finali}  imply that 
$\mu_0^\infty = r_0^\infty = 0$.
\end{remark}

\noindent
{\bf Proof of Theorem \ref{teoremadiriducibilita}.}  
We apply Theorem \ref{thm:abstract linear reducibility}  to the linear operator  
$ {\cal L}_0 := {\cal L}_5$  in  \eqref{mL5}, where 
 $ \mR_0 =  {\cal R }$ defined in \eqref{D0R0} satisfies 
 \be\label{R0tame}
\left|{\cal R}_{0}\right|_{\mathfrak s_{0} + \beta}^{\Lipg} \stackrel{\eqref{stima R 3}}\leq 
\e C(\mathfrak s_0 + \b) \Big(1 + \| u \|_{\mathfrak s_{0} + \s + \b}^{\Lipg}\Big) \stackrel{\eqref{norma bassa u riducibilitˆ}}\leq 
2 \e C(\mathfrak s_0 + \b) \,.
\ee
Then the smallness condition (\ref{piccolezza1}) is implied by \eqref{condizione-kam} 
taking $\delta_0:= \delta_0(\nu)$ small enough.
  
For all $ \l \in \L_\infty^{2\g} \subset 
\cap_{\nu \geq 0} \L_\nu^\g $ (see \eqref{cantorinclu}), 
the operators 
\be\label{conjugnu}
{\cal L}_{\nu} 
\stackrel{\eqref{def:Lj}} =  \o \cdot \partial_{\vphi} + {\cal D}_{\nu} + {\cal R}_{\nu} 
\stackrel{\left|\cdot\right|_{s}^{\Lipg}}{\longrightarrow}\o\cdot \partial_{\vphi} + {\cal D}_{\infty} =: {\cal L}_{\infty} \, ,
\quad
{\cal D}_{\infty} := {\rm diag}_{j \in \Z} \mu_j^{\infty} 
\ee
because 
$$
\left|{\cal D}_{\nu} - {\cal D}_{\infty}\right|_{s}^{\Lipg}  =  \sup_{j \in \Z} 
\left|{\mu}_{j}^{\nu} - \mu_{j}^{\infty}\right|^{\Lipg} \stackrel{\eqref{autovcon}}\leq 
C \left|{\cal R}_{0}\right|_{\mathfrak s_{0}+\beta}^{\Lipg} N_{\nu - 1}^{- \alpha}, \quad 
\left|{\cal R}_{\nu}\right|_{s}^{\Lipg} \stackrel{\eqref{Rsb}} \leq 
\left|{\cal R}_{0}\right|_{s + \beta}^{\Lipg} N_{\nu - 1}^{-\alpha} \,  .
$$
Applying  \eqref{Lnu+1} iteratively we get 
$ {\cal L}_{\nu}  = {{\widetilde \Phi}_{\nu-1}}^{-1} {\cal L}_0 {\widetilde \Phi}_{\nu-1} $ where 
$ {\widetilde \Phi}_{\nu-1} $ is defined by \eqref{Phicompo} 
and 
$  {\widetilde \Phi}_{\nu-1}  \to  {\Phi}_\infty $ in $ | \ |_s $ (Corollary \ref{lem:convPhi}).
Passing to the limit 
we deduce  \eqref{Lfinale}. 
Moreover \eqref{autovcon} and \eqref{R0tame}  
imply \eqref{autofinali}. Then \eqref{Phinftys}, \eqref{stima R 3} (applied to $ {\cal R}_0 = {\cal R} $) 
imply \eqref{stima Phi infty}.  

Estimate \eqref{Phi infty (ph) - I} follows from 
\eqref{interpolazione norme miste} (in $ H^s_x (\T) $), 
Lemma \ref{Aphispace}, and the bound \eqref{stima Phi infty}. 

In the reversible case, since $\Phi_\infty$, $\Phi_{\infty}^{-1}$ are reversibility preserving (see Corollary \ref{lem:convPhi}), and $\mL_0$ is reversible (see Remark \ref{reversibilitˆ step 5} and Lemma \ref{lemma:mostro}), we get that $\mL_\infty$ is reversible too. The eigenvalues $ \mu_j^{\infty} $ are purely imaginary by Lemma \ref{lemma42}.

In the Hamiltonian case, 
$ \mL_0 \equiv \mL_5 $ is Hamiltonian, $\Phi_{\infty}$ is symplectic, and therefore ${\cal L}_{\infty}
= \Phi_{\infty}^{-1} \mL_5 \Phi_{\infty}$ (see \eqref{Lfinale})  
is Hamiltonian, namely $\mD_\infty$ has the structure  $ \mD_\infty = \pa_x \mathcal{B} $, 
where $\mathcal{B} = \mathrm{diag}_{j \neq 0} \{ b_j \}$ is self-adjoint. 
This means that $b_j \in \R$, and therefore $\mu_j^\infty = \ii j b_j$ are all purely imaginary.
\rule{2mm}{2mm}

\subsection{Proof of Theorem \ref{thm:abstract linear reducibility}} \label{subsec:proof of thm abstract linear reducibility}

\noindent
{\sc Proof  of ${\bf ({S}i)}_{0}$, $i=1,\ldots,4$.} 
Properties \eqref{def:Lj}-\eqref{Rsb} in ${\bf({S}1)}_0$
hold by \eqref{L0}-\eqref{D0R0} with $ \mu_j^0 $ defined in \eqref{mu-j-nu} and
 $ r_j^0(\lm) = 0 $ (for \eqref{Rsb} recall  that $ N_{-1} := 1 $, see \eqref{defN}). 
Moreover, since $m_1$, $m_3$ are real functions, $\mu_j^0$ are purely imaginary,  
$\m_j^0 = \overline{{\mu}_{-j}^0}$ and $\mu_j^0 = - \mu_{-j}^0$. In the reversible case, remark \ref{reversibilitˆ step 5}
implies that $\mR_0 := \mR$, $\mL_0 := \mL_5$ are reversible operators.
Then there is nothing else to verify.  

${\bf({S}2)}_0 $ holds extending from  $ \Lambda^\g_0 :=   \Lambda_o $ to 
  $ \Lambda $ the eigenvalues  
$\mu_{j}^0 (\l) $, namely extending the functions 
$ m_1 (\l) $, $ m_3 (\l) $ to $ {\tilde m}_1 (\l) $, $ {\tilde m}_3 (\l) $,  
preserving the sup norm and the Lipschitz semi-norm,  by Kirszbraun theorem. 

${\bf({S}3)}_0 $ follows by \eqref{stima R 2}, for $s = \mathfrak s_0 , \mathfrak s_0 + \b$, and 
\eqref{norma bassa u riducibilitˆ}, \eqref{alpha-beta}.

${\bf({S}4)}_0 $ is trivial because, by definition, $\L_0^\g(u_1) = \L_o = \L_0^{\g-\rho}(u_2)$.

\subsubsection{The reducibility step}\label{the-reducibility-step}

We now describe the generic inductive step, showing how to define $ {\cal L}_{\nu+1 } $ (and $ \Phi_\nu $, $ \Psi_\nu $, etc).
To simplify notations, in this section we drop the index $ \nu $ and we write $ + $ for $ \nu + 1$.
We have
\begin{eqnarray}
{\cal L} \Phi h & = & 
\o \cdot \partial_{\vphi} (\Phi(h)) + {\cal D} \Phi h + {\cal R} \Phi h  \nonumber \\
& = & \o \cdot \partial_{\vphi}  h  + \Psi  \o \cdot \partial_{\vphi} h 
+ (\o \cdot \partial_{\vphi} \Psi ) h +
{\cal D} h + {\cal D} \Psi  h + {\cal R} h + {\cal R} \Psi h \nonumber \\
& = & \Phi \Big( \o \cdot \partial_{\vphi} h + 
{\cal D} h\Big) + \Big(\o \cdot \partial_{\vphi} \Psi + \left[{\cal D}, \Psi \right] + \Pi_{N} {\cal R}\Big) h + 
\Big(\Pi_{N}^{\bot} {\cal R} + {\cal R} \Psi \Big) h \label{1stepNM}
\end{eqnarray}
where $[{\cal D}, \Psi ] := {\cal D} \Psi - \Psi  {\cal D} $ and 
$\Pi_{N}{\cal R}$ is defined in \eqref{SM}. 

\begin{remark}
The application of the smoothing operator $ \Pi_N $ is necessary since we are 
performing a differentiable Nash-Moser scheme. Note also that $ \Pi_N $ regularizes only in time 
(see \eqref{SM}) because
the loss of derivatives of the inverse operator is only in $ \vphi $ (see \eqref{solomo} and the bound on the small divisors \eqref{Omgj}). 
\end{remark}

We look for a solution of the {\it homological} equation  
\be\label{eq:homo}
\o \cdot \partial_{\vphi} \Psi  + \left[{\cal D}, \Psi \right] + \Pi_{N} {\cal R} = [{\cal R}] \qquad {\rm where} 
\qquad [{\cal R}]  := {\rm diag}_{j \in \Z} {\cal R}^{j}_{j}(0) \, .
\end{equation}

\begin{lemma} \label{lemma:redu} {\bf (Homological equation)} 
For all $ \l \in {\L}_{\nu+1}^{\gamma} $, (see \eqref{Omgj}) there exists a unique solution
$ \Psi := \Psi (\vphi) $ of the homological equation \eqref{eq:homo}. The map $ \Psi $ satisfies 
\be\label{PsiR}
\left|\Psi\right|_{s}^{\Lipg} \leq C  N^{2\tau + 1} \gamma^{-1} \left|{\cal R}\right|_{s}^{\Lipg} \, . 
\ee
Moreover if $\g/ 2 \leq \g_1, \g_2 \leq 2\g$ and if $ u_1(\l)$, $  u_2(\l) $   are Lipschitz functions, then $\forall s \in [\mathfrak s_0, \mathfrak s_0 + \b] $, 
$\l \in \L_{\nu + 1}^{\gamma_1}(u_1) \cap \L_{\nu + 1}^{\gamma_2}(u_2)$
\be\label{differenza finita Psi}
| \Delta_{12} \Psi |_s \leq C N^{2\t + 1}\g^{-1}\Big(|\mR(u_2)|_s \Vert u_1 - u_2 \Vert_{\mathfrak s_0 + \s_2} + |\Delta_{12}\mR|_s \Big)  
\ee 
where we define  $ \Delta_{12} \Psi := \Psi ( u_1)  -\Psi (u_2) $.

In the reversible case, $ \Psi $ is reversibility-preserving. 
\end{lemma}

\begin{pf}
Since ${\cal D} := {\rm diag}_{j \in \Z} (\mu_{j})$ we have
$ [{\cal D}, \Psi ]_{j}^{k} =  (\mu_j - \mu_k) \Psi_{j}^{k}(\vphi) $ and
\eqref{eq:homo} amounts to 
$$
\om \cdot \partial_{\vphi} \Psi_{j}^{k}(\vphi) + 
(\mu_{j} - \mu_{k}) \Psi_{j}^{k}(\vphi) + {\cal R}_{j}^{k}(\vphi) =  [{\cal R}]_{j}^{k}  \, , \quad \forall j, k \in \Z \, ,
$$
whose solutions are  $ \Psi_{j}^{k}(\vphi) 
= \sum_{l \in \Z^\nu} \Psi_{j}^{k}(l) e^{\ii l \cdot \vphi} $  
with coefficients 
\be\label{solomo}  
{\Psi}_{j}^{k}(l) := 
\begin{cases}
\dfrac{{\cal R}_{j}^{k}(l)}{ \delta_{ljk}(\l) }  \quad \,  
& 
\text{if} \  (j-k, l) \neq (0,0) \ \ \text{and} \ \  |l |  \leq N \, , \ \ 
\text{where} \ \ \delta_{ljk}(\l)  := \ii \om \cdot l + \mu_j - \mu_k, \\
0 
& \text{otherwise.} 
\end{cases}
\ee
Note that, for all $ \l \in \L_{\nu + 1}^{\g} $, by \eqref{Omgj} and \eqref{omdio}, if  $ j \neq k $  or $  l \neq 0 $ the divisors  
$  \delta_{ljk}(\l) \neq 0 $. Recalling the definition of the $ s $-norm in \eqref{matrix decay norm} 
we deduce by \eqref{solomo},  \eqref{Omgj}, \eqref{omdio}, that 
\be\label{supomo}
| \Psi |_s \leq \g^{-1} N^\t | {\cal R} |_s \, , \quad \forall \l \in \L_{\nu+1}^\g \, .  
\ee
For $ \l_{1}, \l_{2} \in \L_{\nu + 1}^{\gamma} $, 
\begin{equation}\label{italia}
|\Psi_{j}^{k}(l) (\l_{1}) - \Psi_{j}^{k}(l) (\l_{2})| \leq 
\frac{|{\cal R}_{j}^{k}(l) (\l_{1}) - {\cal R}_{j}^{k}(l) (\l_{2})|}{|\delta_{ljk}(\l_{1})|} \, 
+
|{\cal R}_{j}^{k}(l) (\l_{2})| \, \frac{|\delta_{ljk}(\l_{1}) - \delta_{ljk}(\l_{2})|}{|\delta_{ljk}(\l_{1})| |\delta_{ljk}(\l_{2})|} 
\end{equation}
and, since $\om = \lm \bar\om$, 
\begin{eqnarray}
|\delta_{ljk}(\l_{1}) - \delta_{ljk}(\l_{2})| & \stackrel{\eqref{solomo}} = & 
|(\l_{1}-\l_{2})\bar\o\cdot l + (\mu_{j} - \mu_{k})(\l_{1}) - (\mu_{j}-\mu_{k})(\l_{2})| \\
& \stackrel{\eqref{mu-j-nu}} \leq & |\l_{1}-\l_{2}||\bar\o\cdot l | + |m_3(\l_{1}) - m_3(\l_{2})| |j^{3}-k^{3} | + |m_1(\l_{1}) - m_1(\l_{2})| | j-k | \nonumber\\
& & + \,  |r_{j}(\l_{1}) - r_{j}(\l_{2})| + |r_{k}(\l_{1}) - r_{k}(\l_{2})| \nonumber \\
& \lessdot & |\l_{1}-\l_{2}| \Big( | l | + \e\g^{-1}   |j^{3}-k^{3} | + \e\g^{-1}  | j-k | 
+ \e\g^{-1}   \Big)   \label{ultieq} 
\end{eqnarray}
because 
\[ 
\g|m_3|^{\rm lip} 
= \g|m_3 - 1|^{\rm lip} 
\leq |m_3 - 1|^{\Lipg} 
\leq  \e C, 
\quad   
|m_1|^{\Lipg} \leq \e C,
\quad 
|r_{j} |^{\Lipg} \leq \e C 
\quad 
\forall j \in \Z.
\] 
Hence, for $ j \neq k $,  $ \e \g^{-1} \leq 1 $, 
\begin{eqnarray}
\frac{|\delta_{ljk}(\l_{1}) -  \delta_{ljk}(\l_{2})|}{|\delta_{ljk}(\l_{1})||\delta_{ljk}(\l_{2})|}  \!\!\!
& \stackrel{\eqref{ultieq}, \eqref{Omgj}}  \lessdot & \!\!\!
 |\l_{1}-\l_{2}| \Big( | l | +   |j^{3}-k^{3} | \Big) \frac{\left\langle l\right\rangle^{2\tau}}{\gamma^{2}\left|j^{3}-k^{3}\right|^{2}} 
 \lessdot |\l_{1} - \l_{2}|  N^{2\tau + 1} \gamma^{-2} \label{francia}
\end{eqnarray}
for $ |l| \leq N $. 
Finally, recalling  \eqref{matrix decay norm}, the bounds \eqref{italia}, \eqref{francia} and \eqref{supomo} imply 
\eqref{PsiR}.
Now we prove \eqref{differenza finita Psi}. 
By \eqref{solomo}, for any  $ \l \in \L_{\nu + 1}^{\g_1} (u_1) \cap \L_{\nu + 1}^{\g_2} (u_2) $,  $ l \in \Z^{\nu} $, $ j \neq k $, we get 
\begin{equation}\label{delta coefficienti Psi}
\Delta_{12}\Psi_j^k(l) = 
\frac{\Delta_{12}\mR_j^k(l)}{\delta_{ljk}(u_1)} 
- \mR_j^k(l)(u_2) \frac{\Delta_{12}\delta_{ljk}}{\delta_{ljk}(u_1) \delta_{ljk}(u_2)} 
\end{equation}
where
\begin{eqnarray}
|\Delta_{12}\delta_{ljk}| & = & |\Delta_{12}(\mu_j - \mu_k)| \leq |\Delta_{12}m_3 | \, | j^{3} - k^{3}| + 
|\Delta_{12}m_1|\, |j - k| + |\Delta_{12}r_j | + | \Delta_{12}r_k| \nonumber \\
 & \stackrel{\eqref{coefficienti costanti 2},  \eqref{Delta12 rj} } \lessdot & 
 \e |j^3 - k^3| \Vert u_1 - u_2 \Vert_{\mathfrak s_0 + \s_2} \, . \label{cina}
\end{eqnarray}
Then \eqref{delta coefficienti Psi}, \eqref{cina}, $\e \gamma^{-1} \leq 1$, $\g_{1}^{-1}, \g_{2}^{-1} \leq \g^{-1}$ imply
$$
|\Delta_{12}\Psi_j^k(l)| \lessdot N^{2\tau} \gamma^{-1} 
\Big(|\Delta_{12}\mR_j^k (l)| + |\mR_j^k (l)(u_2)| \Vert u_1 - u_2 \Vert_{\mathfrak s_0 + \s_2} \Big)
$$
and so  \eqref{differenza finita Psi} (in fact, \eqref{differenza finita Psi} holds with $2\t$ instead of $2\t+1$).

In the reversible case   
$ \ii \om \cdot l + \mu_j - \mu_k \in \ii\R $, $\overline{{\mu}_{-j}} = \mu_j$ and $ \mu_{-j} = - \mu_j $. 
Hence Lemma \ref{lem:PR} and  \eqref{solomo} imply
$$
\overline{ {\Psi}_{-j}^{-k}(-l)} =  
\frac{ \overline{{\cal R}_{-j}^{-k}(-l)}}{ {-\ii \o \cdot (-l) + \overline{{\mu}_{-j}} - \overline{{\mu}_{-k}}} }
=  \frac{{\cal R}_{j}^{k}(l)}{\ii \o \cdot l + \mu_{j} - \mu_k } =  {\Psi}_j^k(l) 
$$
and so  $ \Psi $ is real, again by Lemma \ref{lem:PR}.
Moreover, since $ {\cal R} : X \to Y $,
$$
\Psi^{-k}_{-j}(-l)  =   \frac{  {\cal R}^{-k}_{-j}(-l)}{\ii \o \cdot (-l) + \mu_{-j} - \mu_{-k} } = 
\frac{- {\cal R}^{k}_{j}(l)}{ \ii\o \cdot (-l) - \mu_{j} + \mu_{k} } 
= {\Psi}^{k}_{j}(l) 
$$
which implies  $ \Psi : X \rightarrow X $ by Lemma \ref{lem:PR}. Similarly we get $ \Psi : Y \rightarrow Y  $.
\end{pf}

\begin{remark}
In the Hamiltonian case $ {\cal R} $ is Hamiltonian and the solution $ \Psi $ in \eqref{solomo} of the homological equation is Hamiltonian, because 
$ \overline{ \d_{l,j,k} } = \d_{-l,k,j} $ 
and, in terms of matrix elements, an operator $G(\ph)$ is self-adjoint if and only if 
$ \overline{ G_j^k(l) } = G_k^j(-l) $.
\end{remark}

Let $ \Psi $ be the solution of the homological equation \eqref{eq:homo} which has been constructed in Lemma \ref{lemma:redu}.
By Lemma \ref{lem:inverti}, if $ C(\mathfrak s_0) | \Psi |_{\mathfrak s_0} < 1 /2 $
then  $ \Phi := I + \Psi $ is invertible and
 by \eqref{1stepNM} (and \eqref{eq:homo}) we deduce that
\be\label{defL+}
{\cal L}_{+} := \Phi^{-1} {\cal L} \Phi =  \o \cdot \partial_{\vphi} + {\cal D}_{+} + {\cal R}_{+} \,,
\ee
where 
\begin{equation}\label{D+}
{\cal D}_{+} := {\cal D} + [{\cal R}] \, , \quad 
{\cal R}_{+} := \Phi^{-1} \Big(\Pi_N^{\bot} {\cal R} + {\cal R} \Psi 
- \Psi  [{\cal R}] \Big). 
\end{equation}
Note that $\mL_+$ has the same form of $ {\cal L} $,  but the remainder $ \mR_+ $ is 
the sum of a quadratic function of $ \Psi, {\cal R} $ and a  remainder supported on high modes.

\begin{lemma}\label{nuovadiagonale} 
{\bf (New diagonal part).}  
The eigenvalues of 
\[
{\cal D}_{+} = {\rm diag}_{j \in \Z}  \{ \mu^{+}_{j}(\l) \}, 
\quad \text{where} \ \ 
\mu^{+}_{j} 
:= \mu_{j} +  {\cal R}^{j}_{j}(0) 
= \mu_j^0 + r_j + \mR_j^j (0) 
= \mu_j^{0} + r_j^+, 
\quad
r_j^+ := r_j + \mR_j^j (0), 
\]
satisfy $\mu_{j}^{+} = \overline{{\mu}^{+}_{-j}}$ and 
\be\label{ultrasim} 
|\mu^{+}_{j} - \mu_{j} |^{\rm lip}
= |r^{+}_{j} - r_{j} |^{\rm lip} 
= |\mR_j^j(0)|^\lip
\leq \left|{\cal R}\right|_{\mathfrak s_0}^{\rm lip},\quad \forall j \in \Z. 
\ee
Moreover if $ u_1 (\l)$, $u_2 (\l)$ are Lipschitz functions, then for all 
$\l\in \L_{\nu}^{\g_1}(u_1) \cap \L_{\nu}^{\g_2}(u_2)$
\be\label{spagna}
|\Delta_{12} r_j^+ - \Delta_{12} r_j| \leq |\Delta_{12}\mR|_{\mathfrak s_0}\,.
\ee 
In the reversible case, all the $\mu_j^{+}$ are purely imaginary and satisfy $\mu^{+}_{j} = - \mu^{+}_{-j}$ for all $j \in \Z$.
\end{lemma}

\begin{pf}
The estimates \eqref{ultrasim}-\eqref{spagna} follow using \eqref{Aii} because
$ |  {\cal R}^{j}_{j}(0) |^{\rm lip} = $ $ |{\cal R}^{(l,j)}_{(l,j)} |^{\rm lip}\leq $  
$ |{\cal R} |_0^{\rm lip}  \leq $ $ |{\cal R} |_{\mathfrak s_0}^{\rm lip} $ and 
$$
|\Delta_{12}r^{+}_{j} - \Delta_{12} r_{j} | 
= | \Delta_{12} {\cal R}^{j}_{j}(0) |
= |\Delta_{12}{\cal R}^{(l,j)}_{(l,j)} | \leq  |\Delta_{12}{\cal R} |_0 \leq 
 |\Delta_{12}{\cal R} |_{\mathfrak s_0} \, .
$$
Since $ {\cal R} $ is real, by Lemma \ref{lem:PR}, 
$$
{\cal R}^{k}_{j}(l)\,=\, \overline{{\cal R}^{-k}_{-j}(-l)} 
\qquad \Longrightarrow \qquad  \mR_{j}^{j}(0) = \overline{{\mR}_{-j}^{-j}(0)} 
$$
and so 
 $ \mu_{j}^+ 
= \overline{{\mu}_{-j}^{+}} $. If  $\mR$ is also reversible, by  Lemma \ref{lem:PR},
$$ 
{\cal R}^{k}_{j}(l) = - {\cal R}^{-k}_{-j}(-l) \, , \quad {\cal R}^{k}_{j}(l) = \overline{{\cal R}^{-k}_{-j}(-l)} 
= - \overline{{ {\cal R}^{k}_{j}(l)}} \, .
$$
We deduce that 
$ {\cal R}^{j}_{j}(0)  = - {\cal R}^{-j}_{-j}(0) $, 
$ {\cal R}^{j}_{j}(0) \in \ii \R $ and therefore,  $  \mu^{+}_{j} = -  \mu^{+}_{-j} $ and 
$\mu_{j}^{+} \in \ii \R$. 
\end{pf}

\begin{remark}\label{rem:Ham10}
In the Hamiltonian case, 
$\mD_\nu $ is Hamiltonian, namely $ \mD_\nu = \pa_x \mathcal{B} $ 
where $\mathcal{B} = \mathrm{diag}_{j \neq 0} \{ b_j \}$ is self-adjoint. 
This means that $b_j \in \R$, and therefore all $\mu_j^\nu = \ii j b_j $ are purely imaginary.
\end{remark}

\subsubsection{The iteration}

Let $\nu \geq 0$, and suppose that the statements ${\bf({S}i)_{\nu}}$ are true. 
We prove  $({\bf Si})_{\nu+1}$, $i=1,\ldots,4$. 
To simplify notations we write $|\cdot|_s$ instead of $|\cdot|_s^{\Lipg}$.

\smallskip

\noindent
{\sc Proof of  $({\bf S1})_{\nu + 1}$}. 
By ${\bf (S1)_\nu} $, the eigenvalues $\mu_j^\nu$ are defined on $\L_\nu^\g$. Therefore 
the set $\L_{\nu+1}^\g$ is well-defined. 
By Lemma \ref{lemma:redu}, for all $ \l \in \L_{\nu+1}^{\g} $ there exists a real solution 
$ \Psi_{\nu} $ of the homological equation \eqref{eq:homo} which satisfies,   
$ \forall s \in [\mathfrak s_0, q- \s - \b ] $, 
\be\label{Psinu}
\left|\Psi_{\nu}\right|_{s} \stackrel{\eqref{PsiR}}{\lessdot}  
N_{\nu}^{2\tau + 1}\left|{\cal R}_{\nu}\right|_{s} \gamma^{-1} 
\stackrel{\eqref{Rsb}} 
{\lessdot} \left|{\cal R}_{0}\right|_{s  + \beta} \gamma^{-1} N_{\nu}^{2\tau + 1}\ N_{\nu-1}^{- \a} 
\ee
which is \eqref{Psinus} at the step $ \nu +1 $. 
In particular, for $ s = \mathfrak s_0 $,  
\be\label{Psinu0}
C(\mathfrak s_0) \left|\Psi_{\nu}\right|_{\mathfrak s_0} \stackrel{\eqref{Psinu}} \leq  
C(\mathfrak s_0) \left|{\cal R}_{0}\right|_{\mathfrak s_0  + \beta} \gamma^{-1} N_{\nu}^{2\tau + 1}\ N_{\nu-1}^{- \a} 
\stackrel{\eqref{piccolezza1}} \leq 1/2   
\ee
for $ N_0 $ large enough.
Then the map  $ \Phi_{\nu} := I + \Psi_\nu $  is invertible and, by \eqref{PhINV}, 
\be\label{Phis0}
\left|\Phi_{\nu}^{-1}\right|_{\mathfrak s_{0}} \leq 2 \,  , \quad  \left|\Phi_{\nu}^{-1}\right|_{s} \leq 1 + C(s) | \Psi_\nu |_s \, . 
\ee
Hence \eqref{defL+}-\eqref{D+} imply
$  {\cal L}_{\nu + 1} := $
$ \Phi_{\nu}^{-1} {\cal L}_{\nu} \Phi_{\nu} = $ $ \o \cdot \partial_{\vphi} + {\cal D}_{\nu + 1} + {\cal R}_{\nu + 1} $
where (see Lemma \ref{nuovadiagonale})
\begin{equation}\label{Dnu+1}
{\cal D}_{\nu + 1} := {\cal D}_{\nu} + [{\cal R}_{\nu}] = {\rm diag}_{j \in \Z} (\mu_j^{\nu + 1}) \, , 
\quad \mu_j^{\nu + 1} := \mu_j^{\nu} + ({\cal R}_{\nu})_{j}^{j}(0)   \, ,  
\end{equation}
with $\mu_j^{\nu + 1} = \overline{\mu_{-j}^{\nu + 1}} $ and
\begin{equation}\label{Rnu+1}
{\cal R}_{\nu+1} := \Phi_\nu^{-1} H_{\nu},\quad H_{\nu}:= \Pi_{N_\nu}^{\bot} {\cal R}_\nu + {\cal R}_\nu \Psi_\nu 
- \Psi_\nu  [{\cal R}_\nu]  \, . 
\end{equation} 
In the reversible case, $\mR_\nu : X \rightarrow Y$, therefore, by Lemma \ref{lemma:redu}, $\Psi_\nu$, $\Phi_\nu$, $\Phi_{\nu}^{-1}$ are reversibility preserving, and then, by formula \eqref{Rnu+1}, also $\mR_{\nu + 1} : X \rightarrow Y$.

Let us prove the estimates \eqref{Rsb} for $ {\cal R}_{\nu + 1}  $.  For all $ s \in [\mathfrak s_0, q - \s - \b ] $ we have
\begin{eqnarray}
|{\cal R}_{\nu + 1} |_{s} & \!\! \!\! \!\!  \stackrel{\eqref{Rnu+1}, \eqref{interpm Lip}} {\leq_s} \!\!\!\!  \!\!  & 
| \Phi_\nu^{-1} |_{\mathfrak s_{0}}   \Big( |\Pi_{N_\nu}^\bot  {\cal R}_\nu |_s + |{\cal R}_\nu |_s |\Psi_\nu |_{\mathfrak s_{0}} +  
|{\cal R}_\nu |_{\mathfrak s_{0}} |\Psi_\nu |_s  \Big) + 
| \Phi_\nu^{-1} |_{s}   \Big( |\Pi_{N_\nu}^\bot  {\cal R}_\nu |_{\mathfrak s_{0}} + |{\cal R}_\nu |_{\mathfrak s_{0}} |\Psi_\nu |_{\mathfrak s_{0}} \Big) 
\nonumber \\
& \!\!\!\! \!\!   \stackrel{\eqref{Phis0}} {\leq_s} \!\!\!\! \!\!  & 2
 \Big(  |\Pi_{N_\nu}^\bot  {\cal R}_\nu |_s + |{\cal R}_\nu |_s |\Psi_\nu |_{\mathfrak s_{0}} +  
|{\cal R}_\nu |_{\mathfrak s_{0}} |\Psi_\nu |_s  \Big)   + (1 + | \Psi_\nu |_s) 
\Big( |\Pi_{N_\nu}^\bot  {\cal R}_\nu |_{\mathfrak s_{0}} + |{\cal R}_\nu |_{\mathfrak s_{0}} |\Psi_\nu |_{\mathfrak s_{0}} \Big) \nonumber \\
& \!\!  \stackrel{\eqref{Psinu0}} {\leq_s} \!\!  &   |\Pi_{N_\nu}^\bot  {\cal R}_\nu |_s + |{\cal R}_\nu |_s |\Psi_\nu |_{\mathfrak s_{0}} +  
|{\cal R}_\nu |_{\mathfrak s_{0}} |\Psi_\nu |_s  \!\!  \stackrel{\eqref{PsiR}} {\leq_s} \!\!    |\Pi_{N_\nu}^\bot  {\cal R}_\nu |_s + N_\nu^{2\t+1} \g^{-1} |{\cal R}_\nu |_s | {\cal R}_\nu |_{\mathfrak s_{0}} \, . \label{Rsgen}
\end{eqnarray}
Hence \eqref{Rsgen} and \eqref{smoothingN} imply
\be
|{\cal R}_{\nu + 1} |_s {\leq_s} 
N_\nu^{-\b}  | {\cal R}_\nu |_{s+\b} +  N_\nu^{2\t+1} \g^{-1} |{\cal R}_\nu |_s | {\cal R}_\nu |_{\mathfrak s_{0}} \label{sch1} 
\ee
which shows that the iterative scheme is quadratic plus a super-exponentially small term. In particular
$$
|{\cal R}_{\nu + 1} |_{s} \!\! 
\stackrel{\eqref{sch1}, \eqref{Rsb}} {\leq_s} \!\!  N_\nu^{-\b} |{\cal R}_0|_{s+\b} N_{\nu-1} +  
N_\nu^{2\t+1} \g^{-1}   |{\cal R}_0|_{s+\b}  |{\cal R}_0|_{\mathfrak s_{0}+\b} N_{\nu-1}^{- 2\a} 
\!\!  \stackrel{\eqref{defbq}, \eqref{alpha-beta}, \eqref{piccolezza1}}\leq \!\!  |{\cal R}_0|_{s+\b} N_{\nu}^{-\a}
$$
($ \chi = 3 / 2 $)
which is the first inequality of \eqref{Rsb} at the step $ \nu +1 $. 
The next key step  is to control the divergence of the high norm $ | {\cal R}_{\nu+1} |_{s+\b} $.
By \eqref{Rsgen} (with $ s + \b $ instead of $ s $) we get
\be
|{\cal R}_{\nu + 1} |_{s+\b} \,
{\leq_{s+\b}} \,  | {\cal R}_\nu |_{s+ \b} +  N_\nu^{2\t+1} \g^{-1} |{\cal R}_\nu |_{s+\b}| {\cal R}_\nu |_{\mathfrak s_{0}}  \label{sch2}
\ee
(the difference with respect to \eqref{sch1} is that we do not apply to 
$ | \Pi_{N_\nu}^\bot {\cal R}_{\nu} |_{s+\b} $ any smoothing). 
Then \eqref{sch2}, \eqref{Rsb}, \eqref{piccolezza1}, \eqref{alpha-beta} imply 
the 
inequality
$$
|{\cal R}_{\nu + 1} |_{s+\b} \leq  C(s+\b) | {\cal R}_\nu |_{s+ \b}, 
$$
whence, iterating,  
$$
|{\cal R}_{\nu + 1} |_{s+\b}  \leq N_{\nu} |{\cal R}_0 |_{s+ \b} 
$$
for $ N_0 := N_0 (s,\b) $ large enough, which is the second inequality of \eqref{Rsb} with index $ \nu +1 $. 

By Lemma \ref{nuovadiagonale}  
the eigenvalues $ \mu_j^{\nu + 1} := \mu_j^0 + r_j^{\nu + 1} $, 
defined on $ \L_{\nu+1}^{\g} $, satisfy  
$\mu_j^{\nu + 1} = \overline{{\mu}_{-j}^{\nu + 1}} $, and, 
in the reversible case, the $\mu_{j}^{\nu + 1}$ are purely imaginary and 
$\mu_j^{\nu + 1} = - \mu_{-j}^{\nu + 1}$. 

It remains only to prove \eqref{rjnu bounded} for $\nu+1$, which is proved below.

\smallskip

\noindent
{\sc Proof of ${\bf({S}2)}_{\nu + 1} $}. By \eqref{ultrasim},  
\be\label{lnu+1lnu}
|\mu_j^{\nu + 1} - \mu_j^{\nu} |^{\Lipg} =  |r_j^{\nu + 1} - r_j^{\nu} |^{\Lipg} 
 \leq |{\cal R}_{\nu} |_{\mathfrak s_{0}}^{\Lipg} 
 \stackrel{\eqref{Rsb}}\leq \left|{\cal R}_{0}\right|^{\Lipg}_{\mathfrak s_{0} + \beta}N_{\nu - 1}^{- \alpha}\,.
\ee
By Kirszbraun theorem, we extend  the function $ \mu_j^{\nu + 1} - \mu_j^{\nu} = r_j^{\nu + 1} - r_j^{\nu} $
to the whole $ \L $, still satisfying \eqref{lnu+1lnu}. 
In this way  we define $ \tilde \mu_j^{\nu + 1}$.
Finally \eqref{rjnu bounded} follows  summing all the terms in \eqref{lnu+1lnu} and using \eqref{stima R 3}. 

\smallskip

\noindent
{\sc Proof of ${\bf({S}3)}_{\nu + 1} $}. 
Set, for brevity,  
 $$
 \mR_{\nu}^{i} := \mR_{\nu}(u_i),\quad \Psi_{\nu - 1}^i := \Psi_{\nu - 1}(u_i),\quad  \Phi_{\nu - 1}^{i} := \Phi_{\nu - 1}(u_i), 
 \quad H_{\nu - 1}^i := H_{\nu - 1} (u_i) \, ,  \quad i:= 1, 2 \, ,
 $$
which are all operators defined for $\l \in \L_{\nu}^{\gamma_1}(u_1) \cap \L_{\nu}^{\gamma_2}(u_2) $. By Lemma \ref{lemma:redu} one can construct $\Psi_{\nu}^{i}:= \Psi_{\nu}(u_i)$, $\Phi_{\nu}^i := \Phi_{\nu}(u_i)$, $i = 1, 2$, 
for all $\l \in \Lambda_{\nu + 1}^{\g_1}(u_1) \cap \Lambda_{\nu + 1}^{\g_2}(u_2)$. 
One has
\begin{align} 
\vert \Delta_{12} \Psi_{\nu} \vert_{\mathfrak s_0} 
& 
\stackrel{\eqref{differenza finita Psi}}{\lessdot} 
 N_\nu^{2\t+1} \g^{-1} \Big( |\mR_\nu(u_2)|_{\mathfrak s_0} \| u_2 - u_1 \|_{\mathfrak s_0 + \s_2} 
+ |\D_{12} \mR_\nu|_{\mathfrak s_0} \Big)
\notag
\\
& 
\stackrel{\eqref{Rsb},\eqref{derivate-R-nu}}{\lessdot} 
N_\nu^{2\t+1} N_{\nu-1}^{-\a} \g^{-1} \big( |\mR_0|_{\mathfrak s_0+\b} + \e \big) 
\| u_2 - u_1 \|_{\mathfrak s_0 + \s_2}
\notag 
\\
& 
\stackrel{\eqref{stima R 3}, \eqref{norma bassa u riducibilitˆ}}{\lessdot} 
 N_\nu^{2\t+1} N_{\nu-1}^{-\a} \e \g^{-1} 
\| u_2 - u_1 \|_{\mathfrak s_0 + \s_2}
\leq  \| u_2 - u_1 \|_{\mathfrak s_0 + \s_2} .
\label{delta Psi nu bassa}
\end{align}
for $ \e \g^{-1} $ small (and \eqref{alpha-beta}). 
By \eqref{derivata-inversa-Phi}, applied to $\Phi:= \Phi_{\nu}$, and \eqref{delta Psi nu bassa}, we get 
 \begin{equation}\label{delta Phi nu s}
 \vert \Delta_{12} \Phi_{\nu}^{-1} \vert_{s} \leq_s  
 \big( \vert \Psi_{\nu}^{1} \vert_s + \vert \Psi_{\nu}^{2} \vert_{s} \big)\Vert u_1 - u_2 \Vert_{\mathfrak s_0 + \s_2} + \vert \Delta_{12} \Psi_{\nu} \vert_{s} 
 \end{equation}
 which implies for $s = \mathfrak s_0$, and using \eqref{Psinus}, \eqref{piccolezza1}, \eqref{delta Psi nu bassa}
 \be\label{derivata-inversa-Phi-bassa}
 \vert \Delta_{12} \Phi_{\nu}^{-1} \vert_{\mathfrak s_0} \lessdot \Vert u_1 - u_2 \Vert_{\mathfrak s_0 + \s_2}.
 \ee
Let us prove the estimates \eqref{derivate-R-nu} for $\Delta_{12}\mR_{\nu + 1}$, which is defined on $\l \in \Lambda_{\nu + 1}^{\g_1}(u_1) \cap \Lambda_{\nu + 1}^{\g_2}(u_2)$.
For all $s \in [ {\mathfrak s}_{0}, {\mathfrak s}_{0}+\b]$, using the interpolation \eqref{interpm} and \eqref{Rnu+1}, 
\begin{equation}\label{derivata-R-nu+1}
|\Delta_{12}{\cal R}_{\nu + 1} |_{s} \! \stackrel{}{\leq_{s}}  \! 
|\Delta_{12}\Phi_{\nu}^{-1} |_{s} |H_{\nu}^{1} |_{\mathfrak s_{0}} + \vert \Delta_{12}\Phi_{\nu}^{-1}\vert_{\mathfrak s_{0}} \vert H_{\nu}^{1} \vert_{s}   \! +   
|(\Phi_{\nu}^2 )^{-1}|_s |\Delta_{12}H_{\nu} |_{\mathfrak s_0} + 
|(\Phi_{\nu}^2 )^{-1} |_{\mathfrak s_{0}} |\Delta_{12}H_\nu |_s \, .
\end{equation}
We estimate the above terms separately.  Set for brevity 
$ A^\nu_{s}  := | \mR_\nu(u_1) |_s + | \mR_\nu(u_2) |_s $. 
By \eqref{Rnu+1} and \eqref{interpm},
\begin{eqnarray}
| \Delta_{12}H_{\nu}|_s  \! \!  \! \!  \! \! 
& \leq_s &  \! \!  \! \!  \! \!  \left|\Pi_{N_{\nu}}^{\bot}\Delta_{12}{\cal R}_{\nu}\right|_{s} + |\Delta_{12}\Psi_{\nu}|_{s}
|{\cal R}_{\nu}^1 |_{\mathfrak s_{0}} + |\Delta_{12}\Psi_{\nu} |_{\mathfrak s_{0}} |{\cal R}_{\nu}^1|_{s} 
+  |\Psi_{\nu}^2 |_{s} |\Delta_{12}{\cal R}_{\nu} |_{\mathfrak s_{0}} + 
|\Psi_{\nu}^2 |_{\mathfrak s_{0}} |\Delta_{12}{\cal R}_{\nu} |_{s} 
\nonumber
\\
& 
\stackrel{\eqref{PsiR},\eqref{differenza finita Psi}}{\leq_{s}} 
&
\left|\Pi_{N_{\nu}}^{\bot}\Delta_{12}{\cal R}_{\nu}\right|_{s} 
+ 
N_{\nu}^{2\tau+1}\gamma^{-1} 
A^\nu_{\mathfrak s_0}
A^\nu_s
\Vert u_1 - u_2\Vert_{\mathfrak s_{0}+\s_2} \nonumber
\\
&  & \, + \, 
N_{\nu}^{2\tau+1}\gamma^{-1}  
A^\nu_{s}
\vert \Delta_{12}{\cal R}_{\nu} \vert_{\mathfrak s_{0}} 
+ 
N_{\nu}^{2\tau+1}\gamma^{-1}  
A^\nu_{\mathfrak s_0}
\vert\Delta_{12}{\cal R}_{\nu} \vert_s \label{stima-derivata-H-nu} \, .
\end{eqnarray}
Estimating the four terms in the right hand side of \eqref{derivata-R-nu+1} in the same way, 
using 
\eqref{delta Phi nu s}, \eqref{Rnu+1}, {\eqref{PsiR}, \eqref{differenza finita Psi}, \eqref{Psinus}, 
\eqref{derivata-inversa-Phi-bassa},  \eqref{Phis0}, \eqref{stima-derivata-H-nu}, \eqref{Rsb}, 
we deduce
\begin{eqnarray}
\vert \Delta_{12}{\cal R}_{\nu+1} \vert_{s} & {\leq_{s}} &  |\Pi_{N_{\nu}}^{\bot} \Delta_{12}\mR_{\nu}|_s + N_{\nu}^{2\tau + 1} \gamma^{-1} A_{s}^{\nu} A_{\mathfrak s_0}^{\nu} \| u_1 - u_2 \|_{\mathfrak s_0 + \s_2} \nonumber\\
& &+ N_{\nu}^{2\tau + 1}\g^{-1} A_{s}^{\nu} |\Delta_{12}\mR_{\nu}|_{\mathfrak s_0} + N_{\nu}^{2\tau + 1} \gamma^{-1} A_{\mathfrak s_0}^{\nu} |\Delta_{12}\mR_{\nu}|_s \label{stima-derivata-R-nu+1-s} \, . 
\end{eqnarray}
Specializing \eqref{stima-derivata-R-nu+1-s} for $ s = 
\mathfrak s_0 $ and using \eqref{stima R 3}, \eqref{smoothingN}, \eqref{Rsb}, \eqref{derivate-R-nu}, 
we deduce
$$
\vert \Delta_{12}{\cal R}_{\nu + 1}\vert_{\mathfrak s_{0}} 
\leq C  ( \e N_{\nu - 1}N_{\nu}^{-\b} + N_{\nu}^{2\tau + 1}N_{\nu - 1}^{-2\alpha}\e^{2}\gamma^{-1} ) 
\Vert u_1 - u_2 \Vert_{\mathfrak s_{0}+\s_2} \leq \e N_{\nu}^{-\alpha} \Vert u_1 - u_2 \Vert_{\mathfrak s_{0}+\s_2} 
$$
for $N_{0}$ large and $\e\gamma^{-1}$ small. Next by \eqref{stima-derivata-R-nu+1-s} with $ s =
\mathfrak s_0 + \b $
\begin{eqnarray}
\vert \Delta_{12}{\cal R}_{\nu} \vert_{\mathfrak s_{0}+\b} & \stackrel{\eqref{Rsb}, \eqref{derivate-R-nu}, \eqref{piccolezza1}} 
{\leq_{\mathfrak s_{0}+\b}} & A_{\mathfrak s_0 + \b}^{\nu} \Vert u_1 -  u_2 \Vert_{\mathfrak s_{0}+\s_2} + \vert \Delta_{12}{\cal R}_{\nu}\vert_{\mathfrak s_{0}+\b}  \nonumber\\
&\stackrel{\eqref{Rsb}\eqref{derivate-R-nu}}{\leq}&C(\mathfrak s_{0}+\b) \e N_{\nu - 1} \Vert u_1 - u_2 \Vert_{\mathfrak s_{0}+\s_2} \leq \e N_{\nu} \Vert u_1 - u_2 \Vert_{\mathfrak s_{0}+\s_2}\nonumber
\end{eqnarray}
for $N_{0}$ large enough. Finally note that \eqref{deltarj12} is nothing but \eqref{spagna}. 

\smallskip

\noindent
{\sc Proof of ${\bf({S}4)}_{\nu + 1} $}.
We have to prove that, if 
$C \e N_{\nu}^{\tau} \Vert u_1 - u_2 \Vert_{\mathfrak s_0 + \s_2} \leq \rho$, then  
$$
\lambda \in \Lambda_{\nu+1}^{\g}(u_1) \quad \Longrightarrow \quad \l \in \Lambda_{\nu+1}^{\g- \rho}(u_2) \, . 
$$
Let $ \lambda \in \Lambda_{\nu+1}^{\g}(u_1) $.
Definition \eqref{Omgj} and ${\bf({S}4)_{\nu}}$ 
(see \eqref{legno})  
imply that 
$ \Lambda_{\nu+1}^{\g}(u_1) \subseteq \Lambda_{\nu}^{\g}(u_1) \subseteq \Lambda_{\nu}^{\g- \rho}(u_2) $. 
Hence $ \l \in \Lambda_{\nu}^{\gamma - \rho}(u_2) \subset \Lambda_{\nu}^{\gamma/2}(u_2)  $. 
Then, by ${\bf({S}1)_{\nu}}$, the eigenvalues $\mu_{j}^{\nu}(\lambda, u_2(\lambda))$ are well defined. 
Now \eqref{mu-j-nu} and the estimates \eqref{coefficienti costanti 2},  \eqref{Delta12 rj} (which holds
because $ \l \in \L_{\nu}^{\g}(u_1) \cap   \L_{\nu}^{\g/2}(u_2) $) imply that
\begin{eqnarray}
|(\mu_j^{\nu} - \mu_k^{\nu})(\l, u_2(\l)) - (\mu_j^{\nu} - \mu_k^{\nu})(\l, u_1(\l))| & \leq &  |(\mu_j^{0} - \mu_k^{0})(\l, u_2(\l)) - (\mu_j^{0} - \mu_k^{0})(\l, u_1(\l))|\nonumber\\
& & + \, 2 \sup_{j \in \Z} |r_{j}^{\nu}(\l, u_2(\l)) - r_{j}^{\nu}(\l, u_1(\l))| \nonumber\\
& \leq & \e C|j^3 - k^3| \Vert u_2 - u_1 \Vert_{\mathfrak s_0 + \s_2}^{\rm sup} \, .  \label{legno3}
\end{eqnarray}
Then we conclude that for all $|l| \leq N_{\nu}$, $j \neq k$, 
using the definition of $\Lambda_{\nu+1}^\g(u_1)$ (which is \eqref{Omgj} with $\nu+1$ instead of $\nu$) and \eqref{legno3},  
\begin{eqnarray}
|\ii  \omega \cdot l + \mu_j^{\nu} (u_2) - \mu_k^{\nu} (u_2) | & \geq & |\ii  \omega \cdot l + \mu_j^{\nu} (u_1) - \mu_k^{\nu} (u_1) | - |(\mu_j^{\nu} - \mu_k^{\nu})(u_2) - (\mu_j^{\nu} - \mu_k^{\nu})(u_1) | \nonumber\\
& 
{\geq}
& \gamma |j^3 - k^3| \langle l \rangle^{-\tau} - C \e |j^3 - k^3| \Vert u_1 - u_2 \Vert_{\mathfrak s_0 + \s_2} \nonumber\\
& \geq & (\gamma - \rho)|j^3 - k^3| \langle l \rangle^{-\tau} \nonumber
\end{eqnarray}
provided $C \e N_{\nu}^{\tau} \Vert u_1 - u_2 \Vert_{\mathfrak s_0 + \s_2} \leq \rho$. 
Hence $ \l \in \L^{\g - \rho}_{\nu+1} ( u_2 ) $. This proves \eqref{legno} at the step $ \nu + 1 $.

\subsection{Inversion of ${\cal L}(u)$}\label{sec:inversione}

In \eqref{Phi 1 2 def} we have conjugated the linearized operator $ \mL$ to $\mL_5$ defined in \eqref{mL5}, 
namely $\mL = \Phi_1 \mL_5 \Phi_2\inv$. 
In Theorem \ref{teoremadiriducibilita} we have conjugated the operator $\mL_5$ 
to the diagonal operator $\mL_{\infty}$ in \eqref{Lfinale}, 
namely $ \mL_5 = \Phi_{\infty} \mL_{\infty} \Phi_{\infty}^{-1}$. 
As a consequence
\be \label{L-coniugato}
\mL = W_1 \mL_\infty W_2\inv, 
\quad W_i :=  \Phi_{i} \Phi_{\infty}, \quad i= 1,2\,.
\ee
We  first  prove that $W_1, W_2 $ and their inverses are linear bijections of $H^{s}$. 
We take 
\begin{equation}\label{gamma-tau}
\gamma  \leq \g_0 / 2  \,, \quad \tau \geq \tau_0\,.
\end{equation}

\begin{lemma}\label{stime-tame-coniugio}
Let $  \mathfrak s_{0} \leq s \leq q - \s - \b -3 $ where $ \b $ is defined in \eqref{defbq} and $ \s $ in \eqref{costanti lemma mostro}. 
Let $u:= u(\l)$ satisfy
$ \Vert u \Vert_{\mathfrak s_0 + \s + \b + 3}^{\Lipg} \leq 1 $, 
and $ \e \g^{-1} \leq \d $ be small enough. Then 
$ W_i $, $ i = 1, 2 $,  satisfy, $ \forall \l \in \L_{\infty}^{2\g}(u) $,  
\begin{equation}\label{W1W2}
\left\| W_i h\right\|_{s} + \left\| W_i^{-1}h\right\|_{s} \leq C(s) \big(\left\|h\right\|_{s } 
+ \left\|u\right\|_{s + \s + \b}\left\|h\right\|_{\mathfrak s_{0}} \big)\, , 
\end{equation}
\begin{equation}\label{tame-Phi1Phi2}
\left\| W_i h\right\|_{s}^{\Lipg} + \left\| W_i^{-1}h\right\|_{s}^{\Lipg} \leq C(s) 
\big(\left\|h\right\|_{s + 3}^{\Lipg} + \left\|u\right\|^{\Lipg}_{s + \s + \b + 3}\left\|h\right\|_{\mathfrak s_{0}+3}^{\Lipg} \big)\,.
\end{equation}
In the reversible case (i.e. \eqref{parity f} holds), $ W_i $, $ W_i^{-1} $, $ i = 1, 2 $ are  reversibility-preserving. 
\end{lemma}

\begin{pf}
The bound \eqref{W1W2}, resp. \eqref{tame-Phi1Phi2}, follows by
\eqref{stima Phi infty}, \eqref{stima Phi 12 nel lemma}, resp. \eqref{stima Lip Phi 12 nel lemma}, 
\eqref{interpolazione norme miste} and Lemma \ref{lemma astratto composizioni}. 
In the reversible case $ W_i^{\pm 1} $  are reversibility preserving because
$ \Phi_i^{\pm1} $, $ \Phi_\infty^{\pm 1} $ are reversibility preserving. 
\end{pf}

By \eqref{L-coniugato} we are reduced to show that, $ \forall \l \in \L^{2\g}_{\infty}(u) $, the operator
$$ 
{\cal L}_\infty := {\rm diag}_{j \in \Z} \{\ii \l \bar \om \cdot  l + \mu_j^\infty (\l)\} \, , \quad \mu_j^\infty (\l) = 
-\ii \big( m_3 (\l) j^3 - m_1(\l) j \big) + r_j^\infty (\l)  
$$ 
is invertible, assuming  \eqref{f = der g}  or the reversibility condition \eqref{parity f}. 

We introduce the following notation:
\be\label{proiez00}
\Pi_C u := \frac{1}{(2\p)^{\nu+1}}\, \int_{\T^{\nu+1}} u(\ph,x) \, d\ph dx, 
\ \ 
\pp u := u - \Pi_C u, 
\ \ 
H^s_{00} := \{ u \in H^s(\T^{\nu+1}) : \Pi_C u = 0 \}.
\ee
If  \eqref{f = der g}  holds, then the linearized operator $ {\cal L} $ in \eqref{mL} satisfies  
\be \label{mL mappa tutto in media nulla}
\mL : H^{s+3} \to H^s_{00}  
\ee
(for $ \mathfrak s_0 \leq s \leq q-1 $).  In the reversible case \eqref{parity f}
\be \label{mL mappa tutto in media nulla REV}
\mL   :  X \cap H^{s+3} \to   Y \cap H^s \subset H^s_{00}   \, .  
\ee
\begin{lemma}\label{lem:modo0}
Assume 
 either \eqref{f = der g} or the reversibility condition \eqref{parity f}. Then the  eigenvalue
\be\label{zero eigenvalue}
\mu_0^\infty (\l) = r^\infty_0 (\l) = 0 \, , \quad \forall \l \in \L_\infty^{2\g} (u) \,  . 
\ee
\end{lemma}

\begin{pf}
Assume  \eqref{f = der g}.  If $ r_0^\infty \neq 0 $ then there exists a solution of $ {\cal L}_\infty w = 1 $, 
which is $ w = 1 / r_0^\infty $.   
Therefore, by \eqref{L-coniugato}, 
$$
 {\cal L} W_2 [1 / r^\infty_0] = {\cal L} W_2 w = W_1 {\cal L}_\infty w = W_1 [1] 
$$
which is a contradiction because 
$  \Pi_C W_1 [1]  \neq 0 $, for $ \e \g^{-1} $ small enough, but 
the average $ \Pi_C {\cal L} W_2 [1 / r^\infty_0]  = 0 $ by \eqref{mL mappa tutto in media nulla}.
In the reversible case $ r^\infty_0 = 0 $ was proved
in remark \ref{r0=0}. 
\end{pf}

As a consequence of \eqref{zero eigenvalue}, the definition of $ \L_\infty^{2 \g}  $ in 
\eqref{Omegainfty} 
(just specializing \eqref{Omegainfty} with $ k = 0 $), and \eqref{omdio} (with $\g$ and $\t$ as in \eqref{gamma-tau}), we deduce   
also the {\it first} order Melnikov non-resonance conditions
\be\label{cantor inversione linearizzato}
\forall \l \in \L_{\infty}^{2 \g} \, , \qquad  
\big|\ii \l\bar\o \cdot l + \mu_j^\infty (\l) \big| \geq 2 \g \frac{ \langle j \rangle^3}{ \langle l \rangle^\t}, \quad \forall (l, j) \neq (0, 0)  \, .
\ee

\begin{lemma} {\bf (Invertibility of ${\cal L}_\infty $)} \label{stima-L-infinito}
For all $ \l \in \L_\infty^{2 \g} (u) $, for all $ g \in H^s_{00} $ the equation $ {\cal L}_\infty w = g $ has the unique
solution with zero average 
\be\label{Linfty inverse}
\mL_{\infty}\inv \, g (\ph,x)
:= \sum_{(l,j) \neq (0,0)} \frac{g_{lj}}{\ii \l\bar\o \cdot l + \mu_j^\infty (\l) }\, e^{\ii(l \cdot \ph + j x)}.
\ee
For all Lipschitz family $ g := g(\l) \in H^s_{00} $  we have
\begin{equation}\label{tame-L-infinito}
\left\|{\cal L}_{\infty}^{-1} g \right\|_{s}^{\Lipg} \leq C \gamma^{-1} \left\| g \right\|_{s + 2\tau + 1}^{\Lipg} 	\, . 
\end{equation}
In the reversible case, if $ g \in Y $ then $ {\cal L}_\infty^{-1} g \in X $. 
\end{lemma}

\begin{pf}
For all $\l \in \L_\infty^{2\g} (u) $, by \eqref{cantor inversione linearizzato}, formula
 \eqref{Linfty inverse} is well defined and 
  \begin{equation}\label{stima-0-inverso-L-infinito}
\left\|{\cal L}_{\infty}^{-1}(\l)g(\l)\right\|_{s} \lessdot  \gamma^{-1}\left\|g(\l)\right\|_{s + \tau}\, .
\end{equation}
Now we prove the Lipschitz estimate. For $ \l_1 , \l_2 \in \L_\infty^{2\g} (u) $
\begin{equation}\label{delta-L-infinito}
{\cal L}_{\infty}^{-1}(\l_{1})g(\l_{1}) - {\cal L}_{\infty}^{-1}(\l_{2})g(\l_{2}) = {\cal L}_{\infty}^{-1}(\l_{1}) [g(\l_1) - g(\l_2)] + 
\big({\cal L}_{\infty}^{-1}(\l_1) - {\cal L}_{\infty}^{-1}(\l_2) \big)g(\l_{2})\,.
\end{equation}
By \eqref{stima-0-inverso-L-infinito}
\begin{equation}\label{primo-pezzo-delta-L-infinito}
\g \|{\cal L}_{\infty}^{-1}(\l_{1}) [g(\l_{1}) - g(\l_{2})] \|_s 
\lessdot  \| g(\l_{1})- g(\l_{2}) \|_{s + \tau} \leq \g^{-1} \Vert g \Vert_{s + \tau}^{\Lipg} |\l_1 - \l_2 | \, . 
\end{equation}
Now we estimate the second term of \eqref{delta-L-infinito}.
We simplify notations writing $ g := g(\l_{2}) $ and $ \d_{lj} := \ii \l \bar \om \cdot l + \mu_j^\infty $.
\begin{equation}\label{secondo-pezzo-delta-L-infinito}
\big({\cal L}_{\infty}^{-1}(\l_{1}) - {\cal L}_{\infty}^{-1}(\l_{2})\big)g = 
\sum_{(l , j)\neq(0,0)}
\frac{\delta_{lj}(\l_{2}) - \delta_{lj}(\l_{1})}{\delta_{lj}(\l_{1})\delta_{lj}(\l_{2})} \, g_{lj} 
e^{\ii(l \cdot \ph + j x)} \, .
\end{equation}
The bound  
\eqref{autofinali} imply
$ \vert \mu_{j}^{\infty} \vert^{\rm lip} \lessdot \e \g^{-1} | j |^{3}  \lessdot | j |^{3} $ and, using also 
\eqref{cantor inversione linearizzato}, 
\begin{eqnarray}
\g\frac{|\delta_{lj}(\l_{2}) - \delta_{lj}(\l_{1}) |}{|\delta_{lj}(\l_{1})| |\delta_{lj}(\l_{2}) |} 
\! \! & \lessdot & \! \!  \frac{( | l | +  | j |^{3}) \langle  l \rangle^{2\tau}}{\gamma \langle j \rangle^{6}} |\l_{2} - \l_{1} | 
\lessdot  \langle l \rangle^{2\tau + 1} \gamma^{-1} | \l_2 - \l_1 | \, . \label{deltalj-o1-o2}
\end{eqnarray}
Then \eqref{secondo-pezzo-delta-L-infinito} and \eqref{deltalj-o1-o2} imply 
$ \g \| ({\cal L}_{\infty}^{-1}(\l_2) - {\cal L}_{\infty}^{-1}(\l_1) )g \|_s 
\lessdot \gamma^{-1} \|g \|_{s + 2\tau + 1}^{\Lipg} |\l_2 - \l_1 | $
that, finally, with  \eqref{stima-0-inverso-L-infinito}, \eqref{primo-pezzo-delta-L-infinito}, 
prove  \eqref{tame-L-infinito}.   The last statement follows by  the property \eqref{reversibilitˆ autovalori finali}.
\end{pf}

In order to solve the equation $ {\cal L} h = f $ we first prove the following lemma.

\begin{lemma} \label{lemma:iso zero mean}
Let $\mathfrak s_0 + \t + 3 \leq s \leq q - \s - \b - 3 $.  Under the assumption \eqref{f = der g} we have 
\be \label{iso W1 Hs0}
W_1 (H^s_{00}) = H^s_{00} \, ,  \quad \ W_1^{-1} (H^s_{00}) = H^s_{00} \, . 
\ee
\end{lemma}

\begin{pf} 
It is sufficient to prove that $ W_1 (H^s_{00}) = H^s_{00} $ because the second equality of \eqref{iso W1 Hs0} follows applying the isomorphism $  W_1^{-1} $. Let us give the proof of  the inclusion 
\begin{equation} \label{inclusion}
W_1 (H^s_{00}) \subseteq H^s_{00}
\end{equation}
(which is essentially algebraic).  For any $ g \in H^s_{00}$, 
let $ w(\ph,x) := {\cal L}_\infty^{-1} g \in H^{s - \t}_{00} $ defined in \eqref{Linfty inverse}.
Then $ h := W_2 w \in H^{s-\t} $ 
satisfies
\[
\mL h 
\stackrel{\eqref{L-coniugato}} = W_1 \mL_\infty W_2\inv h
= W_1 \mL_\infty w 
= W_1 g \, .
\]
By \eqref{mL mappa tutto in media nulla} we deduce 
that $W_1 g = \mL h  \in H^{s - \t - 3}_{00} $. 
Since  $ W_1 g \in H^s $ by Lemma \ref{stime-tame-coniugio}, we conclude 
$ W_1 g \in H^s \cap H^{s - \t - 3}_{00} = H^s_{00}$.   The proof of \eqref{inclusion} is complete.
\smallskip

It remains to prove that $H^s_{00} \setminus W_1(H^s_{00}) = \emptyset$. 
By contradiction, let $ f \in H^s_{00} \setminus W_1(H^s_{00}) $. 
Let $  g := W_1\inv f \in H^s $  by Lemma \ref{stime-tame-coniugio}. 
Since $W_1 g = f \notin W_1(H^s_{00})$, it follows that $ g \notin H^s_{00} $
(otherwise it contradicts \eqref{inclusion}), namely $c := \Pi_C g \neq 0$.
Decomposing $ g = c + \pp g $ (recall \eqref{proiez00}) and applying $W_1$, we get
$ W_1 g = c W_1[1] + W_1 \pp g $. 
Hence 
\[
W_1[1] = c\inv (W_1 g - W_1 \pp g) \in H^s_{00}
\]
because $W_1 g = f \in H^s_{00}$ and $W_1 \pp g \in W_1(H^s_{00}) \subseteq H^s_{00}$ by \eqref{inclusion}. 
However, $ \Pi_C W_1[1] \neq 0 $, a contradiction.
\end{pf}

\begin{remark} 
In the Hamiltonian case (which always satisfies \eqref{f = der g}), 
the $ W_i (\vphi ) $ are maps  of  (a subspace of)  $ H^1_0 $ 
so that  Lemma \ref{lemma:iso zero mean} is automatic, and there is no need
of Lemma  \ref{lem:modo0}.  
\end{remark}

We may now prove the main result  of sections \ref{sec:regu} and \ref{sec:redu}.

\begin{theorem}\label{inversione linearizzato} {\bf (Right inverse of $ {\cal L}$)} 
Let 
\be\label{perdite stima tame linearizzato}
\tau_1 := 2\t + 7, \quad \mu:= 4\t + \s + \b + 14\,,
\ee
where $ \s $, $\b $ are  defined in \eqref{costanti lemma mostro}, \eqref{defbq} respectively. 
Let $ u ( \l ) $, $ \l \in \L_o \subseteq \L $, be a Lipschitz family with 
\be\label{verasmall}
\Vert u \Vert_{\mathfrak s_0 + \mu}^{\Lipg} \leq 1 \, . 
\ee 
Then there exists $ \delta $ (depending on the data of the problem) such that if 
$$ 
\e\g^{-1} \leq \delta \, ,
$$
and condition \eqref{f = der g}, resp. the reversibility condition \eqref{parity f}, holds, 
then for all $ \l \in \L_\infty^{2 \g}(u)$ defined in  \eqref{Omegainfty},  the linearized operator $\mL:= \mL(\l, u(\l))$ 
(see \eqref{mL}) admits a right inverse on $ H^s_{00} $, resp. $ Y \cap H^s $. 
More precisely, for $\mathfrak  s_0 \leq s \leq q - \mu$, for all Lipschitz family $ f(\l) \in H^s_{00} $, resp. $ Y \cap H^s $, 
the function  
\be\label{scelta}
h := {\cal L}^{-1} f  := W_2 \mL_{\infty }\inv \, W_1\inv f 
\ee
is a solution of $ {\cal L} h = f $. In the reversible case,   
$ {\cal L}^{-1} f \in X $. 
Moreover 
\be\label{tame-L}
\Vert \mL^{-1} f \Vert_{s}^{\Lipg} \leq C(s)\g^{-1} \Big(\Vert f \Vert_{s + \t_1}^{\Lipg} + \Vert u \Vert_{s + \mu}^{\Lipg} 
\Vert f \Vert_{\mathfrak s_0}^{\Lipg} \Big) \, . 
\ee
\end{theorem}

\begin{pf} 
Given $f \in H^s_{00}$, resp. $ f \in Y \cap H^s $, with $ s $ like in Lemma \ref{lemma:iso zero mean}, the equation
$ \mL h = f $ 
can be solved for $ h $ because $\Pi_C f = 0 $. 
Indeed, by \eqref{L-coniugato}, the equation $ {\cal L} h = f $ is equivalent to 
$ \mL_\infty W_2\inv h = W_1\inv f $
where $W_1\inv f \in H^s_{00} $  by Lemma \ref{lemma:iso zero mean}, resp. $ W_1\inv f \in Y \cap H^s  $
being $ W_1^{-1} $  reversibility-preserving (Lemma \ref{stime-tame-coniugio}).
As a consequence, by Lemma \ref{stima-L-infinito}, all the solutions of $ {\cal L} h = f $ 
are 
\be\label{tutte le soluzioni}
h = c W_2[1] + W_2 \mL_{\infty}\inv W_1\inv f, \quad c \in \R \, . 
\ee 
The solution \eqref{scelta} is the one 
with $ c = 0 $. In the reversible case, 
the fact that $ {\cal L}^{-1} f \in X $ follows by \eqref{scelta} and the fact that $ W_i $, $ W_i^{-1}$ are reversibility-preserving
and $ {\cal L}_\infty^{-1} : Y \to X $, see Lemma \ref{stima-L-infinito}. 

Finally  
\eqref{tame-Phi1Phi2}, \eqref{tame-L-infinito}, \eqref{verasmall} imply
$$
\Vert \mL^{-1} f \Vert_s^{\Lipg} \leq C(s)\g^{-1} \big( \Vert f \Vert_{s + 2\t + 7}^{\Lipg} + 
\Vert u \Vert_{s + 2\t + \s + \b + 7}^{\Lipg} \Vert f \Vert_{\mathfrak s_0 + 2\t + 7}^{\Lipg} \big) 
$$
and \eqref{tame-L} follows using \eqref{interpolation estremi fine} 
with $ b_0 = \mathfrak s_0 $, $ a_0 :=  \mathfrak s_0 + 2 \t + \s + \b + 7 $, 
$ q = 2 \t + 7 $, $ p = s -  \mathfrak s_0 $.
\end{pf}

In the next section we  apply Theorem \ref{inversione linearizzato}  to deduce
 tame estimates for the inverse linearized operators 
at any step of the Nash-Moser  scheme. 
The approximate solutions along the iteration will 
satisfy \eqref{verasmall}.

\section{The Nash-Moser iteration} \label{sec:NM} 

We define the finite-dimensional subspaces of trigonometric polynomials
$$
H_{n} := \Big\{ u \in L^{2}(\T^{\nu + 1}) : u(\vphi,x)=\sum_{\left|(l , j)\right|\leq N_{n}}u_{lj} e^{\ii(l\cdot\vphi + j x)} \Big\} 
$$
where $ N_n := N_0^{\chi^n}$ (see  \eqref{defN}) and the corresponding orthogonal projectors
$$
\Pi_{n}:=\Pi_{N_{n}} : L^{2}(\T^{\nu + 1}) \rightarrow H_{n} \, , \quad \Pi_{n}^\bot := I - \Pi_{n} \, .
$$
The following smoothing properties hold: for all $\alpha , s \geq 0$, 
\begin{equation}\label{smoothing-u1}
\|\Pi_{n}u \|_{s + \alpha}^\Lipg 
\leq N_{n}^{\alpha} \| u \|_{s}^\Lipg, 
\ \ \forall u(\lm) \in H^{s} \,; 
\quad
\|\Pi_{n}^\bot u \|_{s}^\Lipg 
\leq N_{n}^{-\alpha} \|u \|_{s + \alpha}^\Lipg, 
\ \ \forall u(\lm) \in H^{s + \alpha},
\end{equation}
where the function $u(\lm)$ depends on the parameter $\lm $ in a Lipschitz way.  
The bounds \eqref{smoothing-u1} are the 
classical smoothing estimates for truncated Fourier series, 
which also hold  
with the norm $\| \cdot \|^\Lipg_s $ defined in \eqref{def norma Lipg}.

Let
\be\label{operatorF(u)}
F(u) := F(\l, u) := \l \bar \o \cdot \partial_{\vphi} u + u_{xxx} 
+ \e f(\vphi , x , u, u_{x}, u_{xx}, u_{xxx} ) \, .
\ee
We define the constants 
\begin{equation}\label{sigma-beta-kappa}
\kappa := 28 + 6 \mu, \qquad 
\b_1 := 50 + 11 \mu, \, 
\end{equation}
where $\mu$ is the loss of regularity in \eqref{perdite stima tame linearizzato}. 

\begin{theorem}\label{iterazione-non-lineare} 
{\bf (Nash-Moser)} Assume that $ f \in C^q $, $ q \geq \mathfrak s_0 + \mu + \b_1  $, 
satisfies the assumptions of Theorem \ref{thm:main}
or  Theorem \ref{thm:mainrev}. 
Let 
$ 0 < \gamma \leq {\rm min}\{\g_0, 1/48 \} $, $ \tau > \nu + 1 $. 
Then there exist $ \d  > 0 $, $ C_* > 0  $, $  N_0 \in \N $ (that may depend also on $ \tau $) such that, if
$ \e \g^{-1} < \d $, then, for all $ n \geq 0 $: 
\begin{itemize}
\item[$({\cal P}1)_{n}$] 
there exists  a function 
$u_n : \mG_n \subseteq \Lambda \to H_n$, $\lm \mapsto u_n(\lm)$,  
with $ \| u_{n} \|_{\mathfrak s_0 + \mu}^{\Lipg} \leq 1 $,  $ u_0 := 0$, 
where ${\cal G}_{n} $ are Cantor like subsets of $ \Lambda := [1/2, 3/2] $ defined inductively by: 
$ {\cal G}_{0} := \L $, 
\begin{eqnarray}\label{G-n+1}
{\cal G}_{n+1} & := & 
\Big\{ \l  \in {\cal G}_{n} \, : \, |\ii \o \cdot l + \mu_j^\infty (u_{n}) -
\mu_k^\infty (u_{n})| \geq \frac{2\gamma_{n} |j^{3}-k^{3}|}{\left\langle l\right\rangle^{\tau}}\, , \ \forall j , k \in \Z, \ l \in \Z^{\nu} \Big\}
\end{eqnarray}
where $ \gamma_{n}:=\gamma (1 + 2^{-n}) $. In the reversible case, 
namely \eqref{parity f} holds, then $ u_n (\l ) \in X $. 

The difference $h_n := u_{n} - u_{n-1}$, where, for convenience, $h_0 := 0$, 
satisfy
\be  \label{hn}
\| h_{n} \|_{\mathfrak s_0 + \mu}^\Lipg \leq C_* \e \gamma^{-1} N_{n}^{-\s_1} \,,  \quad \s_1 := 18  + 2 \mu \, . 
\ee
\item[$({\cal P}2)_{n}$]   $ \| F(u_n) \|_{\mathfrak s_{0}}^{\Lipg} \leq C_* \e N_{n}^{- \kappa}$.

\item[$({\cal P}3)_{n}$] \emph{(High norms).} $ \|u_n \|_{\mathfrak s_{0}+ \beta_1}^{\Lipg} 
\leq C_* \e \gamma^{-1} N_{n}^{\kappa}$
and  $ \|F(u_n ) \|_{\mathfrak s_{0}+\beta_1}^{\Lipg} \leq C_*  \e N_{n}^{\kappa}$.

\item[$({\cal P}4)_{n}$] \emph{(Measure).} 
\newcommand{\mB}{\mathcal{B}}
The measure of the Cantor like sets satisfy
\be\label{Gmeasure}
|{\cal G}_0 \setminus {\cal G}_1| \leq C_* \g \, , \quad  \big| {\cal G}_n \setminus {\cal G}_{n+1} \big| 
\leq \g C_* N_{n}^{-1} \, , \ n \geq 1 .
\ee

\end{itemize}
All the Lip norms are defined on $ {\cal G}_{n} $.
\end{theorem}

\begin{pf} 
The proof of Theorem \ref{iterazione-non-lineare} is split into several steps. For simplicity, we denote $ \| \ \|^{\rm Lip}  $ by $ \| \ \|  $.

\smallskip

\noindent
{\sc Step 1:} \emph{prove $(\mP1,2,3)_0$.} $(\mP1)_0$ and the first inequality of $(\mP3)_0$ are trivial because 
$u_0 = h_0 = 0$. $(\mP2)_0$ and the second inequality of $(\mP3)_0$ follow with
$ C_*  \geq $ $ \max\{  \| f(0)\|_{\mathfrak s_0} N_0^\kappa, $ $ \| f(0)\|_{\mathfrak s_0+ \b_1} N_0^{-\kappa}  \} $.

\smallskip

\noindent
{\sc Step 2:} \emph{assume that $(\mP1,2,3)_n$ hold for some $n \geq 0$, and prove $(\mP1,2,3)_{n+1}$.} 
By $(\mP1)_n$ we know that $ \| u_n \|_{\mathfrak s_{0} + \mu} \leq 1 $, namely condition
\eqref{verasmall} is satisfied. Hence, for $ \e \g^{-1}$ small enough,
Theorem \ref{inversione linearizzato} applies. Then, for all $\l \in {\cal G}_{n+1} $ defined in \eqref{G-n+1}, 
the linearized operator 
\[ 
\mL_n(\lm) := {\cal L}(\lm, u_{n}(\lm)) = F'(\lm, u_n(\lm))
\] 
(see \eqref{mL})
admits a right inverse for all $ h \in H^s_{00} $,  if condition \eqref{f = der g} holds, 
respectively for  $ h \in  Y \cap H^s $  if the reversibility condition \eqref{parity f} holds. 
Moreover \eqref{tame-L} gives the estimates
\begin{align} \label{L-1alta}
\| {\cal L}_n^{-1} h \|_s 
& \leq_s \g^{-1} \Big( \| h \|_{s+\t_1} + \| u_n \|_{s+ \mu} \|h \|_{\mathfrak s_0}  \Big) \, ,  
\quad \forall h(\lm), 
\\
\label{L-1s0}
\| {\cal L}_n^{-1} h \|_{\mathfrak s_0} 
& \leq \g^{-1} N_{n+1}^{\t_1}  \| h \|_{\mathfrak s_0}  \, , 
\quad \forall h(\lm) \in H_{n+1} \,,
\end{align}
(use 
\eqref{smoothing-u1} and $ \| u_n \|_{\mathfrak s_{0} + \mu} \leq 1 $), 
for all Lipschitz map $h(\lm)$. 
Then, for all $\l \in {\cal G}_{n+1} $,  we define 
\begin{equation}\label{soluzioni-approssimate}
u_{n+1} := u_{n} + h_{n + 1} \in H_{n+1} \, , \quad 
h_{n + 1}:= - \Pi_{n + 1} {\cal L}_n^{-1} \Pi_{n + 1} F(u_{n}) \, ,
\end{equation}
which is well defined because, if condition
\eqref{f = der g} holds then $  \Pi_{n + 1} F(u_n) \in H^s_{00} $, and, respectively, 
if \eqref{parity f} holds, then $  \Pi_{n + 1} F(u_{n}) \in Y \cap H^s $
(hence in both cases $  {\cal L}_n^{-1} \Pi_{n + 1} F(u_n) $ exists). Note also that
in the reversible case $ h_{n + 1} \in X $ and so $  u_{n + 1} \in X  $. 

Recalling \eqref{operatorF(u)} and that $\mL_n := F'(u_n) $, we write
\be\label{FTay}
F(u_{n + 1}) =  F(u_{n}) + {\cal L}_n h_{n + 1} + \e Q(u_{n}, h_{n + 1})
\ee
where 
$$
Q(u_{n},h_{n + 1}) := {\cal N}(u_{n} + h_{n + 1}) - {\cal N}(u_{n}) - {\cal N}'(u_{n}) h_{n + 1}, 
\quad 
\mN(u) := f(\ph,x,u,u_x, u_{xx}, u_{xxx}).
$$
With this definition, 
\[
F(u) = L_\om u + \e {\cal N}(u), 
\quad 
F'(u) h = L_\om h + \e \mN'(u)h, 
\quad 
L_\om := \ompaph + \pa_{xxx}.
\]
By \eqref{FTay} and \eqref{soluzioni-approssimate} we have
\begin{eqnarray}
F(u_{n + 1}) 
& = & 
F(u_{n}) - \mL_n \Pi_{n + 1} \mL_n^{-1} \Pi_{n + 1} F(u_{n}) + \e Q(u_{n},h_{n + 1}) \nonumber\\
& = & 
\Pi_{n + 1}^{\bot} F(u_{n}) + \mL_n \Pi_{n + 1}^{\bot} \mL_n^{-1} \Pi_{n + 1} F(u_{n}) 
+ \e Q(u_{n},h_{n + 1}) \nonumber\\
& = & \Pi_{n + 1}^{\bot} F(u_{n}) + \Pi_{n + 1}^{\bot} \mL_n  \mL_n^{-1} \Pi_{n + 1} F(u_{n}) 
+ [ \mL_n , \Pi_{n + 1}^{\bot} ] \mL_n^{-1} \Pi_{n + 1} F(u_{n}) 
+ \e Q(u_{n},h_{n + 1})\nonumber \\
& = & \Pi_{n + 1}^{\bot} F(u_{n}) + 
\e [{\cal N}'(u_{n}) , \Pi_{n + 1}^{\bot} ] \mL_n^{-1} \Pi_{n + 1} F(u_{n}) 
+ \e Q(u_{n},h_{n + 1})
\label{F(u-n+1)}
\end{eqnarray}
where we have gained an extra $\e$ from the commutator 
$$
[{\cal L}_n , \Pi_{n + 1}^{\bot} ] 
= [ L_{\o} + \e{\cal N}'(u_{n}) , \Pi_{n + 1}^{\bot} ] 
= \e [{\cal N}'(u_{n}) , \Pi_{n + 1}^{\bot} ] \,. 
$$
\begin{lemma}
Set
\begin{equation}\label{BnB'n}
U_{n}:=\Vert u_{n} \Vert_{\mathfrak s_{0}+\beta_1} + \g^{-1} \Vert F(u_{n}) \Vert_{\mathfrak s_{0}+\beta_1} \,, 
\qquad  
w_n := \g^{-1} \Vert F(u_{n}) \Vert_{\mathfrak s_{0}} \, . 
\end{equation}
There exists $C_0 := C ( \t_1, \mu, \nu, \b_1) >  0 $ such that 
\be\label{quadratico}
w_{n+1} 
\leq C_0 N_{n + 1}^{- \b_1 + \mu'} U_n ( 1 +  w_n ) 
+ C_0  N_{n + 1}^{6 + 2\mu}  w_n^2 ,  
\qquad 
U_{n+1} 
\leq C_0 N_{n + 1}^{9 + 2\mu} ( 1 + w_n )^2 \, U_n \, . 
\ee
\end{lemma}

\begin{pf}
The operators $\mN'(u_n)$ and $Q(u_n,\cdot)$ satisfy the following tame estimates: 
\begin{align} \label{tameQ}
\| Q(u_n , h) \|_s  
& \leq_s  \| h \|_{\mathfrak s_0 + 3}  \Big( \| h \|_{s + 3} 
+ \| u_n \|_{s + 3}  \| h \|_{\mathfrak s_0+3} \Big)
\quad \ \forall h(\lm), 
\\
\label{Qs0}
\| Q(u_n , h) \|_{\mathfrak s_0}  
& \leq N_{n+1}^6   \| h \|_{\mathfrak s_0}^2 \ 
\quad \forall h(\lm) \in H_{n+1} , 
\\
\label{tameN'}
\| \mN'(u_n) h \|_{s} 
& \leq_{s}  \| h \|_{s + 3} + \| u_n \|_{s + 3}  \| h \|_{\mathfrak s_0+3} 
\quad \forall h(\lm),  
\end{align}
where $h(\lm)$ depends on the parameter $\lm$ in a Lipschitz way. 
The bounds \eqref{tameQ} and \eqref{tameN'} follow by 
\ref{lemma:composition of functions, Moser}$(i)$ 
and Lemma \ref{lemma:Lip generale}.
\eqref{Qs0} is simply \eqref{tameQ} at $s = \mathfrak s_0$, using that $\| u_n \|_{\mathfrak s_0 + 3} \leq 1$, 
$u_n, h_{n+1} \in H_{n+1}$ and the smoothing \eqref{smoothing-u1}. 

By \eqref{L-1alta} and \eqref{tameN'}, the term (in \eqref{F(u-n+1)})
$
R_n := [ {\cal N}' (u_n), \Pi_{n+1}^\bot ] {\cal L}_n^{-1} \Pi_{n+1} F(u_n) 
$
satisfies, using also that $  u_n \in H_n $ and \eqref{smoothing-u1}, 
\begin{align} \label{Rn alta}
\| R_n \|_s 
& \leq_s \g^{-1} N_{n+1}^{\mu'} 
\Big( \| F(u_n) \|_s + \| u_n \|_{s} \| F(u_n) \|_{\mathfrak s_0}  \Big), 
\quad \mu' := 3 + \mu, 
\\
\label{RNS0}
\| R_n \|_{\mathfrak s_0} 
& \leq_{\mathfrak s_0 + \b_1} \g^{-1} N_{n+1}^{-\b_1 + \mu'} \Big(\| F(u_n) \|_{\mathfrak s_0+ \b_1} + 
\| u_n \|_{\mathfrak s_0 + \b_1} \| F(u_n) \|_{\mathfrak s_0}  \Big),
\end{align}
because $\mu \geq \t_1 + 3$. 
In proving \eqref{Rn alta} and \eqref{RNS0}, we have simply estimated $\mN'(u_n) \Pi_{n+1}^\perp$ and 
$\Pi_{n+1}^\perp \mN'(u_n)$ separately, without using the commutator structure.

From the definition \eqref{soluzioni-approssimate} of $h_{n+1}$, 
using \eqref{L-1alta}, \eqref{L-1s0} and \eqref{smoothing-u1}, we get
\begin{align} \label{h-n+1-alta}
\| h_{n + 1} \|_{\mathfrak s_{0}+ \beta_1} 
& \leq_{\mathfrak s_0 + \beta_1} \gamma^{-1}  N_{n + 1}^{\mu} \Big( \|F(u_{n}) \|_{\mathfrak s_0+\beta_1} + 
\| u_n \|_{\mathfrak s_0 + \beta_1} \| F(u_{n}) \|_{\mathfrak s_0} \Big), 
\\
\label{h-n+1-bassa}
\|h_{n + 1} \|_{\mathfrak s_{0}} 
& \leq_{\mathfrak s_0} \g^{-1} N_{n + 1}^{\mu} \| F(u_n) \|_{\mathfrak s_0}
\end{align}
because $\mu \geq \t_1$. Then 
\begin{align}
\|u_{n + 1} \|_{\mathfrak s_0 + \beta_1} 
& \stackrel{\eqref{soluzioni-approssimate}} \leq \|u_{n} \|_{\mathfrak s_0 + \beta_1} + \| h_{n + 1} \|_{\mathfrak s_0 + \beta_1} \notag \\
& \stackrel{\eqref{h-n+1-alta}} \leq_{\mathfrak s_0 + \b_1} \| u_n \|_{\mathfrak s_0 + \b_1} \Big( 1 + \g^{-1} N_{n+1}^{\mu} \| F (u_n) \|_{\mathfrak s_0} \Big) 
+ \gamma^{-1} N_{n + 1}^{\mu} \| F(u_n) \|_{\mathfrak s_0 + \beta_1}. 
\label{u-n+1-alta}
\end{align}
Formula \eqref{F(u-n+1)} for $F(u_{n+1})$, and  
\eqref{RNS0}, \eqref{Qs0}, \eqref{h-n+1-bassa}, $\e \g\inv \leq 1$, \eqref{smoothing-u1}, imply
\be\label{stima-induttiva-bassa}
\| F(u_{n + 1}) \|_{\mathfrak s_0} \leq_{\mathfrak s_0 + \b_1}  
N_{n + 1}^{- \b_1 + \mu'} 
\Big( \| F(u_n) \|_{\mathfrak s_0 + \beta_1} + \| u_n \|_{\mathfrak s_0 + \beta_1} \| F(u_n) \|_{\mathfrak s_0} \Big)  
+ \e \g^{-2} N_{n + 1}^{6 + 2 \mu} \| F(u_{n}) \|_{\mathfrak s_{0}}^{2}.
\ee
Similarly, using the ``high norm'' estimates 
\eqref{Rn alta}, \eqref{tameQ}, \eqref{h-n+1-alta}, \eqref{h-n+1-bassa}, $\e \g\inv \leq 1$ and \eqref{smoothing-u1}, 
\begin{equation}\label{stima-induttiva-alta}
\| F(u_{n + 1}) \|_{\mathfrak s_0 + \beta_1} 
\leq_{\mathfrak s_0 + \b_1} 
\Big( \| F(u_n) \|_{\mathfrak s_0 + \b_1} + \| u_n \|_{\mathfrak s_0 + \b_1} \| F(u_n) \|_{\mathfrak s_0} \Big)
\Big( 1 + N_{n+1}^{\mu'} + N_{n+1}^{9 + 2 \mu} \g\inv \| F(u_n) \|_{\mathfrak s_0} \Big).
\end{equation}
By \eqref{u-n+1-alta}, \eqref{stima-induttiva-bassa} and \eqref{stima-induttiva-alta} we deduce 
\eqref{quadratico}. 
\end{pf}

By $ ({\cal P}2)_n $ 
we deduce, for $ \e \g^{-1} $ small, that (recall the definition on $ w_n $ in \eqref{BnB'n})
\be \label{control wU}
w_n \leq  \e \g^{-1} C_* N_{n}^{-\kappa} \leq 1, 
\ee
Then, by the second inequality in \eqref{quadratico}, \eqref{control wU},  $ ({\cal P}3)_n $ 
(recall the definition on $ U_n $ in \eqref{BnB'n}) and the choice of $ \kappa $ in \eqref{sigma-beta-kappa},
we deduce 
$ U_{n+1} \leq C_* \e \g^{-1} N_{n+1}^\kappa $, for 
$ N_0 $ large enough. This proves $ ({\cal P}3)_{n+1} $. 

Next, by the first inequality in \eqref{quadratico}, \eqref{control wU},  $ ({\cal P}2)_n $ 
(recall the definition on $ w_n $ in \eqref{BnB'n}) and  \eqref{sigma-beta-kappa},
we deduce 
$ w_{n+1} \leq C_* \e \g^{-1} N_{n+1}^\kappa $, for 
$ N_0 $ large, $ \e \g^{-1}$ small. This proves $ ({\cal P}2)_{n+1} $.

The bound \eqref{hn} at the step $ n +1$ follows by \eqref{h-n+1-bassa} and $({\cal P}2)_n $ (and 
\eqref{sigma-beta-kappa}).
Then 
$$
\| u_{n+1} \|_{\mathfrak s_0 + \mu} \leq \| u_0 \|_{\mathfrak s_0+ \mu} + \sum_{k=1}^{n+1} \| h_k \|_{\mathfrak s_0 + \mu} 
 \leq \sum_{k=1}^\infty C_* \e \g^{-1} N_k^{-\s_1} \leq 1 
$$
for $ \e \g^{-1} $ small enough. 
As a consequence $(\mP1,2,3)_{n+1}$ hold.

\smallskip

\noindent
{\sc Step 3:} \emph{prove $(\mP4)_n$, $n \geq 0$.} 
For all $n \geq 0$, 
\be\label{natale}
{\cal G}_n \setminus{\cal G}_{n+1} = \bigcup_{l \in \Z^{\nu}, j,k \in \Z} R_{ljk} (u_{n})   
\ee
where 
\begin{eqnarray}
R_{ljk} (u_{n}) 
& := &
\left\{\lambda\in{\cal G}_n : \left|\ii \lambda\bar\o \cdot l + \mu_{j}^{\infty}(\lambda,u_{n}(\lambda)) - \mu_{k}^{\infty}(\lambda,u_{n}(\lambda))\right| < 2\gamma_{n} | j^{3}-k^{3} | \left\langle l\right\rangle^{-\tau}\right\} \, . 
\label{R-ljk(u-n)}
\end{eqnarray}
Notice that, by the definition \eqref{R-ljk(u-n)}, $R_{ljk} (u_{n}) = \emptyset$ for $j = k$. Then we can suppose in the 
sequel that $j \neq k$.
We divide the estimate 
 into some lemmata.

\begin{lemma}\label{risonanti-1}
For  $ \e  \g^{-1}$ small enough,  
for all $n \geq 0$, $|l|\leq N_n$, 
\begin{equation}\label{inclusione-1}
R_{ljk}(u_{n}) \subseteq R_{ljk}(u_{n - 1}).
\end{equation}
\end{lemma}

\begin{pf}
We claim that, for all $ j , k \in \Z $,  
\begin{equation}\label{marco}
|(\mu_{j}^{\infty} - \mu_{k}^{\infty})(u_{n}) 
- (\mu_{j}^{\infty} - \mu_{k}^{\infty})(u_{n-1})| 
\leq C \e |j^{3} - k^{3}| N_n^{-\a} \, , \quad \forall \l \in {\cal G}_n \, , 
\end{equation}
where $\mu_{j}^{\infty}(u_{n}) := \mu_{j}^{\infty}(\lambda, u_{n}(\lambda))$ and $ \a $ is defined in \eqref{alpha-beta}. Before proving
\eqref{marco} we show how it implies \eqref{inclusione-1}. 
For all $ j \neq k$, $|l| \leq N_{n}$,  $\lm \in \mG_n$, by \eqref{marco}  
\begin{align*}
|\ii \l \bar\o \cdot l + \mu_{j}^{\infty}(u_{n}) - \mu_{k}^{\infty}(u_{n}) | 
& \geq 
|\ii \l \bar\o \cdot l + \mu_{j}^{\infty}(u_{n-1}) - \mu_{k}^{\infty}(u_{n-1}) |  
- |(\mu_{j}^{\infty}- \mu_{k}^{\infty})(u_{n}) - (\mu_{j}^{\infty}- \mu_{k}^{\infty})(u_{n - 1}) | 
\\ 
& \geq  
2\gamma_{n - 1} |j^3 - k^3| \langle l \rangle^{-\tau} - C \e  |j^3 - k^3|  N_{n}^{-\a}   
\geq   2\gamma_{n} |j^3 - k^3|  \langle l \rangle^{-\tau} 
\end{align*}
for 
$C \e \g\inv  N_{n}^{\t -\a}\,  2^{n+1} \leq 1$
(recall that $\gamma_n := \g (1 + 2^{-n})$), which implies  \eqref{inclusione-1}.
\\[1mm]
{\sc Proof of \eqref{marco}.} 
By \eqref{espressione autovalori}, 
\begin{align}
(\mu_{j}^{\infty}- \mu_{k}^{\infty})(u_{n}) - (\mu_{j}^{\infty}- \mu_{k}^{\infty})(u_{n - 1})
& = 
-\ii \big[ m_3(u_{n}) - m_{3}(u_{n - 1}) \big] (j^{3} - k^{3}) 
+\ii \big[ m_1(u_{n}) - m_1(u_{n-1}) \big]  (j - k)  \nonumber
\\
& \quad 
+ r_{j}^{\infty}(u_{n}) - r_{j}^{\infty}(u_{n - 1}) 
- \big( r_{k}^{\infty}(u_{n}) - r_{k}^{\infty}(u_{n-1}) \big) \label{vici}
\end{align}
where $ m_3 (u_{n}) := m_3(\lambda, u_{n}(\lambda))$ and similarly for 
$ m_1, r_{j}^{\infty}$.
We first apply Theorem \ref{thm:abstract linear reducibility}-${\bf (S4)_{\nu}}$ 
with $ \nu = n + 1 $, $ \g = \g_{n-1} $, $ \g - \rho = \g_n $, and $ u_1 $, $ u_2 $, replaced, respectively,   
by $ u_{n-1} $, $ u_n $, 
 in order to conclude that
\be\label{primoste}
\L_{n+1}^{\g_{n-1}} ( u_{n-1}) \subseteq \L_{n+1}^{\g_n} ( u_n ) \, . 
\ee
The smallness condition in \eqref{legno} is satisfied because 
$\s_2 < \mu$  
(see definitions \eqref{alpha-beta},  
\eqref{perdite stima tame linearizzato})  and so
$$
\e C N_n^{\tau} \Vert u_n - u_{n - 1} \Vert_{\mathfrak s_0 + \s_2} \leq 
\e C N_n^{\tau} \Vert u_n - u_{n - 1} \Vert_{\mathfrak s_0 + \mu} \stackrel{\eqref{hn}}{\leq} 
\e^2 \g^{-1} 
C C_* N_{n}^{\tau - \sigma_1} \leq \gamma_{n-1} - \gamma_{n} =: \rho = \g 2^{-n } 
$$ 
for $\e \gamma^{-1}$ small enough, because $ \s_1  > \t   $ (see \eqref{hn}, \eqref{perdite stima tame linearizzato}). 
Then, by the definitions \eqref{G-n+1} and \eqref{Omegainfty}, we have 
$$
{\cal G}_{n} := {\cal G}_{n-1} \cap \L_{\infty}^{2 \g_{n-1}} (u_{n-1}) 
\stackrel{\eqref{cantorinclu}} \subseteq \bigcap_{\nu \geq 0} \Lambda_{\nu}^{\gamma_{n - 1}}(u_{n - 1}) 
\subset \Lambda_{n+1}^{\g_{n-1}}(u_{n-1})
\stackrel{\eqref{primoste}} \subseteq \Lambda_{n+1}^{\g_n}(u_n).
$$
Next, for all $ \l \in {\cal G}_n  \subset \Lambda_{n+1}^{\g_{n-1}}(u_{n-1}) \cap 
 \Lambda_{n+1}^{\g_n}(u_n)  $ both $ r_j^{n+1} (u_{n-1}) $ and $ r_j^{n+1} (u_{n}) $ are well defined, and
 we deduce by Theorem \ref{thm:abstract linear reducibility}-${\bf (S3)}_\nu $
 with $ \nu = n+1 $,  that 
\be\label{vicin+1}
| r^{n+1}_j (u_n) - r_j^{n+1} (u_{n-1})| \stackrel{\eqref{Delta12 rj}} 
\lessdot \e  \| u_{n-1} - u_n \|_{\mathfrak s_0 + \s_2} \, .
\ee
Moreover \eqref{autovcon}  (with $ \nu = n+1 $) and \eqref{stima R 1} imply that 
\begin{eqnarray}\label{diffrkn}
 | r_j^{\infty}(u_{n -1}) - r_j^{n + 1}(u_{n - 1})| +
|r_j^{\infty}(u_{n}) - r_j^{n + 1}(u_{n})| & \lessdot &
\e (1 + \| u_{n-1} \|_{\mathfrak s_0 +  \b + \s}+  \| u_n \|_{\mathfrak s_0 +  \b + \s}) N_{n}^{-\a} \nonumber \\ 
& \lessdot &  \e 
N_{n}^{-\a} 
\end{eqnarray}
because $ \s + \b < \mu $ and $ \| u_{n-1} \|_{\mathfrak s_0 +  \mu} + $ $  \| u_n \|_{\mathfrak s_0 +  \mu} \leq 2 $ by 
$ {\bf (S1)}_{n-1} $ and $ {\bf (S1)}_n $.
Therefore, for all $\lambda \in {\cal G}_{n}$,  $ \forall j \in \Z $,   
\begin{align}
\big| r_j^{\infty}(u_{n}) - r_j^{\infty}(u_{n - 1}) \big|
& \leq 
\big| r_j^{n + 1}(u_{n}) - r_j^{n + 1}(u_{n-1}) \big| 
+ | r_j^{\infty}(u_{n}) - r_j^{n + 1}(u_{n})| + | r_j^{\infty}(u_{n -1}) - r_j^{n + 1}(u_{n - 1})| 
\nonumber\\
& \stackrel{\eqref{vicin+1}, \eqref{diffrkn}} \lessdot
\e \| u_n - u_{n - 1} \|_{\mathfrak s_0+ \s_2}  + \e  N_{n}^{-\a}
\stackrel{\eqref{hn}} \lessdot  \e N_{n}^{-\a}  \label{variazione autovalori finali in u}
\end{align}
because  
$ \s_1 > \a $ (see \eqref{alpha-beta}, \eqref{hn}).  
Finally  \eqref{vici}, \eqref{variazione autovalori finali in u}, \eqref{coefficienti costanti 2}, 
$\| u_n \|_{\mathfrak s_0 + \mu}\leq 1$,   imply \eqref{marco}. 
\end{pf}

By definition, $ R_{ljk}(u_n) \subset {\cal G}_n $ (see  \eqref{R-ljk(u-n)}) and, by \eqref{inclusione-1},
for all $ |l| \leq N_n $, we have 
$ R_{ljk} (u_n) \subseteq R_{ljk} (u_{n-1}) $.
On the other hand $ R_{ljk}(u_{n-1}) \cap {\cal G}_n = \emptyset $, see \eqref{G-n+1}. 
As a consequence, $ \forall |l| \leq N_n $, 
$ R_{ljk} (u_n) = \emptyset $, and 
\be\label{parametri cattivi}
{\cal G}_{n} \setminus{\cal G}_{n+1}  \stackrel{\eqref{natale}}  
\subseteq \bigcup_{|l|> N_{n}, j,k \in \Z} R_{ljk}(u_{n}) \, , 
\quad \forall n \geq 1.   
\ee

\begin{lemma}\label{risonanti-2}
Let $n \geq 0$. 
If $R_{ljk}(u_{n}) \neq \emptyset$, then $|j^{3}-k^{3}| \leq 8 |\bar\o \cdot l|$.  
\end{lemma}

\begin{pf}
If $R_{ljk}(u_{n})\,\neq\,\emptyset$ then there exists $\lambda\in \Lambda$ such that
$  |\ii \lambda\bar\o\cdot l + \mu_{j}^{\infty}(\lambda,u_{n}(\lambda))- 
\mu_{k}^{\infty}(\lambda,u_{n}(\lambda)) | < $ $ 2 \gamma_{n} |j^{3}-k^{3} | \langle l \rangle^{-\tau} $ and, 
therefore, 
\begin{equation}\label{lambda-j-lambda-k}
|\mu_{j}^{\infty}(\lambda,u_{n}(\lambda)) - \mu_{k}^{\infty}(\lambda,u_{n}(\lambda)) | 
< 2\gamma_{n} |j^{3}-k^{3}| \langle l \rangle^{-\tau}\, + 2 |\bar\o\cdot l|.
\end{equation}
Moreover, by \eqref{espressione autovalori}, \eqref{coefficienti costanti 1}, \eqref{autofinali}, for $\e$ small enough, 
\be\label{lower-bound-jk}
|\mu_{j}^{\infty} - \mu_{k}^{\infty}| 
 \geq |m_3| |j^{3}-k^{3}| - |m_1| |j-k| - |r_j^{\infty}| - |r_k^{\infty}|   \geq 
\frac{1}{2} |j^{3}-k^{3}| - C \e |j - k| - C \e 
\geq \frac{1}{3} |j^{3}-k^{3}| 
\ee
if  $j \neq k$. 
Since $\g_n \leq 2\g$ for all $n \geq 0$, $\g \leq 1/ 48$, by \eqref{lambda-j-lambda-k} and \eqref{lower-bound-jk} we get 
$$
2 |\bar\o \cdot l| 
\geq \Big(\frac13 -\frac{4\gamma}{\langle l \rangle^{\tau}} \Big) |j^{3}-k^{3}|\ 
\geq \frac14|j^{3}- k^3| 
$$
proving the Lemma. 
\end{pf}

\begin{lemma}\label{risonanti-3}
For all $n \geq 0$, 
\begin{equation}\label{stima-risonanti}
 |R_{ljk}(u_{n})| \leq C \gamma \left\langle l\right\rangle^{-\tau}.
\end{equation}
\end{lemma}

\begin{pf}
Consider the function $ \phi : \L \to \C$ defined by
\begin{eqnarray} 
\phi(\l) & := &  \ii \l \bar\o\cdot l +\mu_{j}^{\infty}(\l)- \mu_{k}^{\infty}(\l)
\nonumber
\\
& \stackrel{\eqref{espressione autovalori}} = &\ii \l \bar\o\cdot l - \ii {\tilde m}_3(\lambda)(j^{3}-k^{3}) + 
 \ii {\tilde m}_1(\lambda)(j-k) + r_{j}^{\infty}(\lambda)- r_{k}^{\infty}(\lambda)   \nonumber 
\end{eqnarray}
where  $ {\tilde m}_3 (\l) $,   
$ {\tilde m}_1 (\l) $, $ r^{\infty}_j (\l) $, $ \mu_{j}^{\infty}(\l)$, 
are defined for all $ \l \in \L$ and satisfy \eqref{autofinali}
by  $ \| u_n \|^{\Lipg}_{\mathfrak s_0 + \mu, \mG_n} \leq 1$ (see $({\cal P}1)_n$).   
Recalling $ | \cdot |^\lip \leq \g\inv | \cdot |^\Lipg $    
and using \eqref{autofinali}  
\be\label{31}
| \mu_{j}^{\infty} - \mu_{k}^{\infty} |^{\rm lip} \leq |{\tilde m}_3|^{\rm lip} |j^{3} - k^{3}| 
+ | {\tilde m}_1|^{\rm lip} |j - k| + |r_{j}^{\infty}|^{\rm lip} + |r_{k}^{\infty}|^{\rm lip} \leq C \e \g^{-1} |j^{3} - k^{3}|\,.
\ee
Moreover Lemma \ref{risonanti-2} implies that, $\forall \l_1, \l_2 \in \L $, 
$$
|\phi(\l_1) - \phi(\l_2)|  \geq  \big( |\bar{\o} \cdot l| - |\mu_j^{\infty} - \mu_k^{\infty}|^{\rm lip} \big) |\l_1 - \l_2| 
\stackrel{\eqref{31}} \geq \big(\frac18 - C \e\g^{-1} \big)|j^3 - k^3| |\l_1 - \l_2| \geq \frac{|j^3 - k^3|}{9} |\l_1 - \l_2| 
$$
for $\e\gamma^{-1}$ small enough.
Hence 
$$
|R_{ljk}(u_n)| \leq \frac{4 \g_n|j^{3}-k^3|}{\langle l  \rangle^{\t}} \frac{9}{|j^3 - k^3|} 
\leq \frac{72 \g}{\la l \ra^{\t}} \,, 
$$
which is \eqref{stima-risonanti}.
\end{pf}

Now we prove $(\mP4)_0$. 
We observe that, for each fixed $l$, all the indices $j,k$ such that $R_{ljk}(0) \neq \emptyset$ are confined in the ball $j^2 + k^2 \leq 16 |\bar\om| |l|$, because 
$$
|j^3 - k^3| = |j-k| |j^2+jk+k^2| 
\geq j^2 + k^2 - |jk| 
\geq \frac12\, (j^2 + k^2) \, , 
\quad \forall j,k \in \Z, \ j \neq k,
$$
and $|j^{3}-k^{3}| \leq 8 |\bar\o| |l|$ by Lemma \ref{risonanti-2}. 
As a consequence 
$$
| {\cal G}_0 \setminus {\cal G}_1 | \stackrel{\eqref{natale}} = 
\Big|  \bigcup_{l,j,k} R_{ljk}(0) \Big| 
\leq \sum_{l \in \Z^\nu} \sum_{j^2 + k^2 \leq 16 |\bar\om| |l|}  |R_{ljk}(0)| 
\stackrel{\eqref{stima-risonanti}} 
\lessdot \sum_{l \in \Z^\nu}  \g \la l \ra^{-\tau+1}  
= C \g
$$
if $\t > \nu + 1$. 
Thus the first estimate in \eqref{Gmeasure} is proved, taking a larger $C_*$ if necessary.

Finally, $(\mP4)_n$ for $n \geq 1$, follows by  
\begin{eqnarray*}
|{\cal G}_{n } \setminus {\cal G}_{n+1}| 
& \stackrel{\eqref{parametri cattivi}} \leq 
& \sum_{|l|> N_{n} |j|,|k|\leq C |l|^{1/2}} |R_{ljk}(u_{n})| 
\stackrel{\eqref{stima-risonanti}} \lessdot 
 \sum_{|l|> N_{n} |j|,|k|\leq C|l|^{1/2}} \g \langle l \rangle^{-\t} 
\\
& \lessdot
& \sum_{|l| > N_n} \g \langle l \rangle^{-\t + 1} 
\lessdot \g N_{n}^{-\t + \nu}
\leq C \g N_{n}^{-1}
\end{eqnarray*}
and \eqref{Gmeasure} is proved. 
The proof of Theorem \ref{iterazione-non-lineare} is complete.
\end{pf}

\subsection{Proof of Theorems  \ref{thm:main}, 
\ref{thm:mainH}, \ref{thm:mainrev}, \ref{thm:reducibility}
and \ref{cor:stab}}\label{sec:proof}

\textsc{Proof of Theorems \ref{thm:main}, \ref{thm:mainH}, \ref{thm:mainrev}.} 
Assume that $ f \in C^q $ satisfies the assumptions in Theorem \ref{thm:main}
or in Theorem \ref{thm:mainrev} with a smoothness exponent $ q  := q(\nu)  \geq \mathfrak s_0 + \mu + \b_1 $ which depends only on  $ \nu $ once we have fixed $ \t := \nu + 2 $ 
(recall that $ \mathfrak s_0 := (\nu + 2 ) \slash 2 $, 
$ \b_1 $  is defined in \eqref{sigma-beta-kappa} and $ \mu $ in  \eqref{perdite stima tame linearizzato}). 

For $ \g = \e^a $, $ a \in (0,1) $ the smallness condition 
$ \e \g^{- 1} = \e^{1- a}  < \d  $ of Theorem  \ref{iterazione-non-lineare} is satisfied.  
Hence on the Cantor set ${\cal G}_{\infty} := \cap_{n \geq 0} {\cal G}_{n} $, 
the sequence $ u_{n}(\l) $ is well defined and 
converges in norm $ \|\cdot \|_{\mathfrak s_{0}+\mu, \mG_\infty}^{\Lipg}$ 
(see  \eqref{hn}) to a solution $u_{\infty}(\l)$ of 
$$
F(\lm, u_\infty(\lm)) = 0 \quad {\rm with}  \quad \sup_{\l \in {\cal G}_\infty} \| u_\infty (\l) \|_{\mathfrak s_0 + \mu} \leq C \e \g^{-1} 
= C \e^{1-a} \, ,
$$
namely $u_{\infty}(\l)$ is a solution of the 
perturbed KdV equation \eqref{eq:invqp} with $ \om = \l \bar \om $.    
Moreover, 
by \eqref{Gmeasure},  the measure of the complementary set satisfies 
\[
|\Lambda \setminus {\cal G}_{\infty}| 
\leq \sum_{n \geq 0}|{\cal G}_{n } \setminus {\cal G}_{n+1} | 
\leq C \g + \sum_{n \geq 1} \g C N_{n}^{-1} \leq C \g = C \e^{a} \, ,
\qedhere 
\]
proving \eqref{Cmeas}. 
The proof of Theorem \ref{thm:main} is complete. 
In order to finish the proof of Theorems \ref{thm:mainH} or 
\ref{thm:mainrev}, it remains to prove the linear stability of the solution, namely 
Theorem 
\ref{cor:stab}. 

\smallskip

\noindent
\textsc{Proof of Theorem \ref{thm:reducibility}.} 
Part $(i) $  
follows by \eqref{L-coniugato}, Lemma \ref{stime-tame-coniugio}, 
Theorem \ref{teoremadiriducibilita} (applied to the solution $ u_\infty (\l) $)
with the exponents $ \bar \s := \s + \b + 3 $, $ \L_\infty (u) := \L_\infty^{2\g} (u) $, see \eqref{Omegainfty}.
Part ($ii$) follows by the dynamical interpretation of the conjugation procedure, as explained in 
section \ref{sec: dyn redu}. 
Explicitely, in sections \ref{sec:regu} and \ref{sec:redu},  we 
have  proved that
\[
\mL = {\cal A} B \rho W \mL_\infty W^{-1} B^{-1} {\cal A}^{-1}, 
\quad W := {\cal M}  {\cal T} \mS \Phi_\infty \, .
\]
By the arguments in Section \ref{sec: dyn redu} we deduce that  
a curve $h(t)$ in the phase space $H^s_x$ is a solution of the dynamical system \eqref{KdV:lin}
if and only if the transformed curve 
\be\label{vh}
v(t) := W^{-1}(\om t) B^{-1} {\cal A}^{-1}(\om t) h(t) 
\ee
(see notation \eqref{notationA}, 
Lemma \ref{lemma:stime stabilita Phi 12}, \eqref{Phi infty (ph) - I}) 
is a solution of the constant coefficients dynamical system \eqref{Lin: Red}.

\smallskip

\noindent
\textsc{Proof of Theorem \ref{cor:stab}.}
If all $ \mu_j $ are purely imaginary, the Sobolev norm of the solution $ v(t) $ of \eqref{Lin: Red} 
is constant in time, see
\eqref{constant v}. 
We now show that also the Sobolev norm of the solution $ h(t) $ in \eqref{vh} does not grow in time. 
For each $t \in \R$, $ {\cal A}(\om t) $ and  
$W(\om t)$ are transformations of the phase space $H^s_x $ that depend quasi-periodically on time, and
satisfy,  by \eqref{A(ph)}, \eqref{Phi mM mS (ph)}, \eqref{Phi infty (ph) - I},
\be \label{nuova carla}
\| {\cal A}^{\pm 1}(\om t) g \|_{H^s_x} 
+ \| W^{\pm 1}(\om t) g \|_{H^s_x} 
\leq C(s) \| g \|_{H^s_x} \, , \quad \forall t \in \R, \ \forall g = g(x) \in H^s_x,
\ee
where the constant $C(s)$ depends on $\| u \|_{s + \s + \b + \mathfrak s_0} < + \infty $. 
Moreover, the transformation $B$ is a quasi-periodic reparametrization of the time variable 
(see \eqref{time repar}), namely
\be \label{B t tau}
Bf(t) = f(\psi(t)) = f(\t), \quad 
B^{-1}f(\t) = f(\psi^{-1}(\t)) = f(t) \quad \forall f : \R \to H^s_x,
\ee
where $\t = \psi(t) := t + \a(\om t)$, $t = \psi^{-1}(\t) = \t + \tilde \a(\om \t)$ and $\a$, $\tilde\a$ are defined in Section \ref{step-2}.
Thus

\begin{align*}
\| h(t) \|_{H^s_x} 
& \stackrel{\eqref{vh}} =  
\| {\cal A}(\om t) B W(\om t) v (t) \|_{H^s_x} 
\stackrel{\eqref{nuova carla}}  \leq 
C(s) \| B W(\om t) v (t) \|_{H^s_x} 
\stackrel{\eqref{B t tau}}  = C(s) \| W(\om \t) v (\t) \|_{H^s_x} 
\\ & 
\stackrel{\eqref{nuova carla}}   \leq C(s) \| v (\t) \|_{H^s_x} 
\stackrel{\eqref{constant v}}   = C(s) \| v (\t_0) \|_{H^s_x} 
\stackrel{\eqref{vh}} = C(s) \| W^{-1}(\om \t_0) B^{-1} {\cal A}^{-1}(\om \t_0) h(\t_0) \|_{H^s_x} 
\\ & 
\stackrel{\eqref{nuova carla}} \leq  C(s) \| B^{-1} {\cal A}^{-1}(\om \t_0) h(\t_0) \|_{H^s_x} 
\stackrel{\eqref{B t tau}}  = C(s) \| {\cal A}^{-1}(0) h(0) \|_{H^s_x} 
\stackrel{\eqref{nuova carla}}   \leq C(s) \| h(0) \|_{H^s_x}
\end{align*} 
having chosen  $\t_0 := \psi(0) = \a(0)$ 
(in the reversible case, $\a$ is an odd function, and so $\a(0) = 0$). Hence \eqref{stability s} is proved. 
To prove \eqref{stability epsilon}, we collect the estimates 
\eqref{A(ph)-I}, \eqref{Phi mM mS (ph) - I}, \eqref{Phi infty (ph) - I} into 
\be \label{nuova rosa}
\| ({\cal A}^{\pm 1}(\om t) - I) g \|_{H^s_x} 
+ \| (W^{\pm 1}(\om t) - I) g \|_{H^s_x} 
\leq \e \g^{-1} C(s) \| g \|_{H^{s+1}_x} \,, \quad \forall t \in \R, \ \forall g \in H^s_x,
\ee
where the constant $C(s)$ depends on $\| u \|_{s + \s + \b + \mathfrak s_0}$. 
Thus  
\begin{eqnarray*}
\| h(t) \|_{H^s_x} 
& \!\!\!\!\!\!  \stackrel{\eqref{vh}}  = \!\!\!\!\!\!  & \| {\cal A}(\om t) B W(\om t) v (t) \|_{H^s_x}  \leq 
\| B W(\om t) v (t) \|_{H^s_x}  + \| ({\cal A}(\om t)-I) B W(\om t) v (t) \|_{H^s_x} \\ 
& \!\!\!\!\!\!  \stackrel{\eqref{B t tau}\eqref{nuova rosa}}  \leq \!\!\!  \!\!\! & \| W(\om \t) v(\t) \|_{H^s_x}
+ \e \g^{-1} C(s) \| B W(\om t) v (t) \|_{H^{s+1}_x}  \\ 
& \!\!\! \!\!\!  \stackrel{\eqref{B t tau}}  = \!\!\!\!\!\!  & \| W(\om \t) v(\t) \|_{H^s_x} 
+ \e \g^{-1} C(s) \| W(\om \t) v (\t) \|_{H^{s+1}_x}  \\ 
& \!\!\!\!\!\!  \stackrel{\eqref{nuova carla}} \leq \!\!\! \!\!\! & \| v(\t) \|_{H^s_x} + \| (W(\om \t) - I) v(\t) \|_{H^s_x} 
+ \e \g^{-1} C(s) \| v (\t) \|_{H^{s+1}_x}  \\ 
& \!\!\! \!\!\! \stackrel{\eqref{nuova rosa}}  \leq \!\!\!\!\!\!  &   \| v(\t) \|_{H^s_x} + \e \g^{-1} C(s) \| v(\t) \|_{H^{s+1}_x}  
\stackrel{\eqref{constant v}}  = \| v(\t_0) \|_{H^s_x} + \e \g^{-1} C(s) \| v(\t_0) \|_{H^{s+1}_x}  \\ 
& \!\!\!\!\!\!  \stackrel{\eqref{vh}}  = \!\!\!\!\!\!  & \| W^{-1}(\om \t_0) B^{-1} {\cal A}^{-1}(\om \t_0) h(\t_0) \|_{H^s_x} 
+ \e \g^{-1} C(s) \| W^{-1}(\om \t_0) B^{-1} {\cal A}^{-1}(\om \t_0) h(\t_0)  \|_{H^{s+1}_x} \, .
\end{eqnarray*}
Applying the same chain of inequalities at $ \t = \t_0 $, $ t = 0 $, we get that the last term is
$$  
\leq \| h(0) \|_{H^s_x} + \e \g^{-1} C(s) \| h(0) \|_{H^{s+1}_x} \, ,
$$
proving the second inequality in \eqref{stability epsilon} with $ \mathtt a := 1 - a $. 
The first one follows similarly.

\section{Appendix A. General tame and Lipschitz estimates}

In this Appendix we present standard tame and Lipschitz estimates 
for composition of functions and changes of variables
which are used in the paper. Similar material is contained in 
 \cite{Hormander-geodesy}, \cite{Ioo-Plo-Tol}, \cite{Berti-Bolle-Ck-Nodea}, 
 \cite{Baldi-Benj-Ono}.  

We  first remind classical embedding, algebra, interpolation and tame 
estimates in the Sobolev spaces  $ H^s := H^{s}(\T^d,\C) $ and $ W^{s, \infty} := W^{s, \infty}(\T^d,\C) $, $ d \geq 1$ .

\begin{lemma} \label{lemma:standard Sobolev norms properties}
 Let 
 $ s_0 > d/2$. Then 
\\
$(i)$ {\bf Embedding.} $\| u \|_{L^\infty} \leq C(s_0) \| u \|_{s_0}$ for all $u \in H^{s_0} $.
\\
$(ii)$ {\bf Algebra.} $\| uv \|_{s_0} \leq C(s_0) \| u \|_{s_0} \| v \|_{s_0}$ for all 
$u, v \in H^{s_0}$.
\\
$(iii)$ {\bf Interpolation.} For $0 \leq s_1 \leq s \leq s_2$, $s = \lm s_1 + (1-\lm) s_2$, 
\begin{equation} \label{interpolation GN}
\| u \|_{s} \leq  \| u \|_{s_1}^\lm \| u \|_{s_2}^{1-\lm} \, , 
\quad \forall u \in H^{s_2} \, .
\end{equation}
Let $ a_0, b_0 \geq 0$ and $ p,q >  0 $. For all $ u \in H^{a_0 + p + q} $, $  v \in H^{b_0 + p + q} $, 
\be \label{interpolation estremi fine}
\| u \|_{a_0 + p} \| v \|_{b_0 + q}
\leq 
\| u \|_{a_0 + p + q} \| v \|_{b_0}
+ \| u \|_{a_0} \| v \|_{b_0 + p + q} \, .
\ee
Similarly, for the $|u|_{s, \infty} := \sum_{|\b| \leq s} | D^\b u |_{L^\infty} $ norm,  
\begin{equation} \label{interpolation GN-s}
| u |_{s, \infty} \leq C(s_1, s_2) | u |_{s_1, \infty}^\lm | u |_{s_2, \infty}^{1-\lm} \, , 
\quad \forall u \in W^{s_2, \infty} \, .
\end{equation}
and $ \forall u \in W^{a_0 + p + q, \infty} $, $  v \in W^{b_0 + p + q, \infty} $,
\be \label{interpolation estremi fine-s}
| u |_{a_0 + p, \infty} | v |_{b_0 + q, \infty}
\leq C(a_0, b_0, p,q)\big(  
| u |_{a_0 + p + q, \infty} | v |_{b_0, \infty}
+ | u |_{a_0, \infty} | v |_{b_0 + p + q, \infty}\big)  \, .
\ee
\noindent
$(iv)$ {\bf Asymmetric tame product.} For $s \geq s_0$, 
\begin{equation} \label{asymmetric tame product}
\| uv \|_s \leq C(s_0) \|u\|_s \|v\|_{s_0} + C(s) \|u\|_{s_0} \| v \|_s \, , 
\quad \forall u,v \in H^s \, .
\end{equation}

\noindent
$(v)$ {\bf Asymmetric tame product in $W^{s,\infty}$.} For $s \geq 0$, $s \in \N$,
\begin{equation} \label{tame product infty}
| uv |_ {s, \infty}  \leq \tfrac32 \, | u |_ {L^\infty} |v|_ {s, \infty}    + C(s) |u|_ {s, \infty}  |v|_ {L^\infty}  \, , 
\quad \forall u,v \in W^{s,\infty} \, .
\end{equation}

\noindent
$(vi)$ {\bf Mixed norms asymmetric tame product.} For $s \geq 0$, $s \in \N$,
\begin{equation} \label{mixed norms tame product}
\| uv \|_s \leq \tfrac32 \,  |u|_ {L^\infty}  \| v \|_s + C(s)| u |_ {s, \infty} \| v \|_0  \, , 
\quad \forall u \in W^{s,\infty} \, , \ v \in H^s   \, .
\end{equation}
If $u := u(\l)$ and $v := v(\l)$ depend in a lipschitz way on $\l \in \L \subset \R $, 
all the previous statements hold if we replace the norms 
$\Vert \cdot \Vert_s$, $\vert \cdot \vert_ {s, \infty}  $ with the norms $\Vert \cdot \Vert_s^{\Lipg}$, $\vert \cdot \vert_ {s, \infty}^{\Lipg}$.

\end{lemma} 

\begin{pf} 
The interpolation estimate \eqref{interpolation GN} 
for the Sobolev norm \eqref{Hs1} follows by H\"older inequality, see 
also \cite{Moser-Pisa-66}, page 269. 
Let us prove \eqref{interpolation estremi fine}. Let
$ a = a_0 \lm + a_1 (1-\lm) $,  $ b = b_0 (1-\lm) + b_1 \lm $,  $ \lm \in [0,1] $. 
Then  \eqref{interpolation GN} implies
\be \label{ulti2}
\| u \|_a \| v \|_b
\leq 
 \big( \| u \|_{a_0} \| v \|_{b_1} \big)^{\lm}
\big( \| u \|_{a_1} \| v \|_{b_0} \big)^{1-\lm} \leq \l \| u \|_{a_0} \| v \|_{b_1} + (1- \l) \| u \|_{a_1} \| v \|_{b_0}
\ee
by Young inequality. 
Applying \eqref{ulti2} with
$ a = a_0 + p $,  $ b = b_0 + q $,  $ a_1 = a_0 + p + q $,  $ b_1 = b_0 + p + q $, 
then $ \lm = q \slash (p+q) $ and  we get \eqref{interpolation estremi fine}. 
Also the interpolation estimates
 \eqref{interpolation GN-s} are classical (see e.g. \cite{Hormander-geodesy}, \cite{Berti-Bolle-Procesi-AIHP-2010})
 and \eqref{interpolation GN-s} implies  \eqref{interpolation estremi fine-s} as above.

$(iv)$: see the Appendix of \cite{Berti-Bolle-Procesi-AIHP-2010}. 
$(v)$: we write, in the standard multi-index notation, 
\begin{equation}  \label{derivate pure}
D^\a(uv) = \sum_{\b+\g = \a} C_{\b, \g} (D^\b u) (D^\g v) = 
u D^\a v  + \sum_{\b+\g = \a, \b \neq 0} C_{\b, \g} (D^\b u) (D^\g v)  \, .
\end{equation}
Using $ |(D^\b u)(D^\g v)|_ {L^\infty}  \leq |D^\b u|_ {L^\infty}  |D^\g v|_ {L^\infty} \leq |u|_{|\b|, \infty} |v|_{|\g|, \infty} $, 
and the interpolation inequality \eqref{interpolation GN-s} for every $ \b \neq 0$   
with $\lm := |\b| / |\a| \in (0,1]$ (where $ |\a| \leq s $),  we get, for any $K > 0$, 
\begin{align}
C_{\b, \g} |D^\b u|_ {L^\infty}  |D^\g v|_ {L^\infty}  
& 
\leq C_{\b, \g}  C(s) 
\big( |v|_{L^\infty}  |u|_ {s, \infty}  \big)^{\lm} \big( |v|_ {s, \infty}  |u|_ {L^\infty}  \big)^{1 - \lm} \nonumber \\
& = 
\frac{C(s)}{K} \big[ (K C_{\b, \g})^{\frac{1}{\lm}} |v|_ {L^\infty}  |u|_ {s, \infty}  \big]^{\lm} 
\big( |v|_ {s, \infty}  |u|_ {L^\infty}  \big)^{1 - \lm} \nonumber 
\\ 
& \leq \frac{C(s)}{K} \, 
\big\{ (K C_{\b, \g})^{\frac{|\a|}{|\b|}} |v|_ {L^\infty}  |u|_ {s, \infty}  \, + \, |v|_ {s, \infty}  |u|_ {L^\infty}  \big\}. \label{riga3}
\end{align}
Then \eqref{tame product infty} follows by \eqref{derivate pure}, \eqref{riga3} taking $ K := K(s) $ large enough.  
$(vi)$: same proof as $(v)$, using the elementary inequality 
$\|(D^\b u)(D^\g v) \|_0 \leq |D^\b u|_{L^\infty}  \| D^\g v \|_0 $.
\end{pf}

We now recall classical tame estimates for composition of functions, see \cite{Moser-Pisa-66}, 
section 2, pages 272--275, and \cite{Rabinowitz-tesi-1967}-I, Lemma 7 in the Appendix, pages 202--203.

A function $ f : \T^d \times B_1 \to \C $, where $B_1 := \{ y \in \R^m : |y| < 1\} $, induces 
the composition operator 
\be\label{comp}
\tilde f(u)(x) := f(x,u(x),Du(x),\ldots,D^p u(x)) 
\ee
where $D^k u(x)$ denotes the partial derivatives $\pa_x^\a u(x)$ of order $|\a|=k$ 
(the number $ m $ of $ y $-variables depends on $p, d $).

\begin{lemma} {\bf (Composition of functions)}
\label{lemma:composition of functions, Moser}
Assume $ f \in C^r (\T^d \times B_1)$. Then

$(i)$ 
For all $ u \in H^{r+p} $ such that  $ |u |_{p, \infty}  < 1 $, 
the composition operator \eqref{comp} is well defined and 
\[
\| \tilde f(u) \|_r
\leq C \| f \|_{C^r} (\|u\|_{r+p} + 1) 
\]
where the constant $C $ depends on $ r,d,p $.  If $ f \in C^{r+2} $,  
then, for all $ |u|_{p, \infty} $, $  | h |_{p, \infty}  < 1 / 2 $,  
\begin{align*}
\big\| \tilde f(u+h) - \tilde f (u) \big\|_r 
& \leq C \| f \|_{C^{r+1}} \, ( \| h \|_{r+p} + | h |_{p,\infty} \| u \|_{r+p}) \, ,  
\\
\big\| \tilde f(u+h) - \tilde f (u) - \tilde f'(u) [h] \big\|_r 
& \leq C \| f \|_{C^{r+2}} \, | h |_{p,\infty}  ( \| h \|_{r+p} + | h |_{p,\infty} \| u \|_{r+p}) \, . 
\end{align*}
$(ii)$ The previous statement also holds replacing 
$\| \ \|_r$  
with the norms $| \ |_{r, \infty} $.  
\end{lemma}

\begin{lemma} {\bf (Lipschitz estimate on parameters)}  \label{lemma:Lip generale}
Let $d \in \N$, $d/2 < s_0 \leq s$, $p \geq 0$, $\g>0$.
Let $ F $ be a $ C^1 $-map satisfying the tame estimates: $ \forall \| u \|_{s_0+p} \leq 1 $, $ h \in H^{s+p} $, 
\begin{align} 
\label{aux tame 1}
\| F(u) \|_s 
& \leq C(s) (1 + \| u \|_{s+p}) \, ,  \\
\label{aux tame 2}
\| \pa_u F(u)[h] \|_s 
& \leq C(s) (\| h \|_{s+p} + \| u \|_{s+p} \| h \|_{s_0 + p} ) \, .  
\end{align}
For $\Lambda \subset \R $, 
let $ u(\lm) $ be a Lipschitz family of functions with 
$ \| u \|_{s_0 +p}^{\Lipg} \leq 1  $ (see \eqref{def norma Lipg}).
Then 
\[
\| F(u) \|_s^\Lipg \leq C(s) \big( 1 + \| u \|_{s+p}^\Lipg \big).
\]
The same statement also holds when all the norms $\| \ \|_s $ are replaced by  $| \ |_{s, \infty} $. 
\end{lemma}

\begin{pf}
By \eqref{aux tame 1} we get $ \sup_\lm \| F(u(\lm)) \|_s 
\leq 
C(s) ( 1 + \| u \|_{s+p}^{\Lipg}) $. Then, 
denoting $ u_1 := u(\lm_1)$ and $h := u(\lm_2) - u(\lm_1)$, we have 
\begin{align*}
\| F(u_2) - F(u_1) \|_s 
& 
\leq \int_0^1 \| \pa_u F(u_1 + t (u_2-u_1))[h] \, \|_s \, dt 
\\ & \stackrel{\eqref{aux tame 2}} \leq_s \| h \|_{s+p} 
+ \| h \|_{s_0 + p} \int_0^1 \big( (1-t) \| u(\lm_1) \|_{s+p} + t \| u(\lm_2) \|_{s+p} \big) \, dt
\end{align*}
whence 
\begin{align*}
\g \, \sup_{\begin{subarray}{c} \lm_1, \lm_2 \in \Lambda \\ \lm_1 \neq \lm_2 \end{subarray}} 
\frac{\| F(u(\lm_1)) - F(u(\lm_2)) \|_{s}}{|\lm_1 - \lm_2|} \, 
& 
\leq_s \| u \|_{s+p}^\Lipg 
+  \| u \|_{s_0 + p}^\Lipg \sup_{\lm_1, \lm_2} \big( \, \| u(\lm_1) \|_{s+p} + \, \| u(\lm_2) \|_{s+p} \big)
\\ 
& \leq_s \| u \|_{s+p}^\Lipg
+ \| u \|_{s_0 + p}^\Lipg \| u \|_{s+p}^\Lipg
\leq C(s) \| u \|_{s+p}^\Lipg \, , 
\end{align*}
because $ \| u \|_{s_0 +p}^{\Lipg} \leq 1  $, and  the lemma follows.
\end{pf}

\smallskip

The next lemma is also classical, see for example \cite{Hormander-geodesy}, Appendix, and \cite{Ioo-Plo-Tol}, Appendix G. 
The present version is proved in \cite{Baldi-Benj-Ono}, adapting Lemma 2.3.6 on page 149 of \cite{Hamilton}, except for the part on the Lipschitz dependence on a parameter, which is proved here below.

\begin{lemma} {\bf (Change of variable)}  \label{lemma:utile} 
Let $p:\R^d \to \R^d$ be a $2\p$-periodic function in $W^{s,\infty}$, 
$ s \geq 1$, with 
$ |p|_{1, \infty} \leq 1/2 $. Let $f(x) = x + p(x)$. Then:

$(i)$ $f$ is invertible, its inverse is $f\inv(y) = g(y) = y + q(y)$ where $q$ is $ 2 \pi $-periodic, $q \in W^{s,\infty}(\T^d,\R^d)$, and 
 $|q|_{s, \infty}  \leq C |p|_{s, \infty} $. More precisely,
\begin{equation} \label{stime-q}
| q |_{L^\infty} = | p |_{L^\infty}, \quad   
| Dq |_{L^\infty} \leq 2 | Dp |_{L^\infty} , \quad
| Dq |_{s-1, \infty} \leq C | Dp |_{s-1, \infty}.
\end{equation}
where the constant $C$ depends on $d, s$.

Moreover, suppose that $p = p_\lm$ depends in a Lipschitz way by a parameter $\lambda \in \L \subset \R $, 
and suppose, as above, that 
$|D_x p_\lm|_ {L^\infty}  \leq 1/2$ for all $\lm$. 
Then $q = q_\lm$ is also Lipschitz in $\lambda$, and
\begin{equation}\label{stime-lipschitz-q}
|q|_{s, \infty}^{\Lipg}
\leq  C \Big( |p|_ {s, \infty} ^{\Lipg} + \big\{ \sup_{\lm \in \L} |p_\lm|_{s+1, \infty} \big\} \, |p|_{L^\infty}^{\Lipg} \Big)
 \leq  C  |p|_{s+1, \infty}^{\Lipg}, 
\end{equation}
The constant $C$ depends on $d, s $ (and is independent on $\g$).

$(ii)$ If $u \in H^s (\T^d,\C)$, then $u\circ f(x) = u(x+p(x))$ is also in $H^s $, and, with the same $C$ as in $(i)$, 
\begin{align}
\| u \circ f \|_s 
& \leq  C (\|u\|_s + |Dp|_{s-1, \infty} \|u\|_1),
\label{tame-cambio-di-variabile}
\\
\| u \circ f - u \|_s 
& \leq C \big( | p |_{L^\infty} \| u \|_{s + 1}  + |p|_{s, \infty} \| u \|_{2} \big) ,
\label{cambio di variabile meno identita} \\
\| u \circ f \|_{s}^{\Lipg}
& \leq C  \, 
\big( \| u \|_{s+1}^{\Lipg} + |p|_{s, \infty}^{\Lipg} \| u \|_2^{\Lipg} \big). 
\label{tame-lipschitz-cambio-di-variabile}
\end{align}
\eqref{tame-cambio-di-variabile}, \eqref{cambio di variabile meno identita}
\eqref{tame-lipschitz-cambio-di-variabile} also hold for $u \circ g$ .

$(iii)$ Part $(ii)$ also holds with $\| \cdot \|_k$ replaced by 
$| \cdot |_{k, \infty}$, and $\Vert \cdot \Vert_{s}^{\Lipg}$ replaced by $\vert \cdot \vert_{s, \infty}^{\Lipg}$,
namely 
\begin{align}  \label{composizione infty}
| u \circ f |_{s, \infty} 
& \leq C (|u|_{s, \infty} + |Dp|_{s-1, \infty} |u|_{1, \infty}), 
\\
\label{composizione infty Lip}
| u \circ f |_{s, \infty}^\Lipg
& \leq C (|u|_{s+1, \infty}^\Lipg + |Dp|_{s-1, \infty}^\Lipg |u|_{2, \infty}^\Lipg).
\end{align}
\end{lemma}

\begin{pf} The bounds 
\eqref{stime-q}, \eqref{tame-cambio-di-variabile} and \eqref{composizione infty} are proved in \cite{Baldi-Benj-Ono}, Appendix B. Let us  prove \eqref{stime-lipschitz-q}. 
Denote $p_\lm(x) := p(\lm,x)$, and similarly for $q_\lm, g_\lm, f_\lm$. 
Since $y = f_\lm(x) = x + p_\lm(x)$ if and only if $x = g_\lm(y) = y + q_\lm(y)$, one has 
\begin{equation}\label{q-p-lambda}
q_\lm(y) + p_\lm(g_\lm(y)) = 0 \, , \quad \forall \lm \in \L, \ y \in \T^d. 
\end{equation}
Let $\lambda_{1}, \lambda_{2} \in \L $, and denote, in short, $q_1 = q_{\lm_1}$, $q_2 = q_{\lm_2}$, and so on.
By \eqref{q-p-lambda},  
\begin{align} 
q_1 - q_2 
& = p_2 \circ g_2 - p_1 \circ g_1  
= (p_2 \circ g_2 - p_1 \circ g_2)
+ (p_1 \circ g_2 - p_1 \circ g_1)
\notag \\ 
& = A_2\inv (p_2 - p_1) 
+ \int_{0}^{1} A_t\inv (D_{x} p_1) \, dt \, (q_2 - q_1) \label{contrazione Lip}
\end{align}
where  $ A_2\inv h := h \circ g_2 $, 
$ A_t\inv h := h \circ \big( g_1 + t [g_2 - g_1] \big)$, $ t \in [0,1]$. 
By \eqref{contrazione Lip},  the $ L^\infty $ norm  
of $(q_2 - q_1)$ satisfies
\[
|q_2 - q_1|_{L^\infty} 
 \leq  
|A_2\inv (p_2 - p_1)|_{L^\infty} 
+ \int_{0}^{1} |A_t\inv (D_{x} p_1)|_{L^\infty} \, dt \, |q_2 - q_1|_{L^\infty} 
\leq  
|p_2 - p_1|_{L^\infty}  
+ \int_{0}^{1} |D_{x} p_1|_{L^\infty}  dt \, |q_2 - q_1|_{L^\infty} 
\]
whence, using the assumption  $|D_{x} p_1|_{L^\infty} \leq 1/2$,   
\begin{equation} \label{norma 0 q2-q1}
|q_2 - q_1|_{L^\infty} \leq 2 |p_2 - p_1|_{L^\infty} \, . 
\end{equation}
By \eqref{contrazione Lip},  using \eqref{tame product infty}, 
the $W^{s,\infty}$ norm of $(q_2 - q_1)$, for $s \geq 0$, satisfies
\[
|q_1 - q_2|_{s, \infty}
\leq |A_2\inv (p_2 - p_1)|_{s, \infty}  
+ \frac32\, \int_{0}^{1} |A_t\inv (D_{x} p_1)|_{L^\infty} \, dt \, |q_2 - q_1|_{s, \infty} 
+ C(s) \int_{0}^{1} |A_t\inv (D_{x} p_1)|_{s, \infty} \, dt \, |q_2 - q_1|_{L^\infty}.
\]
Since $|A_t\inv (D_{x} p_1)|_{L^\infty} = |D_x p_1|_{L^\infty} \leq 1/2$, 
\[
\Big( 1 - \frac34 \Big) |q_2 - q_1|_{s, \infty} 
\leq |A_2\inv (p_2 - p_1)|_{s, \infty}  + C(s) \int_{0}^{1} |A_t\inv (D_{x} p_1)|_{s, \infty} \, dt \, |q_2 - q_1|_{L^\infty}.
\]
Using \eqref{norma 0 q2-q1}, \eqref{composizione infty}, \eqref{interpolation estremi fine-s} and 
\eqref{stime-q},
\[
|q_2 - q_1|_{s, \infty} 
\leq C(s) \Big( |p_2 - p_1|_{s, \infty} + \big\{ \sup_{\lm \in \Lambda} |p_\lm|_{s+1, \infty} \big\} |p_2 - p_1|_{L^\infty} \Big)
\]
and \eqref{stime-lipschitz-q} follows.

{\it Proof of \eqref{cambio di variabile meno identita}}. 
We have $ u \circ f - u = \int_{0}^{1} A_t (D_x u)\,d t\, p  $
where $A_t u (x) := u (x + t p(x)) $, $t \in [0, 1]$.
Then, by  \eqref{mixed norms tame product} and \eqref{tame-cambio-di-variabile},
\begin{eqnarray*}
\Big\| \int_{0}^{1} A_t (D_x u)\,d t\, p \Big\|_{s} & 
\leq_s 
& \int_{0}^{1} \Vert A_t (D_x u) \Vert_{s} \, d t\, |p|_{L^{\infty}} + \int_{0}^{1} \Vert A_t (D_x u) \Vert_{0} \, d t\, |p|_{s, \infty}   
\\
& 
\leq_s 
& \Vert u \Vert_{s + 1} |p|_{L^{\infty}} + |p|_{s, \infty} |p|_{L^\infty} \Vert u \Vert_{2} +  |p|_{s,\infty} \Vert u \Vert_{1} \,,
\end{eqnarray*}
which implies \eqref{cambio di variabile meno identita}. 

{\it Proof of \eqref{tame-lipschitz-cambio-di-variabile}}. 
With the same notation as above, 
\[
u_2 \circ f_2 - u_1 \circ f_1 
= (u_2 \circ f_2 - u_2 \circ f_1) 
+ (u_2 \circ f_1 - u_1 \circ f_1) 
= \int_{0}^{1} A_t (D_{x} u_2) \, dt \, (f_2 - f_1) 
+ A_1 (u_2 - u_1),
\]
where $A_1 h = h \circ f_1$ and $A_t h = h\circ (f_1 + t[f_2 - f_1])$. 
Using \eqref{mixed norms tame product} and \eqref{tame-cambio-di-variabile},
\[
\Big\| \int_{0}^{1} A_t (D_{x} u_2) \, dt \, (f_2 - f_1) \Big\|_s 
\leq_s  
\Big( \| D_x u_2 \|_{s} + \big(  \sup_\lm |D_x p_\lm|_{s-1, \infty} \big) \| D_x u_2 \|_1 \Big) 
|p_2 - p_1|_{L^\infty}
+ \| D_x u_2 \|_{0} |p_2 - p_1|_{s, \infty} 
\]
and 
$ \| A_1 (u_2 - u_1) \|_s \leq_s  \| u_2 - u_1 \|_s + |D_x p_1|_{s-1, \infty} \| u_2 - u_1 \|_1 $.  
Therefore 
\begin{eqnarray*}
\| u_2 \circ f_2 - u_1 \circ f_1 \|_s  & \leq_s &   
|p_2 - p_1|_{L^\infty}  
\Big( \sup_\lm \| u_\lm \|_{s+1} + 
\big( \sup_\lm |p_\lm|_{s, \infty} \big)  \big( \sup_\lm \| u_\lm \|_{2} \big) \Big) \\ 
& & + |p_2 - p_1|_{s, \infty}  \big( \sup_\lm \| u_\lm \|_1 \big) 
+ \| u_2 - u_1 \|_s + \big( \sup_\lm |p_\lm|_{s, \infty} \big) \, \| u_2 - u_1 \|_1 
\end{eqnarray*}
whence \eqref{tame-lipschitz-cambio-di-variabile} follows.
The proof of \eqref{composizione infty Lip} is the same as for \eqref{tame-lipschitz-cambio-di-variabile}, 
replacing all norms $\| \cdot  \|_s$ with $| \cdot  |_{s, \infty}$. 
\end{pf}

\begin{lemma} {\bf (Composition)} \label{lemma astratto composizioni}
Suppose that  for all $\Vert u \Vert_{s_{0}+ \mu_i} \leq 1$ the operator ${\cal Q}_i(u)$ satisfies
\begin{equation}\label{tame Phi-i}
\Vert {\cal Q}_i h \Vert_{s} \leq C(s) \big(\Vert h \Vert_{s + \tau_i} + \Vert u \Vert_{s + \mu_i} \Vert h \Vert_{s_{0}+\tau_i} \big), \quad i = 1, 2.
\end{equation}
Let $\t := {\rm max}\{ \t_1, \t_2 \}$, $\mu := {\rm max}\{ \mu_1, \mu_2 \}$.
Then, for all 
\begin{equation}\label{u-s0-tau-mu}
 \Vert u \Vert_{s_{0}+ \tau + \mu} \leq 1 \, , 
\end{equation}
the composition operator ${\cal Q} := {\cal Q}_{1} \circ {\cal Q}_2$ satisfies the tame estimate
\begin{equation} \label{stima tame astratta composizioni}
\Vert {\cal Q} h \Vert_{s} \leq C(s) \big( \Vert h \Vert_{s + \tau_1 + \tau_2} + \Vert u \Vert_{s + \tau + \mu} \Vert h \Vert_{s_0 + \tau_1 + \tau_2}\big).
\end{equation}
Moreover, if ${\cal Q}_1$, ${\cal Q}_2$, $u$ and $h$ depend in a lipschitz way on a parameter $\lambda$, 
then \eqref{stima tame astratta composizioni} also holds with $\Vert \cdot \Vert_s$ replaced by $\Vert \cdot \Vert_{s}^{\Lipg}$.
\end{lemma}

\begin{pf} 
Apply the estimates for \eqref{tame Phi-i} to $\Phi_1$ first, then to $\Phi_2$, using condition \eqref{u-s0-tau-mu}.  
\end{pf}

\section{Appendix B: proof of Lemmata \ref{lemma:mostro} and \ref{lemma:stime stabilita Phi 12}} 
\label{sec:proofs}

The proof is elementary. 
It is based on a repeated use of the tame estimates of the Lemmata of the Appendix A. 
For  convenience,  we split it into many points.
We remind that $ \mathfrak s_0 := (\nu+2) / 2 $ is fixed 
(it plays the role of the constant $ s_0 $ in Lemma \ref{lemma:standard Sobolev norms properties}).

\smallskip

{\bf Estimates in Step $ 1 $.}

1. --- We prove that $b_3 = b $ defined in  \eqref{c} satisfies the tame estimates
\begin{align} \label{stima b3}
\| b_3 - 1 \|_{s} 
& \leq  \e \, C(s) \big( 1 + \| u \|_{s + 3} \big),
\\
\label{stima derivata b3}
\| \pa_u b_3(u)[h] \|_{s}
& \leq  \e \, C(s) \big( \| h \|_{s+3} + \| u \|_{s+3} \| h \|_{ \mathfrak s_0+3} \big),
\\
\label{stima Lip b3}
\| b_3-1 \|_{s}^\Lipg 
& \leq  \e \, C(s) \big( 1 + \| u \|_{s+3}^\Lipg \big).
\end{align}
{\it Proof of \eqref{stima b3}}.
Write $b_3 = b $ (see  \eqref{c}) as
\begin{equation}  \label{b as}
b_3 - 1 = \psi\big( M [ g(a_3) - g(0) ] \, \big) - \psi(0), 
\quad 
\psi(t) := (1 + t)^{-3}, \quad 
Mh := \frac{1}{2\p}\, \int_\T h \, dx, \quad 
g(t) := (1 + t)^{-\frac13}.
\end{equation}
Thus, for $ \e $ small, 
\[ 
\| b_3 - 1 \|_{s} \leq C(s) \|  M [ g(a_3) - g(0) ] \, \|_s 
\leq C(s) \| g(a_3) - g(0) \|_s 
\leq C(s) \| a_3 \|_s .
\]
In the first inequality we have applied Lemma \ref{lemma:composition of functions, Moser}$(i)$ to the 
function $\psi$, with $u=0$, $ p = 0 $, $h=M[g(a_3)-g(0)]$. 
In the second inequality we have used the trivial fact that $ \| Mh \|_{s} \leq \| h \|_{s} $ for all $h$. 
In the third inequality we have applied again Lemma \ref{lemma:composition of functions, Moser}$(i)$ to the function $g$, with $u=0$, $ p = 0 $, $h=a_3$. 
Finally we estimate  $a_3 $ by \eqref{stima coeff ai 1} with $ s_0 = \mathfrak s_0 $, which holds for $s+2 \leq q$.

\noindent
{\it Proof of \eqref{stima derivata b3}}. 
Using \eqref{b as}, the derivative of $b_3$ with respect to $u$ in the direction $h$ is 
\[
\pa_u b_3(u)[h] = 
\psi' \big( M [ g(a_3) - g(0) ] \big) \, M \big( g'(a_3) \pa_u a_3[h] \, \big).
\]
Then use \eqref{asymmetric tame product}, Lemma \ref{lemma:composition of functions, Moser}$(i)$ applied to the functions $\psi'$ and $g'$, and \eqref{stima coeff ai 2}.

\noindent
{\it Proof of \eqref{stima Lip b3}}. 
It follows from \eqref{stima b3}, \eqref{stima derivata b3} and Lemma \ref{lemma:Lip generale}.

\smallskip

2. --- Using the definition \eqref{primaequazione} of $\rho_0$, 
estimates \eqref{stima b3}, \eqref{stima derivata b3}, \eqref{stima Lip b3} for $b_3$ and estimates \eqref{stima coeff ai 1}, \eqref{stima coeff ai 2}, \eqref{stima coeff ai 3} for $a_3$, one proves that 
$\rho_0$ also satisfies the same estimates \eqref{stima b3}, \eqref{stima derivata b3}, \eqref{stima Lip b3} as $(b_3 - 1)$. 
Since $\b = \pa_x\inv \rho_0$ (see \eqref{defb1b0}), by Lemma \ref{lemma:standard Sobolev norms properties}($i$)
we get 
\begin{align} \label{stima beta}
| \b |_{s, \infty} 
& \leq C(s) \| \b \|_{s + {\mathfrak s}_0} \leq C(s) \| \rho_0 \|_{s + {\mathfrak s}_0} 
\leq  \e \, C(s) \big( 1 + \| u \|_{s + {\mathfrak s}_0 + 3} \big),
\\
\intertext{and, with the same chain of inequalities,}
\label{stima derivata beta}
| \pa_u \b(u)[h] |_{s, \infty} 
& 
\leq  \e \, C(s) \big( \| h \|_{s + {\mathfrak s}_0 +3} + \| u \|_{s + {\mathfrak s}_0 +3} \| h \|_{{\mathfrak s}_0 +3} \big) \, . 
\end{align}
Then Lemma \ref{lemma:Lip generale} implies 
\be \label{stima Lip beta}
| \b |_{s, \infty}^\Lipg
\leq  \e \, C(s) \big( 1 + \| u \|_{s+ {\mathfrak s}_0 +3}^\Lipg \big),
\ee
for all $s + {\mathfrak s}_0 + 3 \leq q$. 
Note that $x \mapsto x + \b(\ph,x)$ is a well-defined diffeomorphism if $|\b|_{1, \infty} \leq 1/2$, and, 
by \eqref{stima beta}, this condition is satisfied provided 
$ \e \, C \big( 1 + \| u \|_{{\mathfrak s}_0 + 4} \big) \leq 1/2$. 

Let $(\ph,y) \mapsto (\ph, y + \tilde\b(\ph,y))$ be the inverse diffeomorphism of  
$(\ph,x) \mapsto (\ph, x + \b(\ph,x))$. 
By Lemma \ref{lemma:utile}($i$) on the torus $\T^{\nu+1}$, $\tilde\b$ satisfies 
\begin{equation} \label{stima beta tilde}
| \tilde\b |_{s, \infty}  \leq C | \b |_{s, \infty}  \stackrel{\eqref{stima beta}} \leq  \e \, C(s) \big( 1 + \| u \|_{s + 3 + {\mathfrak s}_0} \big).
\end{equation}
Writing explicitly the dependence on $u$,  we have $ \tilde\b(\ph,y;u) + \b\big( \ph, \, y + \tilde\b(\ph,y;u) ; u \big) = 0 $. 
Differentiating the last equality with respect to $u$ in the direction $h$ gives
\[
(\pa_u \tilde\b)[h] = 
-  \inv \Big( \frac{ \pa_u \b [h] }{ 1+\b_x } \Big),
\]
therefore, applying Lemma \ref{lemma:utile}$(iii)$ to deal with ${\cal A}\inv$, 
\eqref{tame product infty} for the product $(\pa_u \b [h]) (1 + \b_x)\inv$, 
the estimates \eqref{stima beta}, \eqref{stima derivata beta}, \eqref{stima Lip beta} for $\b$, 
and \eqref{interpolation estremi fine} (with $ a_0 = \mathfrak s_0 + 3 $, $ b_0 =  \mathfrak s_0 + 4 $, $ p = 1 $, 
$ q = s - 1$),  we obtain (for  $s + {\mathfrak s}_0 + 4 \leq q$)
\be\label{stima derivata beta tilde}
| \pa_u \tilde\b(u)[h] |_{s, \infty} 
\leq  \e \, C(s) \big( \| h \|_{s+3 + {\mathfrak s}_0} + \| u \|_{s+ 4 + {\mathfrak s}_0} \| h \|_{3 + {\mathfrak s}_0} \big) \, . 
\ee
Then, using  Lemma \ref{lemma:Lip generale} with $ p = 4 + \mathfrak s_0 $, the bounds 
\eqref{stima beta tilde}, \eqref{stima derivata beta tilde} imply 
\be
\label{stima Lip beta tilde}
| \tilde\b |_{s, \infty}^\Lipg 
\leq  \e \, C(s) \big( 1 + \| u \|_{s+ 4 + {\mathfrak s}_0}^\Lipg \big).
\ee

\smallskip

3. --- {\sc Estimates of ${\cal A}(u)$ and ${\cal A}(u)^{-1}$. } 
By \eqref{tame-cambio-di-variabile}, \eqref{stima beta} and \eqref{stima beta tilde},
\begin{equation}
\Vert {\cal A}(u)h \Vert_{s} + \| {\cal A}(u)\inv h \|_{s} 
\leq \, C(s) \big(\| h \|_{s} + \| u \|_{s + {\mathfrak s}_0 + 3} \| h \|_1 \big). 
\label{stima A e sua inversa} 
\end{equation}
Moreover, by \eqref{tame-lipschitz-cambio-di-variabile},  
\eqref{stima Lip beta} and \eqref{stima Lip beta tilde},
\begin{equation}\label{stima A e sua inversa lipschitz}
\|{\cal A}(u) h \|_{s}^\Lipg + \|{\cal A}(u)\inv h \|_{s}^\Lipg 
\leq C(s) \big(\|h\|_{s + 1}^\Lipg + \| u \|_{s + {\mathfrak s}_0 + 4}^\Lipg \|h \|_{2}^\Lipg \big).
\end{equation}

Since ${\cal A}(u)g(\vphi , x) =  g(\vphi , x + \beta(\vphi , x; u))$, 
the derivative of ${\cal A}(u)g$ with respect to $u$ in the direction $h$ is the product
$ \partial_u \big( {\cal A}(u)g \big)[h]\,=\, ( {\cal A}(u)g_{x} ) \, \partial_u\beta(u)[h] $. 
Then, by \eqref{mixed norms tame product}, 
\eqref{stima derivata beta} and \eqref{stima A e sua inversa}, 
\begin{equation}\label{stima derivata A}
\| \partial_u ({\cal A}(u)g) [h] \, \|_s 
\leq \e C(s) \Big(\| g \|_{s + 1} \| h \|_{ {\mathfrak s}_0 + 3} + \| g \|_{2} \| h \|_{s + {\mathfrak s}_0 + 3}  
+ \| u \|_{s + {\mathfrak s}_0 + 3}  \| g \|_{2} \| h \|_{{\mathfrak s}_0 + 3}\Big).
\end{equation}
Similarly 
$\partial_u ({\cal A}(u)\inv g ) [h] = ({\cal A}(u)\inv g_x) \, \partial_u \tilde{\beta}(u)[h]$, therefore \eqref{mixed norms tame product}, \eqref{stima derivata beta tilde}, \eqref{stima A e sua inversa} imply that
\begin{equation}\label{stima derivata inversa A}
\| \partial_u ( {\cal A}^{-1}(u)g ) [h] \, \|_s 
\leq \e C(s) \Big(\| g \|_{s + 1} \| h \|_{{\mathfrak s}_0 + 3} + \| g \|_{2} \| h \|_{s + {\mathfrak s}_0 + 3}  + 
\| u \|_{s + {\mathfrak s}_0 + 4}\| g \|_{2} \| h \|_{{\mathfrak s}_0 + 3} \Big).
\end{equation}

\smallskip

4. ---  
The coefficients $b_0, b_1, b_2$ are given in 
\eqref{b1 b3}, \eqref{b0 b2}. 
By \eqref{mixed norms tame product}, \eqref{stima A e sua inversa}, \eqref{palla Lip di sicurezza}, \eqref{stima beta} and 
\eqref{stima coeff ai 1},  
\begin{equation}\label{stime bi}
\| b_i \|_s \leq \e C(s) ( 1 + \| u \|_{s + {\mathfrak s}_0 + 6} ), \quad i = 0 , 1 , 2.
\end{equation}
Moreover, in analogous way, by \eqref{mixed norms tame product}, 
\eqref{stima A e sua inversa lipschitz}, \eqref{palla Lip di sicurezza}, \eqref{stima Lip beta} and \eqref{stima coeff ai 3},
\begin{equation}\label{stime bi lipschitz}
\| b_i \|_{s}^\Lipg \leq \e C(s) (1 + \| u \|_{s + {\mathfrak s}_0 + 7}^\Lipg ), \quad i = 0 , 1 , 2.
\end{equation}
Now we estimate the derivative with respect to $u$ of $b_1$. 
The estimates for $b_0$ and $b_2$ are analogous.
By \eqref{b1 b3} we write $b_1(u) = {\cal A}(u)\inv b_1^*(u)$ 
where
$ b_1^* := $ $ \o\cdot\partial_{\vphi}\b + $ $ (1 + a_3)\beta_{xxx} + $ $ a_2 \beta_{xx} + $ $ a_1(1 + \beta_x) $. 
The bounds \eqref{stima coeff ai 2}, \eqref{stima derivata beta}, \eqref{stima beta}, \eqref{palla Lip di sicurezza}, 
and \eqref{mixed norms tame product} imply that
\begin{equation}\label{stima b1*}
\| \partial_u b_{1}^{*}(u)[h]\|_{s} \leq \e C(s) \big(\| h \|_{s + {\mathfrak s}_0 + 6} + 
\| u \|_{s + {\mathfrak s}_0 + 6} \| h \|_{{\mathfrak s}_0 + 6} \big)\,.
\end{equation}
Now, 
\begin{equation}\label{marcellino}
\partial_u b_1(u)[h] 
= \partial_u \big({\cal A}(u)\inv b_1^{*}(u) \big) [h] 
= (\partial_u {\cal A}(u)\inv) (b_1^*(u)) [h] + {\cal A}(u)^{-1} (\partial_{u}b_1^{*}(u)[h]).
\end{equation}
Then \eqref{asymmetric tame product}, 
\eqref{marcellino}, 
\eqref{stima A e sua inversa}, 
\eqref{stima derivata inversa A}, \eqref{interpolation estremi fine} 
(with $ a_0 = \mathfrak s_0 +4 $, $ \b_0 = \mathfrak s_0  + 6 $, $ p = s- 1$, $ q = 1 $) 
\eqref{stima b1*} imply
\begin{eqnarray}
\| \partial_u {\cal A}(u)\inv (b_1^*(u)) [h] \|_s 
& \leq & \e C(s) \big( \| h \|_{s + {\mathfrak s}_0 + 3} + \| u \|_{s + {\mathfrak s}_0 + 7} \| h \|_{{\mathfrak s}_0 + 3} \big) 
\label{primo pezzo derivata b1} \\
\| {\cal A}(u)^{-1} \partial_{u} b_1^{*}(u)[h] \|_{s} 
& \leq & \e C(s) \big( \| h \|_{s + {\mathfrak s}_0 + 6} + \| u \|_{s + {\mathfrak s}_0 + 6} \| h \|_{ {\mathfrak s}_0 + 6} \big). 
\label{secondo pezzo derivata b1} 
\end{eqnarray}
Finally \eqref{marcellino}, \eqref{primo pezzo derivata b1} and \eqref{secondo pezzo derivata b1} imply 
\begin{equation}\label{stima finale b1}
\| \partial_u b_1(u)[h] \|_s \leq \e C(s) \big(\| h \|_{s + {\mathfrak s}_0 + 6} + \| u \|_{s + {\mathfrak s}_0 + 7} 
\| h \|_{{\mathfrak s}_0 + 6}  \big),
\end{equation}
which holds for all $s + {\mathfrak s}_0 + 7 \leq q$.

\smallskip

{\bf Estimates in  Step 2. }

5. --- We prove that the coefficient $\muff_3$, defined in \eqref{mu 3}, satisfies the following estimates:
\begin{eqnarray}
|\muff_3 - 1 | \, ,  
 |\muff_3 - 1 |^\Lipg &\leq& \e C  \label{stima mu3-1 Lip}\\
|\partial_u \muff_3(u)[h]| & \leq & \e C \Vert h \Vert_{{\mathfrak s}_0 + 3} .
\label{stima derivata mu 3}
\end{eqnarray}
Using \eqref{mu 3} \eqref{stima b3}, \eqref{palla Lip di sicurezza}
$$
|\muff_3 - 1|  \leq \frac{1}{(2\pi)^{\nu}} \int_{\T^{\nu}}|b_3 - 1| \, d\vphi 
\leq C \Vert b_3 - 1\Vert_{{\mathfrak s}_0}  \leq \e C.
$$
Similarly we get the Lipschitz part of \eqref{stima mu3-1 Lip}. 
The estimate \eqref{stima derivata mu 3} follows by  \eqref{stima derivata b3}, since
$$
| \partial_u \muff_3(u)[h] \, | 
\leq \frac{1}{(2\pi)^{\nu}} \int_{\T^{\nu}}\vert\partial_{u}b_3(u)[h] \vert \, d\vphi 
\leq C \Vert \partial_{u}b_3(u)[h]\Vert_{{\mathfrak s}_0} 
\leq \e C \Vert h \Vert_{{\mathfrak s}_0 + 3}.
$$
6. --- \textsc{Estimates of $\alpha$. } The function
$\alpha(\ph)$, defined in \eqref{alpha}, satisfies 
\begin{eqnarray}
\vert \alpha \vert_{s, \infty} & \leq & \e \gamma_0^{-1} \, C(s) 
\big(1 + \Vert u \Vert_{s + \tau_0 + {\mathfrak s}_0 + 3} \big) \label{stima  alpha}
\\
 \vert \alpha \vert_{s, \infty}^\Lipg
& \leq & \e \gamma_0^{-1} \, C(s) 
\big(1 + \Vert u \Vert_{s + \tau_0 + {\mathfrak s}_0 + 3}^\Lipg \big) \label{stima Lip alpha}
\\
\vert \partial_u \alpha(u)[h] \vert_{s, \infty} 
& \leq & \e \gamma_0^{-1}  \, C(s) \big(\Vert h \Vert_{s + \tau_0 + {\mathfrak s}_0 + 3} + \Vert u \Vert_{s + \tau_0 
+ {\mathfrak s}_0 + 3}\Vert h \Vert_{ {\mathfrak s}_0 + 3} \big). 
\label{stima derivata alpha}
\end{eqnarray}
Remember that $\om = \lm \bar\om$, and $|\bar\om \cdot l| \geq 3 \g_0 |l|^{-\t_0}$, $ \forall l \neq 0$, see \eqref{omdio}. 
By \eqref{stima b3} and \eqref{stima mu3-1 Lip},
$$
| \alpha |_{s, \infty} 
\leq \| \alpha \|_{s + \mathfrak s_0} \leq
C \gamma_0^{-1} \Vert b_3 - \muff_3 \Vert_{s + \mathfrak s_0 + \tau_0} 
\leq C(s) \gamma_0^{-1} \e (1 + \Vert u \Vert_{s +\tau_0 + \mathfrak s_0  + 3} )
$$
proving \eqref{stima  alpha}. Then \eqref{stima Lip alpha} holds similarly using \eqref{stima Lip b3} and
$ (\ompaph)\inv = \lm^{-1} \, (\bar\om \cdot \pa_\ph)\inv $. 
Differentiating formula \eqref{alpha} with respect to $u$ in the direction $h$ gives
$$
\partial_{u}\alpha(u)[h] = (\lambda \bar\om \cdot \partial_{\vphi})^{-1} \Big(\frac{\partial_u b_3(u)[h] \muff_3 - b_3 \partial_{u}\muff_3(u)[h]}{\muff_3^{2}}\Big)
$$
then, the standard Sobolev embedding,  \eqref{stima b3}, \eqref{stima derivata b3}, \eqref{stima mu3-1 Lip}, 
\eqref{stima derivata mu 3}  imply 
\eqref{stima derivata alpha}. 
Estimates \eqref{stima Lip alpha} and \eqref{stima derivata alpha} hold for $s + \t_0 + {\mathfrak s}_0 + 3 \leq q$.
Note that \eqref{cambio2} is a well-defined diffeomorphism if $|\a|_{1, \infty} \leq 1/2$, and, by \eqref{stima Lip alpha}, this holds by
\eqref{palla di sicurezza}. 

\smallskip

7. --- {\sc Estimates of $\tilde{\alpha}$. } 
Let $\th \rightarrow \th + \o \tilde{\alpha}(\th)$ be the inverse change of variable of \eqref{cambio2}.   
The following estimates hold: 
\begin{eqnarray}
\vert \tilde{\alpha} \vert_{s, \infty} 
& \leq & \e\gamma_0^{-1} \, C(s) \big(1 + \Vert u \Vert_{s + \tau_0 + {\mathfrak s}_0 + 3}  \big) \label{stima alpha tilde}\\
\vert \tilde{\alpha} \vert_{s, \infty}^\Lipg 
& \leq & \e\gamma_0^{-1} \, C(s) \big(1 + \Vert u \Vert_{s + \tau_0 + {\mathfrak s}_0 + 4}^\Lipg \big) \label{stima Lip alpha tilde}\\
\vert \partial_{u}\tilde{\alpha}(u)[h] \vert_{s, \infty}
& \leq & \e\gamma_0^{-1} \, C(s) \big(\Vert h \Vert_{s + \tau_0 + {\mathfrak s}_0 + 3} + \Vert u \Vert_{s + \tau_0 
+ {\mathfrak s}_0 + 4} \Vert h \Vert_{\tau_0 + {\mathfrak s}_0 + 3} \big). 
\label{stima derivata alpha tilde}
\end{eqnarray}
The bounds \eqref{stima alpha tilde}, \eqref{stima Lip alpha tilde} follow  by  
\eqref{stime-q}, \eqref{stima alpha}, and  \eqref{stime-lipschitz-q}, \eqref{stima Lip alpha}, respectively.
To estimate the partial derivative of $\tilde{\alpha}$ with respect to $u$ we differentiate the identity 
 $ \tilde{\alpha}(\vartheta ; u) + \alpha(\vartheta + \o \tilde{\alpha}(\vartheta; u); u) = 0 $, 
which gives 
$$
 \partial_u \tilde{\alpha}(u)[h] = - B^{-1} \Big(\frac{\partial_u \alpha[h]}{1 + \o\cdot\partial_{\vphi}\alpha} \Big).
 $$
Then applying Lemma \ref{lemma:utile}$(iii)$ to deal with $B^{-1}$, 
\eqref{tame product infty} for the product $\partial_u \alpha[h] \, (1 + \o\cdot \partial_{\vphi}\alpha )^{-1}$, and estimates \eqref{stima Lip alpha}, \eqref{stima derivata alpha},  \eqref{interpolation estremi fine},  
we obtain \eqref{stima derivata alpha tilde}.
 
\smallskip

8. ---  
The transformations $B(u)$ and $B(u)^{-1}$, defined in \eqref{operatore2} resp. \eqref{B-1}, satisfy the following estimates:
\begin{eqnarray}
\Vert B(u) h \Vert_s + \Vert B(u)^{-1} h \Vert_s 
\!\!\!\! & \leq & \!\!  \!\!C(s) \big(\Vert h \Vert_s + \Vert u \Vert_{s + \tau_0 + {\mathfrak s}_0 + 3} \Vert h \Vert_1 \big) \label{stima B e sua inversa}
\\
\Vert B(u) h \Vert_s^\Lipg + \Vert B(u)^{-1} h \Vert_s^\Lipg  
\!\! \!\!&\leq &  \!\!\!\! C(s) \big(\Vert h \Vert_{s + 1}^\Lipg + \Vert u \Vert_{s + \tau_0 + {\mathfrak s}_0 + 4}^\Lipg \Vert h \Vert_2^\Lipg \big) \label{stima Lip B e sua inversa}
\\
\Vert \partial_u (B(u)g) [h] \Vert_s 
\!\! \!\! &  \leq & \!\!\!\!
C(s) \big(\Vert g\Vert_{s + 1} \Vert h \Vert_{\s_0} 
+ \Vert g\Vert_{1} \Vert h \Vert_{s + \s_0} 
+ \Vert u \Vert_{s + \s_0} \Vert g\Vert_{2} \Vert h \Vert_{ \s_0} \big) \label{stima derivata B}
\\
\Vert \partial_u ( B(u)^{-1}g ) [h] \Vert_s 
\!\!\!\! &\leq& \!\!\!\!
C(s) \big( \Vert g\Vert_{s + 1} \Vert h \Vert_{\s_0} 
+ \Vert g \Vert_{1} \Vert h \Vert_{s +\s_0} 
+ \Vert u \Vert_{s + \s_0 + 1} \Vert g\Vert_{2} \Vert h \Vert_{\s_0} \big)
\label{stima derivata B  inversa}
\end{eqnarray}
where $ \s_0 := \t_0 + {\mathfrak s}_0 +3 $. 
Estimates \eqref{stima B e sua inversa} and \eqref{stima Lip B e sua inversa} follow 
by Lemma \ref{lemma:utile}$(ii)$ and  
\eqref{stima alpha}, \eqref{stima alpha tilde}, \eqref{stima Lip alpha},  \eqref{stima Lip alpha tilde}.
The derivative of $B(u)g$ with respect to $u$ in the direction $h$ is the product $f z$ where 
$f := B(u)(\o\cdot\partial_{\vphi}g)$ and $z := \partial_{u}\alpha(u)[h]$. 
By \eqref{mixed norms tame product}, $\| fz \|_s \leq C(s) ( \| f \|_s |z|_{L^\infty} + \| f \|_0 |z|_{s, \infty})$.
Then \eqref{stima derivata alpha}, \eqref{stima B e sua inversa} imply \eqref{stima derivata B}. 
In analogous way, \eqref{stima derivata alpha tilde} and \eqref{stima B e sua inversa} give \eqref{stima derivata B inversa}.

\smallskip

9. --- {\sc estimates of $\rho$. }
The function $\rho$ defined in \eqref{anche def rho}, namely $\rho = 1 + B\inv (\ompaph \a)$, satisfies 
\begin{eqnarray}
\vert \rho -1 \vert_{s, \infty} 
& \leq & \e \g_0 \inv \, C(s) ( 1 + \Vert u \Vert_{s +\tau_0 + {\mathfrak s}_0 + 4} ) \label{stima rho} 
\\
\vert \rho -1 \vert_{s, \infty}^\Lipg 
& \leq & \e \g_0 \inv \, C(s) (1 + \Vert u \Vert_{s + \tau_0 + {\mathfrak s}_0 + 5}^\Lipg ) \label{stima Lip rho}
\\
\Vert \partial_u \rho(u)[h] \, \Vert_{s} 
& \leq & 
\e \g_0 \inv \, C(s ) \big( \Vert h \Vert_{s + \tau_0 + {\mathfrak s}_0 + 4} + \Vert u \Vert_{s + \tau_0 + {\mathfrak s}_0 + 5} \Vert h \Vert_{ \tau_0 + {\mathfrak s}_0 + 4} \big). 
\label{stima derivata rho}
\end{eqnarray}
The bound \eqref{stima rho} follows by 
\eqref{anche def rho}, \eqref{composizione infty}, \eqref{stima alpha}, \eqref{palla di sicurezza}. 
Similarly \eqref{stima Lip rho} follows by \eqref{composizione infty Lip},  
\eqref{stima Lip alpha} and \eqref{palla Lip di sicurezza}.
Differentiating \eqref{anche def rho} with respect to $u$ in the direction $h$ we obtain
 $$
  \partial_{u}\rho(u)[h]\,=\,\partial_{u}B(u)^{-1}(\o\cdot\partial_{\vphi}\alpha) [h] + B(u)^{-1}\big(\o\cdot\partial_{\vphi}(\partial_u\alpha(u)[h]) \big).
 $$
By \eqref{stima derivata B inversa}, \eqref{stima alpha}, and \eqref{palla di sicurezza}, we get
  \begin{equation}\label{stima primo pezzo derivata rho}
 \Vert \partial_{u}B(u)^{-1}(\o\cdot\partial_{\vphi}\alpha) [h] \Vert_{s} 
 \leq \e \g_0 \inv \, C(s) \big( \Vert h \Vert_{s + \tau_0 + {\mathfrak s}_0 + 3} + \Vert u \Vert_{s + \tau_0 + {\mathfrak s}_0 + 5} \Vert h \Vert_{ \tau_0 + {\mathfrak s}_0 + 3} \big).
  \end{equation}
Using \eqref{stima B e sua inversa}, \eqref{stima derivata alpha}, \eqref{palla di sicurezza},
 and applying \eqref{interpolation estremi fine}, one has
  \begin{equation}\label{stima secondo pezzo derivata rho}
  \Vert B(u)^{-1}\big(\o\cdot\partial_{\vphi}(\partial_u\alpha(u)[h]) \big) \Vert_{s} 
  \leq \e \g_0 \inv \, C(s) \big(\Vert h \Vert_{s + \tau_0 + {\mathfrak s}_0 + 4} + \Vert u \Vert_{s + \tau_0 + {\mathfrak s}_0 + 4}\Vert h \Vert_{ \tau_0 + {\mathfrak s}_0 + 4} \big)\,.
    \end{equation}
Then \eqref{stima primo pezzo derivata rho} and \eqref{stima secondo pezzo derivata rho} imply \eqref{stima derivata rho}, for all $s + \t_0 + {\mathfrak s}_0 + 5 \leq q$.

\smallskip

10. ---   
The coefficients $c_0$, $c_1$, $c_2$ defined in \eqref{coefficienti mL2} satisfy the following estimates: for $ i = 0,1,2 $, $ s 
\geq \mathfrak s_0 $, 
\begin{eqnarray}
\Vert c_i \Vert_{s} & \leq & \e C(s) \big(1 + \Vert u \Vert_{s + \tau_0 + {\mathfrak s}_0 + 6} \big), 
\label{stime ci}
\\
\Vert c_i \Vert_{s}^\Lipg & \leq & \e C(s) \big(1 + \Vert u \Vert_{s + \tau_0 + {\mathfrak s}_0 + 7}^\Lipg \big),
 \label{stime Lip ci}
\\
\Vert \partial_u c_i [h] \Vert_s & \leq & \e C(s ) 
\big(\Vert h \Vert_{s + \tau_0 + {\mathfrak s}_0 + 6} + 
\Vert u \Vert_{s + \tau_0 + {\mathfrak s}_0 + 7} \Vert h \Vert_{  \tau_0 + 2{\mathfrak s}_0 + 6} \big) \, . 
\label{stime derivate ci}
\end{eqnarray}
The definition of $c_i$ in \eqref{coefficienti mL2}, \eqref{mixed norms tame product}, 
\eqref{palla di sicurezza}, \eqref{stima B e sua inversa}, \eqref{stima rho}, \eqref{stime bi} 
and $\e \g_0\inv < 1 $, imply \eqref{stime ci}.  
Similarly \eqref{palla Lip di sicurezza}, \eqref{stima Lip B e sua inversa}, 
\eqref{stima Lip rho} and \eqref{stime bi lipschitz} imply \eqref{stime Lip ci}. 
Finally \eqref{stime derivate ci} follows from differentiating the formula of $c_i(u)$ and using 
\eqref{palla di sicurezza},
\eqref{stime bi}, 
\eqref{stima derivata B inversa}, 
\eqref{stima B e sua inversa},
 \eqref{asymmetric tame product}-\eqref{mixed norms tame product},
\eqref{stima rho}, 
\eqref{stima derivata rho}. 

\smallskip

  {\bf Estimates in the step 3.}
  
\smallskip

11. ---  
The function $v$ defined in \eqref{descent ordine zero} satisfies the following estimates:
\begin{eqnarray}
\Vert v - 1 \Vert_s & \leq & \e C(s)\big(1 + \Vert u \Vert_{s + \tau_0 + {\mathfrak s}_0 + 6} \big) 
\label{stima v}
\\
\Vert v - 1 \Vert_{s}^\Lipg & \leq & \e C(s) \big(1 + \Vert u \Vert_{s + \tau_0 + {\mathfrak s}_0 + 7}^\Lipg \big) 
\label{stima 
Lip v}
\\
\Vert \partial_u v[h] \Vert_s & \leq & \e C(s) \big(\Vert h \Vert_{s + \tau_0 + {\mathfrak s}_0 + 6} + \Vert u \Vert_{s + \tau_0 + {\mathfrak s}_0 + 7} \Vert h \Vert_{ \tau_0 + 2 {\mathfrak s}_0 + 6} \big) 
\label{stima derivata v}
\end{eqnarray}
In order to prove \eqref{stima v}  
we apply the Lemma 
\ref{lemma:composition of functions, Moser}$(i)$ with $f(t) := \exp(t) $ (and $ u = 0 $, $ p = 0 $): 
$$
\Vert v - 1 \Vert_s 
= \Big\| f\Big(- \frac{\partial_y^{-1}c_2}{3 \muff_3}\Big) - f(0) \Big\|_s 
\stackrel{\eqref{stima mu3-1 Lip}} \leq C \Vert c_2 \Vert_s \stackrel{\eqref{stime ci}} \leq  
\e C(s)\big(1 + \Vert u \Vert_{s + \tau_0 + {\mathfrak s}_0 + 6} \big) \, . 
$$
Similarly \eqref{stima Lip v} follows.
Differentiating formula \eqref{descent ordine zero} we get  
   $$
   \partial_u v [h] = - f'\Big(- \frac{\partial_y^{-1}c_2}{3 \muff_3} \Big)
\left\{\frac{1}{3 \muff_3}\partial_u \Big(\partial_y^{-1} c_2 \Big)[h] - \frac{\partial_y^{-1}c_2 \partial_u \muff_3[h]}{3 \muff_3^{2}} \right\} .
   $$
Then using \eqref{palla di sicurezza}, \eqref{asymmetric tame product}, 
Lemma \ref{lemma:composition of functions, Moser}$(i)$ applied to $f' = f$, 
and the estimates \eqref{stime ci}, 
\eqref{stime derivate ci}, 
\eqref{stima mu3-1 Lip} 
and \eqref{stima derivata mu 3} we get \eqref{stima derivata v}.

\smallskip

12. --- The multiplication operator $ {\cal M}  $ defined in \eqref{cambio3} and its inverse $ {\cal M}  \inv$ (which is
the multiplication operator by $v^{-1}$) both satisfy
\begin{align} 
\label{stima Phi}
\| {\cal M} ^{\pm 1} h \|_s 
& \leq C(s) \big( \| h \|_s + \|u \|_{s + \tilde\s} \| h \|_{\mathfrak s_0} \big), 
\\
\label{stima Lip Phi}
\| {\cal M} ^{\pm 1} h \|_s^\Lipg 
& \leq C(s) \big( \| h \|_s^\Lipg + \|u \|_{s + \tilde\s +1}^\Lipg \| h \|_{\mathfrak s_0}^\Lipg \big), 
\\
\label{stima derivata Phi}
\| \pa_u {\cal M} ^{\pm 1} (u) g [h] \|_s 
& \leq 
\e C(s) \big( \Vert g\Vert_{s} \Vert h \Vert_{\mathfrak s_0 + \tilde\s} 
+ \Vert g \Vert_{\mathfrak s_0} \Vert h \Vert_{s + \tilde\s} 
+ \Vert u \Vert_{s + \tilde\s +1} \Vert g\Vert_{\mathfrak s_0} \Vert h \Vert_{\mathfrak s_0 + \tilde\s} \big),
\end{align}
with  $\tilde\s := \t_0 + {\mathfrak s}_0 + 6$.

The inequalities  \eqref{stima Phi}-\eqref{stima derivata Phi}
follow by \eqref{palla di sicurezza}, \eqref{palla Lip di sicurezza}, \eqref{asymmetric tame product}, 
\eqref{stima v}-\eqref{stima derivata v}.

\smallskip

13. --- The coefficients $d_1, d_0$, defined in \eqref{mL3}, satisfy, for $ i = 0,1 $
\begin{align} 
\| d_i \|_s 
& \leq \e C(s) (1 + \| u \|_{s + \t_0 + {\mathfrak s}_0 + 9}), 
\label{stima di}
\\
\| d_i \|_s^\Lipg 
& \leq \e C(s) (1 + \| u \|_{s + \t_0 + {\mathfrak s}_0 + 10}^\Lipg), 
\label{stima Lip di}
\\
\| \pa_u d_i(u)[h] \|_s 
& \leq \e C(s) \big(\| h \|_{s + \tau_0 + {\mathfrak s}_0 + 9} 
+ \| u \|_{s + \tau_0 + {\mathfrak s}_0 + 10} \| h \|_{ \tau_0 +2 {\mathfrak s}_0 + 9} \big), 
\label{stima derivata di}
\end{align}
by \eqref{asymmetric tame product},
\eqref{palla di sicurezza}, \eqref{palla Lip di sicurezza},  \eqref{stime ci}-\eqref{stime derivate ci} and 
\eqref{stima v}-\eqref{stima derivata v}. 

\smallskip

{\bf Estimates in the Step 4.}

\smallskip

14. --- The constant $\muff_1$ defined in  \eqref{m-p} satisfies 
\begin{equation}\label{stime mu1}
\vert \muff_1 \vert + \vert \muff_1 \vert^\Lipg  \leq \e C, \quad 
\vert \partial_u \muff_1(u)[h] \vert \leq \e C \Vert h \Vert_{ \tau_0 + 2 {\mathfrak s}_0 + 9} \, , 
\end{equation}
by 
\eqref{palla Lip di sicurezza},
\eqref{stima di}-\eqref{stima derivata di}.

 \smallskip

15. --- The function $p(\th)$ defined in \eqref{def:p} satisfies the following estimates: 
 \begin{eqnarray}
 \vert p \vert_{s, \infty} & \leq & \e\gamma_{0}^{-1} C(s) ( 1 + \Vert u \Vert_{s + 2\tau_0 + 2{\mathfrak s}_0 + 9})\label{stima p}\\
 \vert p \vert_{s, \infty}^\Lipg & \leq & \e\gamma_{0}^{-1} C(s) ( 1 + \Vert u \Vert_{s + 2\tau_0 + 2 {\mathfrak s}_0 + 10}^\Lipg)\label{stima Lip p}\\
 \vert \partial_u p(u)[h] \vert_{s, \infty} & \leq & \e\gamma_{0}^{-1} C(s) \big(\Vert h \Vert_{s + 2\tau_0 + 2{\mathfrak s}_0 + 9} + \Vert u \Vert_{s + 2\tau_0 + 2{\mathfrak s}_0 + 10} \Vert h \Vert_{\tau_0 + 2{\mathfrak s}_0 + 9} \big). 
 \label{stima derivata p}
 \end{eqnarray}
which follow by
\eqref{stima di}-\eqref{stima derivata di} and \eqref{stime mu1}
applying the same argument  used in the proof of \eqref{stima Lip alpha}. 
 
\smallskip

16. --- The operators $ {\cal T} $, ${\cal T}^{-1} $  defined in \eqref{MM-1} satisfy  
\begin{eqnarray}
\Vert {\cal T}^{\pm1} h \Vert_{s} 
& \leq & C(s) \big(\Vert h \Vert_s + \Vert u \Vert_{s + \bar\s} \Vert h \Vert_1 \big)
\label{stima M e sua inversa}
\\
\Vert {\cal T}^{\pm1} h \Vert_{s}^\Lipg 
& \leq & C(s) \big(\Vert h \Vert_{s + 1}^\Lipg + \Vert u \Vert_{s + \bar\s +1}^\Lipg \Vert h \Vert_2^\Lipg \big)
\label{stima Lip M e sua inversa}
\\
\Vert \partial_u({\cal T}^{\pm1} (u)g)[h] \Vert_s 
& \leq &
\e\gamma_0^{-1} \, C(s)\big(\Vert g \Vert_{s+1}\Vert h \Vert_{\bar\s} + \Vert g \Vert_1\Vert h \Vert_{s + \bar\s} + \Vert u \Vert_{s + \bar\s + 1}\Vert g \Vert_2 \Vert h \Vert_{ \bar\s} \big),
\label{stima derivata M e sua inversa}
\end{eqnarray}
with  $ \bar\s:= 2\tau_0 + 2{\mathfrak s}_0 + 9 $.  
The estimates \eqref{stima M e sua inversa} and \eqref{stima Lip M e sua inversa} follow by \eqref{tame-cambio-di-variabile}, 
\eqref{tame-lipschitz-cambio-di-variabile} and using \eqref{stima p} and \eqref{stima Lip p}.
The derivative $ \partial_u( {\cal T} (u)g)[h] $ is the product
$ ({\cal T}(u) g_y)\,\partial_u p(u)[h]$. Hence \eqref{mixed norms tame product},
\eqref{stima M e sua inversa} and \eqref{stima derivata p} imply \eqref{stima derivata M e sua inversa}.

\smallskip

17. --- The coefficients $e_0$, $e_1$, 
defined in \eqref{e0 e1}, satisfy the following estimates: for  $ i = 0, 1 $
\begin{eqnarray}
\Vert e_i \Vert_s & \leq & \e C(s)  (1 + \Vert u \Vert_{s + 2\tau_0 + 2{\mathfrak s}_0 + 9}),
\label{stima ei}
\\
\Vert e_i \Vert_s^\Lipg & \leq & \e C(s) (1 + \Vert u \Vert_{s + 2\tau_0 + 2{\mathfrak s}_0 + 10}^\Lipg), 
\label{stima Lip ei}
\\
\Vert \partial_u e_i(u)[h] \Vert_s & \leq & \e C(s)\big(\Vert h \Vert_{s + 2\tau_0 + 2{\mathfrak s}_0 + 9} + \Vert u \Vert_{s + 2\tau_0 + 2{\mathfrak s}_0 + 10} \Vert h \Vert_{ 2\tau_0 + 2{\mathfrak s}_0 + 9} \big) \, . 
\label{stima derivata ei}
\end{eqnarray}
The estimates \eqref{stima ei}, \eqref{stima Lip ei} follow 
 by \eqref{palla di sicurezza}, \eqref{palla Lip di sicurezza}, \eqref{equazione omologica step 4}, \eqref{stima di}, 
\eqref{stima Lip di}, \eqref{stima M e sua inversa} and \eqref{stima Lip M e sua inversa}. 
The estimate \eqref{stima derivata ei} follows 
differentiating the formulae of $e_0$ and $e_1$ in \eqref{e0 e1}, and applying \eqref{stima di}, \eqref{stima derivata di}, 
\eqref{stima M e sua inversa} and \eqref{stima derivata M e sua inversa}.

\smallskip

{\bf Estimates in the Step 5.}

\smallskip

18. --- The function $w$ defined in \eqref{w} satisfies the following estimates:
\begin{eqnarray}
\Vert w \Vert_s & \leq & \e C(s)  (1 + \Vert u \Vert_{s + 2\tau_0 + 2{\mathfrak s}_0 + 9})\label{stima w}\\
\Vert w \Vert_s^\Lipg & \leq & \e C(s) (1 + \Vert u \Vert_{s + 2\tau_0 + 2{\mathfrak s}_0 + 10}^\Lipg)\label{stima Lip w}\\
\Vert \partial_u w(u)[h] \Vert_s & \leq & \e C(s)\big(\Vert h \Vert_{s + 2\tau_0 + 2{\mathfrak s}_0 + 9} + \Vert u \Vert_{s + 2\tau_0 + 2{\mathfrak s}_0 + 10} \Vert h \Vert_{ 2\tau_0 + 2{\mathfrak s}_0 + 9} \big) \label{stima derivata w}
\end{eqnarray}
which follow  by 
\eqref{stima mu3-1 Lip}, 
\eqref{stima derivata mu 3}, 
\eqref{stime mu1}, 
\eqref{stima ei}-\eqref{stima derivata ei}, \eqref{palla di sicurezza}, \eqref{palla Lip di sicurezza}.

\smallskip

19. --- 
The operator $\mS = I + w \pa_x\inv$, defined in \eqref{mS}, and its inverse $\mS\inv$ both satisfy the following estimates 
(where the $s$-decay norm $| \cdot |_s $ is defined in \eqref{matrix decay norm}):
\begin{eqnarray}
| \mS^{\pm1} - I |_s 
& \leq & \e C(s)(1 + \Vert u \Vert_{s + 2\tau_0 + 2{\mathfrak s}_0 + 9}),
\label{stima mS}
\\
| \mS^{\pm1} - I |_{s}^\Lipg 
& \leq & \e C(s)(1 + \Vert u \Vert_{s + 2\tau_0 + 2{\mathfrak s}_0 + 10}^\Lipg),
\label{stima Lip mS}
\\
\big| \partial_{u}\mS^{\pm1}(u) [h]\big|_{s} 
& \leq & \e C(s) \big(\Vert h \Vert_{s + 2\tau_0 + 2{\mathfrak s}_0 + 9} + \Vert u \Vert_{s + 2\tau_0 + 2{\mathfrak s}_0 + 10} \Vert h \Vert_{ 2\tau_0 + 3{\mathfrak s}_0 + 9} \big).
\label{stima derivata mS}
\end{eqnarray}
Thus \eqref{stima mS}-\eqref{stima derivata mS} for $\mS$ follow  by 
\eqref{stima w}-\eqref{stima derivata w} and the fact that 
the matrix decay norm $| \pa_x\inv  |_s \leq 1$, $s \geq 0$, 
using \eqref{multiplication}, \eqref{multiplication Lip}, \eqref{algebra}, \eqref{algebra Lip}.  
The operator $\mS^{-1}$ satisfies the same 
bounds \eqref{stima mS}-\eqref{stima Lip mS} 
by Lemma \ref{lem:inverti}, which may be applied thanks to
\eqref{stima mS}, \eqref{palla di sicurezza},  \eqref{palla Lip di sicurezza}
and $\e$  small enough. 

Finally \eqref{stima derivata mS} for $\mS\inv$ follows by
$$
\partial_u \mS^{-1}(u) [h] = - \mS^{-1}(u) \, \partial_u \mS(u)[h] \, \mS^{-1}(u) \, , 
$$
and \eqref{interpm}, \eqref{stima mS} for $\mS\inv$, 
and \eqref{stima derivata mS} for $\mS$.

\smallskip

20. --- The operatpr $\mR$, defined in \eqref{mL5} , 
where $r_0$, $r_{-1}$ are defined in \eqref{r0}, \eqref{r-1},
 satisfies the following estimates:
\begin{eqnarray}
\big| \mR \big|_s & \leq & \e C(s)(1 + \Vert u \Vert_{s + 2\tau_0 + 2{\mathfrak s}_0 + 12})\label{stima mR}\\
\big| \mR \big|_s^\Lipg& \leq & \e C(s)(1 + \Vert u \Vert_{s + 2\tau_0 + 2{\mathfrak s}_0 + 13}^\Lipg)\label{stima Lip mR}\\
\big| \partial_u \mR(u)[h] \big|_s & \leq & \e C(s)\big(\Vert h \Vert_{s + 2\tau_0 + 2{\mathfrak s}_0 + 12} + \Vert u \Vert_{s + 2\tau_0 + 2{\mathfrak s}_0 + 13} \Vert h \Vert_{2\tau_0 + 3{\mathfrak s}_0 +12} \big).
\label{stima derivata mR}
\end{eqnarray}
Let $T := r_0 + r_{-1}\partial_x^{-1}$.
By \eqref{multiplication}, \eqref{multiplication Lip}, \eqref{asymmetric tame product}, 
\eqref{stima w}, \eqref{stima Lip w}, 
\eqref{stima ei}, \eqref{stima Lip ei}, \eqref{stime mu1}, \eqref{stima mu3-1 Lip}, 
and using the trivial fact that 
$|\partial_x^{-1}|_s \leq 1$ and $|\pi_0|_s \leq 1$ for all $s \geq 0$,
we get 
\begin{eqnarray}
\big| T \big|_{s}  & \leq & \e C(s)(1 + \Vert u \Vert_{s + 2\tau_0 + 2{\mathfrak s}_0 + 12})\label{stima T}\\ 
\big|T \big|_{s}^\Lipg  & \leq & \e C(s)(1 + \Vert u \Vert_{s + 2\tau_0 + 2{\mathfrak s}_0 + 13}^\Lipg).
\label{stima Lip T}
\end{eqnarray}
Differentiating $T$ with respect to $u$, and using \eqref{multiplication},  \eqref{asymmetric tame product}, 
\eqref{stima derivata w}, \eqref{stima derivata ei}, \eqref{stime mu1}, 
\eqref{stima mu3-1 Lip} and \eqref{stima derivata mu 3}, 
one has
\begin{equation}\label{stima derivata T}
\big| \partial_u T(u)[h] \big|_s \leq \e C(s)\big(\Vert h \Vert_{s + 2\tau_0 + 2{\mathfrak s}_0 + 12} + \Vert u \Vert_{s + 2\tau_0 + 2{\mathfrak s}_0 +13} \Vert h \Vert_{ 2\tau_0 + 3{\mathfrak s}_0 +12} \big).
\end{equation}
Finally  \eqref{interpm}, \eqref{interpm Lip}
\eqref{stima mS}-\eqref{stima derivata mS}, 
\eqref{stima T}-\eqref{stima derivata T}
imply the estimates \eqref{stima mR}-\eqref{stima derivata mR}.

\smallskip

21. --- 
Using Lemma \ref{lemma astratto composizioni}, \eqref{palla di sicurezza} and all the previous estimates on 
${\cal A}, B, \rho, {\cal M} , {\cal T}, \mS$, 
the operators $\Phi_1 = {\cal A} B \rho {\cal M}  {\cal T} \mS$ and 
$\Phi_2 = {\cal A} B {\cal M}  {\cal T} \mS $, defined in \eqref{Phi 1 2 def}, satisfy \eqref{stima Phi 12 nel lemma}
(note that $ \s >  2\t_0 + 2{\mathfrak s}_0 + 9$).
Finally, if the condition \eqref{palla Lip di sicurezza} holds, we get  
the estimate \eqref{stima Lip Phi 12 nel lemma}. 

The other estimates \eqref{coefficienti costanti 1}-\eqref{stima R 3} 
follow by 
 \eqref{stima mu3-1 Lip}, \eqref{stima derivata mu 3}, \eqref{stime mu1}, \eqref{stima mR}-\eqref{stima derivata mR}. 
The proof of Lemma \ref{lemma:mostro} is complete.

\smallskip

\noindent 
\textbf{Proof of Lemma \ref{lemma:stime stabilita Phi 12}.} 
For each fixed $\ph \in \T^\nu$, ${\cal A}(\ph)h(x) := h(x+\b(\ph,x))$. 
Apply \eqref{tame-cambio-di-variabile} to the change of variable $\T \to \T$, $x \mapsto x + \b(\ph,x)$:
\[
\| {\cal A}(\ph) h \|_{H^s_x} 
\leq C(s) \big( \| h \|_{H^s_x} + | \b(\ph, \cdot) |_{W^{s, \infty}(\T)} \| h \|_{H^1_x}  \big).
\]
Since $| \b(\ph, \cdot) |_{W^{s, \infty}(\T)} \leq | \b |_{s, \infty}$ for all $\ph \in \T^\nu$, 
by \eqref{stima beta} we deduce \eqref{A(ph)}. 
Using \eqref{cambio di variabile meno identita}, \eqref{palla di sicurezza}, and \eqref{stima beta},
\[
\| ({\cal A}(\ph) - I) h \|_{H^s_x} 
\leq_s  | \b |_{L^\infty} \| h \|_{H^{s+1}_x}  
+ |\b|_{s,\infty} \| h \|_{H^2_x}  
\leq_s 
\e \big( \| h \|_{H^{s+1}_x} + \| u \|_{s + \mathfrak s_0 + 3} \| h \|_{H^2_x} \big). 
\]
By \eqref{stima beta tilde}, estimates \eqref{A(ph)} and \eqref{A(ph)-I}  
also hold for ${\cal A}(\ph)^{-1} = {\cal A}^{-1}(\ph)$ $: h(y) \mapsto h(y + \tilde \b(\ph,y))$. 

The multiplication operator $ {\cal M} (\ph) : H^s_x \to H^s_x$, $\, {\cal M} (\ph) h := v(\ph, \cdot) h$ satisfies 
\begin{multline}
\| ({\cal M} (\ph) - I)h \|_{H^s_x} 
= \| (v(\ph, \cdot) - 1) h \|_{H^s_x}
\leq_s
\| v(\ph, \cdot) - 1 \|_{H^s_x} \| h \|_{H^1_x} 
+ \| v(\ph, \cdot) - 1 \|_{H^1_x} \| h \|_{H^s_x} 
\\
\leq_s
\| v-1 \|_{s + \mathfrak s_0}  \| h \|_{H^1_x} 
+ \| v-1 \|_{1 + \mathfrak s_0} \| h \|_{H^s_x}
\leq_s 
\e \big( \| h \|_{H^s_x} + \| u \|_{s + \tau_0 + 2 {\mathfrak s}_0 + 6}  \| h \|_{H^1_x} \big)
\label{Phi(ph)-I}
\end{multline}
by \eqref{asymmetric tame product}, \eqref{multiplication}, Lemma \ref{Aphispace}, 
\eqref{stima v} and \eqref{palla di sicurezza}. 
The same estimate also holds for $ {\cal M} (\ph)^{-1} = {\cal M} ^{-1}(\ph)$, 
which is the multiplication operator by $v^{-1}(\ph,\cdot)$.
The operators $ {\cal T}^{\pm 1}(\ph)h(x) = h(x \pm p(\ph))$ satisfy
\be \label{mM(ph)-I}
\| {\cal T}^{\pm 1}(\ph) h \|_{H^s_x} = \| h \|_{H^s_x}, \quad
\| ({\cal T}^{\pm 1}(\ph) - I) h \|_{H^s_x} \leq \e \g_0^{-1} C \| h \|_{H^{s+1}_x}, \quad
\ee
by \eqref{cambio di variabile meno identita}, \eqref{palla di sicurezza}, \eqref{stima p} 
and by the fact that $p(\ph)$ is independent on the space variable.

By \eqref{interpolazione norme miste}, \eqref{stima mS}, \eqref{palla di sicurezza} 
and Lemma \ref{Aphispace}, 
the operator $\mS(\ph) = I + w(\ph,\cdot) \pa_x^{-1}$ and its inverse satisfy
\be \label{mS(ph)-I}
\| (\mS^{\pm 1}(\ph) - I)h \|_{H^s_x} 
\leq_s \e
\big( \| h \|_{H^s_x} + \| u \|_{s + 2\tau_0 + 3{\mathfrak s}_0 + 9} \| h \|_{H^1_x} \big).
\ee

Collecting estimates \eqref{Phi(ph)-I}, \eqref{mM(ph)-I}, \eqref{mS(ph)-I} we get
\eqref{Phi mM mS (ph)} and \eqref{Phi mM mS (ph) - I}.
Lemma \ref{lemma:stime stabilita Phi 12} is proved.

\noindent
Massimiliano Berti, Pietro Baldi, Dipartimento di Matematica e Applicazioni ``R. Caccioppoli",
Universit\`a degli Studi Napoli Federico II,  Via Cintia, Monte S. Angelo, 
I-80126, Napoli, Italy,  {\tt m.berti@unina.it},  {\tt pietro.baldi@unina.it}.
\\[1mm]
\noindent
Riccardo Montalto,  SISSA,  Via Bonomea 265, 34136, Trieste, Italy,  {\tt riccardo.montalto@sissa.it}.
\\[2mm]
\indent
This research was supported by the European Research Council under FP7  
and partially by the  
PRIN2009 grant ``Critical Point Theory and Perturbative Methods for Nonlinear Differential Equations".

\end{document}